\documentclass[a4paper]{article}

\usepackage{amssymb}
\usepackage{amsthm}
\usepackage{amsmath}
\usepackage{chemformula}
\usepackage{multirow, makecell}
\usepackage{thm-restate}
\usepackage{algorithm}
\usepackage[noend]{algpseudocode}
\usepackage[scientific-notation=true]{siunitx}
\usepackage{mathtools}
\usepackage{caption}
\usepackage{subcaption}
\usepackage{booktabs}
\DeclarePairedDelimiter\ceil{\lceil}{\rceil}
\usepackage{url}
\newcommand*{\vertbar}{\rule[0ex]{0.5pt}{2.5ex}}
\usepackage{hyperref}
\usepackage{geometry}
\newcommand{\ionp}{r}
\newcommand{\vect}{l}
\newcommand{\energy}{\Phi}

\newcommand{\rvect}{k}
\newcommand{\lattice}{{L_n}}
\newcommand{\rlattice}{{G_m}}
\newcommand{\allionp}{R}
\newcommand{\cell}{\mathcal{C}}
\newcommand{\erf}[1]{\ensuremath{\operatorname{erf}}\left(#1\right)}
\newcommand{\erfc}[1]{\ensuremath{\operatorname{erfc}}\left(#1\right)}

\newtheorem{assumption}{Assumption}
\newtheorem{definition}{Definition}
\newtheorem{prop}{Proposition}

\makeatletter
\newcommand{\pushright}[1]{\ifmeasuring@#1\else\omit\hfill$\displaystyle#1$\fi\ignorespaces}
\makeatother
\title{First Order Methods for Geometric Optimization of Crystal Structures}

\usepackage{authblk}
\author[1]{Antonia Tsili}
\author[2]{Matthew Dyer}
\author[3]{Vladimir Gusev}
\author[4]{Piotr Krysta}
\author[5]{Rahul Savani}

\affil[1]{Department of Computer Science, University of Liverpool, Liverpool, UK. {\tt a.tsili@liverpool.ac.uk}.}
\affil[2]{Department of Chemistry, University of Liverpool, Liverpool, UK. {\tt msd30@liverpool.ac.uk}.}
\affil[3]{Department of Computer Science, University of Liverpool, Liverpool, UK. {\tt vladimir.gusev@liverpool.ac.uk}.}
\affil[4]{Department of Computer Science, University of Liverpool, Liverpool, UK. {\tt pkrysta@liverpool.ac.uk}.}
\affil[5]{Department of Computer Science, University of Liverpool, Liverpool, UK. {\tt rahul.savani@liverpool.ac.uk}.}
\date{}

\begin{document}
\maketitle
\thispagestyle{empty}


\begin{abstract}
The geometric optimization of crystal structures is a procedure widely used in Chemistry that changes the geometrical placement of the particles inside a structure. It is called structural relaxation and constitutes a local minimization problem with a non-convex objective function whose domain complexity increases according to the number of particles involved. In this work we study the performance of the two most popular first order optimization methods in structural relaxation. Although frequently employed, there is a lack of their study in this context from an algorithmic point of view. We run each algorithm in combination with a constant step size, which provides a benchmark for the methods' analysis and direct comparison. We also design dynamic step size rules and study how these improve the two algorithms' performance. Our results show that there is a trade-off between convergence rate and the possibility of an experiment to succeed, hence we construct a function to assign utility to each method based on our respective preference. The function is built according to a recently introduced model of preference indication concerning algorithms with deadline and their run time. Finally, building on all our insights from the experimental results, we provide algorithmic recipes that best correspond to each of the presented preferences and select one recipe as the optimal for equally weighted preferences. 

Alongside our results we present our Python software \href{https://github.com/lrcfmd/veltiCRYS.git}{``\textit{veltiCRYS}''}, which was used to perform the geometric optimization experiments. Our implementation, can be easily edited to accommodate other energy functions and is especially targeted for testing different methods in structural relaxation.
\end{abstract}

\pagenumbering{arabic}
\section{Introduction}
\label{intro}

The study of crystal structures is a major component of materials discovery. Their properties are the subject of examination in many applications, ranging from pharmaceuticals to signal processing and beyond. Crystal structures are periodic formations, meaning they can be represented as tilings expanding towards all 3 dimensions of Euclidean space. Since each tile is identical to the rest, it suffices to define one ``central'' tile to which changes are applied, so as to amend the whole structure. Each such tile is called a unit cell and can be represented by a parallelepiped. The unit cell comprises an arrangement of ions which determines the properties of the crystal and, given the ions' number and element type, we seek to find their optimal geometrical placement in the $\mathbb{R}^3$ space spanned by the unit cell. This is a hard optimization problem. In fact, it is a local minimization problem of a function $\energy$, called energy potential, with $3N+9$ variables, when $N$ is the number of the ions in the unit cell and $\energy$ is the objective function. We search for the function's approximate local minimum, which corresponds to bringing the crystal to an energy equilibrium. In practice this is achieved by procedures such as heat application, hence, this geometric optimization is also called structural relaxation in a Chemistry context and constitutes a frequently employed procedure with many applications in Computational Chemistry. For example, it is a particularly important part of Crystal Structure Prediction, for which it can take up to $90\%$ of the computation time in experiments. The problem of geometric optimization of a crystal structure can be defined as
\begin{equation}
 \begin{gathered}
    \tag{P}
    \min{\energy(x)} \\
    x = (\ionp_1, \ionp_2, ..., \ionp_N, \vect_1, \vect_2, \vect_3),\ \ionp_i, \vect_m \in \mathbb{R}^3, i\in [N], m\in \{1,2,3\} \\ 
    \text{ s.t. } r_{i,j} > 0 \text{ where } r_{i,j}=\|\ionp_i - \ionp_j\|, \ j\in [N]
    \label{eqn:Problem}
\end{gathered}   
\end{equation}
in which $\ionp_i$ denotes the position of ion $i$, $N$ is the fixed number of the ions in a unit cell and $\vect_m$ denotes a lattice vector. $x$ can be more concisely written as $x=(R,L)$, where $R$ is the matrix comprised of the ion positions and $L$ is the matrix comprised of the lattice vectors. We call $x$ an approximate local minimiser of $\energy$ when $g(x)<\epsilon$ for some small $\epsilon>0$, in which $g$ is the average component value of the gradient norm $\|\nabla\energy\|$.

$\energy$ is a nonlinear, non-convex function with a complicated domain called Potential Energy Surface (PES). The complexity of PES increases along with the number of ions $N$ included in the calculation, as the number of local minima also increases. $\energy$ is locally $C^2$-smooth for $r_{i,j}>0$\cite{lindbo2011spectral}. More specifically, the function is not continuous in areas where ions $i,j$ are separated by a pairwise distance approaching zero. We study the application of two standard unconstrained continuous optimization algorithms, Gradient Descent and Conjugate Gradient, in finding an approximate the local minimum of the energy potential function $\energy$ through structural relaxation. Although widely employed, they are not systematically optimised in this context and parameter configurations such as the choice of step size are kept on default. We assume that the minimization path remains in a feasible area where the pairwise distance is a positive number. However, there can be cases for which the discontinuity is approached -- this is called Buckingham catastrophe. Our work shows that the optimization process can be correspondingly adapted to reach the minimiser of $\energy$ without falling into such a case. However, there is a trade-off between a method's robustness and convergence speed.

\begin{figure}
    \centering
    \includegraphics[width=0.65\textwidth]{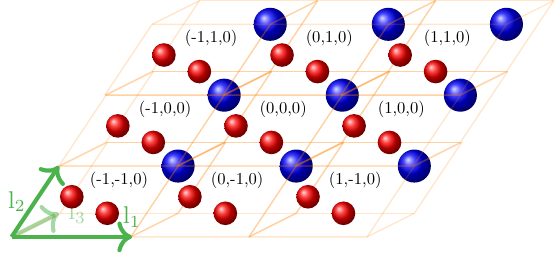}
    \caption{\textbf{Representation of a 3D crystal lattice.}A mock representation of a crystal comprising ions of two kinds of elements (red and blue) in the 3-dimensional Euclidean space. The green arrows represent the lattice vectors and all ions' positions can be represented as their linear combination.}
     \label{fig:crystal_vects}
\end{figure}

The calculation of the energy of crystal structures is a long studied subject that has been investigated since the first attempts to understand materials' properties.  Various computer programs, with GULP~\cite{Gale2003TheGULP} and LAMMPS~\cite{LAMMPS} among the most popular, have been developed for Chemistry applications such as the study of Molecular Dynamics, which offer this functionality. However, we found that they do not accommodate the necessary freedom to amend, embed and test various optimization algorithms along with the parameters that affect their efficiency. 

\subparagraph*{Our contributions.}
In this paper we investigate the derivation of the energy function $\energy$ and the forces $-\nabla\energy$ and present their most numerically robust form. The energy $\energy$ is a potentially infinite summation, therefore it is crucial in which way it is calculated and in which order the terms are summed up. We implement the energy and forces and systematically  design a set of experiments to test two popular first order local optimization algorithms frequently used by Chemists, Gradient Descent and Conjugate Gradient, using our implementation. More specifically:
\begin{itemize}
    \item  We carefully construct the energy potential model and present the analytical process of evaluating the energy function. Towards this goal, we prove some useful propositions, critical for the Computational Chemistry related background theory. They can be found in  Appendix~\ref{app1}. We propose a new geometric method called Inflated Cell Truncation which is also a part of the energy evaluation.
    \item We provide a thorough derivation of the forces, the energy function's first derivatives, with formal proofs in Appendix~\ref{app3} and we explain what are the parameters with respect to which they are evaluated. Because of treating a parameter of the energy function as a constant, some terms of the derivatives are usually overlooked in literature. Our derivation, specifically targeted for structural relaxation, includes these terms and highlights this difference. 
    \item All derived formulae and the two algorithms have been implemented with our software \textit{veltiCRYS} found at \url{https://github.com/lrcfmd/veltiCRYS.git}.
    \item We conduct extensive experiments to formally compare the two algorithms in the setting of structural relaxation. The setting of the experiments is presented in Section~\ref{eexperiments}.
    \item To the best of our knowledge, this is the first paper studying in detail the algorithms applied to the problem in question. As such, we provide a performance benchmark for first order methods using each algorithm with constant step in Section~\ref{results}.
    \item We improve the algorithms' benchmarked performance by designing scheduling rules for the step size in Section~\ref{results}.
    \item  We provide a thorough analysis of the experiments and a tool to formally evaluate the algorithms' suitability to our preferences in Section~\ref{results}. This tool is a utility function that can help to decide which of the algorithms will be more useful in certain applications.
\end{itemize}

The road-map of the paper is as follows. In Section~\ref{energy} we introduce the necessary background that formulates the problem we discuss, including the results from all our derivations. In Section~\ref{eexperiments} we present the technical setting and describe our experiments. In Section~\ref{results} we present the results of the experiments, which are separated into two parts; results from optimization with constant step size and results from optimization with adaptive step size. The results are accompanied by our analysis, which also contains usage of the tool that we propose for formally determining the utility of crystal relaxation methods according to the user's preferences. Finally, we elaborate on the utility scores of each method used in our experiments and state our conclusions.

\section{The energy model}
\label{energy}

\subsection{Background and notation} \label{notation}
Gradient Descent and Conjugate Gradient algorithms are no strangers to geometric optimization in Chemistry applications. Any generic geometric optimization review will give credit to these two simple algorithms in unconstrained minimization settings, as the work of Schlegel shows~\cite{Schlegel2011GeometryOptimization}. Publications as early as Catlow and Mackrodt's~\cite{Catlow1982} study function minimisation through ion displacement and lattice deformation. The work of Payne et al.~\cite{Payne1992IterativeMinTecniques} describes the use of Conjugate Gradient for energy minimisation of the crystal structure but also in energy functional minimisation. Many approaches to geometric optimization have been proposed since then~\cite{Bitzek2006StructuralSimple} with machine learning getting increasingly more attention~\cite{DelRio2019LocalStructuresb,Born2021GeometryApproaches}, as in any application related to Computer Science. However, there is a lack of analysis of the algorithmic aspects of the aforementioned simple first order algorithms. These are still being used today for structural relaxation~\cite{Fan2023EffectNa2Bi2SeO33F2},\cite{Olaniyan2023DetectionAdsorption}, but little analysis has been provided in this context. The recent study of Salih and Faraj~\cite{Taha2021ComparisonMatlab} investigates the algorithms' performance and compares them on the basis of 3 simple nonlinear function applications. Our work extends this comparison to a much more complex setup and is focused on the efficiency of their direction selection by excluding the line search.

For the rest of the article we refer only to optimization problems whose goal is function minimisation. The specific objective function will be introduced in detail. Before continuing further, we will list some useful definitions and notations that will be frequently used. In general, given a vector $u\in\mathbb{R}^3$, we denote as $u_\lambda$ one of its 3 components in the Euclidean space so that $\lambda \in \{x,y,z\}$, where $\lambda$ is any Greek letter.

\begin{figure}
    \centering
    \includegraphics[width=0.4\textwidth]{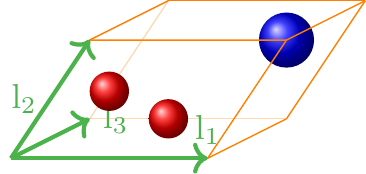}
    \caption{\textbf{The unit cell of the crystal structure depicted in Figure~\ref{fig:crystal_vects}.} The vectors $\vect_1$,$\vect_2$ and $\vect_3$ represent the 3 lattice vectors.}
    \label{fig:unitcell_vects}
\end{figure}

The periodicity of the crystal structure allows for amending a single tile of its pattern in order to simulate changes to its entire formation. A tile is called unit cell $\cell_n$ and is a parallelepiped built on a set of three vectors $L = (\vect_1, \vect_2, \vect_3)$ such that $\vect_t \in \mathbb{R}^3 ,\ \forall t \in \{1,2,3\}$ and contains N positions $\ionp_1,\ionp_2,...,\ionp_N$, $\ionp_i \in \mathbb{R}^3, \ i \in [N]$, where the ions are placed. Formally, it is defined as follows:
\begin{definition}[unit cell]
The  smallest  fundamental  arrangement  of  the  ions  positions’  that  reflects  the  crystal’s  symmetry and structure is called unit cell.
\end{definition} 
\noindent An example is given in Figure~\ref{fig:unitcell_vects}, which depicts the unit cell of the example structure in Figure~\ref{fig:crystal_vects}.  By defining a ``central'' unit cell $\cell_{(0,0,0)}$ for reference, any unit cell copy in the crystal can be represented by using translations of the points in  $\cell_{(0,0,0)}$.
\begin{definition}[lattice]\label{def:lattice}
Given the translation vectors $\lattice=n^TL$ with $L = (\vect_1, \vect_2, \vect_3)$, $n\in\mathbb{N}$, lattice is the set $\mathcal{D}_\lattice \subset \mathbb{R}^3 $ of mathematical points that correspond to the infinitely repeated positions of ions which form the crystal structure.
\end{definition}
\noindent  The points of the lattice can be defined using the lattice vectors $L$, such that $\forall \ionp_{i_1},\ionp_{i_2} \in \mathcal{D}_\lattice$ we have $\ionp_{i_2} = \ionp_{i_1} + n_1\vect_1 + n_2\vect_2 + n_3\vect_3$ and $n_1,n_2,n_3$ are arbitrary integers. Another important set of points that we will introduce is called the reciprocal lattice.
\begin{definition}[reciprocal lattice]
The Fourier transform of the lattice of Definition~\ref{def:lattice}. It is formed by the reciprocal lattice vectors  $G=(\rvect_1, \rvect_2, \rvect_3)$ with $k_t = 2\pi \cdot (\vect_{t\text{mod}3+1}\times \vect_{(t+1)\text{mod}3+1})/\vect_t^T(\vect_{(t\text{mod}3+1}\times \vect_{(t+1)\text{mod }3+1})$, $t\in {1,2,3}$.
\end{definition}
\noindent The reciprocal lattice is used along with the real space lattice (Definition~\ref{def:lattice}) in crystallography for periodic structures. The points can be defined in the same way as for the real space lattice using the vectors $\rlattice = m^TG$, such that $\forall \ionp_{i_1},\ionp_{i_2} \in \mathcal{D}_\rlattice$ we have $\ionp_{i_2} = \ionp_{i_1} + m_1\rvect_1 + m_2\rvect_2 + m_3\rvect_3$ and $m_1,m_2,m_3$ are arbitrary integers. Finally, our problem's objective function, henceforth $\energy : \mathbb{R}^{3N+9}\rightarrow\mathbb{R}$, is the energy that exists on account of the ions' positions, the potential energy. More formally:
\begin{definition}[potential energy]
Potential is the energy stored in a structure as a result of the relative positions $\ionp_i\in\mathbb{R}^3, \ i\in [N]$ of the ions and the forces that one exerts to another.
\end{definition}

Since our research is based on ionic structures, we model our problem using the Buckingham-Coulomb energy potential. More specifically, we use structures that include ions of Sr, Ti and O randomly placed in the unit cell, so that the expected minimum of the function is represented by the structure \ch{Sr_3Ti_3O_9}, a commonly used crystal for benchmarking in Crystal Structure Prediction. The total energy comprises two summation schemes built by repeating the unit cell $\cell_{(0,0,0)}$ in all dimensions of the Euclidean space and by evaluating the distance of all pairs of ions. The Coulomb energy potential addresses the electrostatic energy and the Buckingham energy potential accounts for the Pauli repulsion and van der Waals interactions between two atoms. These interactions depend the interatomic distance between the atoms, thus, the resulting energy is basically a function of the separation vector between two ions. In the case of the Buckingham potential, its two terms per summand represent repulsion and attraction respectively, while the Buckingham constants $A,C$ and $\rho$ have been experimentally determined in literature~\cite{Collins2018ThePrediction} and differ among the combination of elements in the atoms pairs. The respective energy functions are given in Equations~(\ref{eq:ECoul}),(\ref{eq:EBuck}) and they are combined as shown in Equation~(\ref{eq:energyall}).
\begin{align}
\energy_{Coul}(\allionp, L) = & \ \frac{k_e}{2}\sum_{i,j}^{N'}\sum_{n\in\mathcal{N}}\frac{q_iq_j}{\|\ionp_{i,j,n}\|} \label{eq:ECoul} \\
\energy_{Buck}(\allionp, L) = & \ \frac{1}{2}\sum_{i,j}^{N'}\sum_{n\in\mathcal{N}} A\exp(-\frac{\|\ionp_{i,j,n}\|}{\rho}) - \frac{C}{\|\ionp_{i,j,n}\|^6} \label{eq:EBuck} 
\end{align}
\begin{gather}
    \label{eq:energyall}
    \energy = \ \energy_{Coul} + \energy_{Buck}
\end{gather}
The independent variables of both functions $\energy_{Coul}:\mathbb{R}^{3N+9}\rightarrow\mathbb{R}$, $\energy_{Buck}:\mathbb{R}^{3N+9}\rightarrow\mathbb{R}$ are the ion positions $\allionp = (\ionp_1,\ionp_2,...,\ionp_N)$ and the lattice vectors $L=(\vect_1,\vect_2,\vect_3)$. Our energy calculations involve particles subjected to periodic boundary conditions (PBCs), so that each ion $i$ of the central reference unit cell $\cell_{(0,0,0)}$ interacts with every other ion $j$ of the same unit cell and the neighbouring ions residing in surrounding images of $\cell_{(0,0,0)}$. Each ion $i$ is separated from images of ion $j$ by the pairwise distances $\|\ionp_{i,j,n}\| = \|\ionp_i + \lattice - \ionp_j\|$, where the vector $\lattice$ creates the periodic translations of each ion $j$ of the central unit cell. More specifically, $n$ is a tuple, so that $n = (n_1, n_2, n_3), \ n_1, n_2, n_3\in \mathbb{Z}\cup\{0\}$. It contains three integers with each tuple component standing for a dimension in the 3-dimensional Euclidean space and the tuple itself defines a relative position to $\cell_{(0,0,0)}$, as shown in Figure~\ref{fig:crystal_vects}. In other words, the tuple $n = (1, 0, 0)$ corresponds to the image $\cell_{(1,0,0)}$ of $\cell_{(0,0,0)}$, which is its exact adjacent copy in the direction parallel to the x-axis. Since the crystal structure is a formation theoretically expanding infinitely in space, then each component $n_t\in(-\infty,\infty)$. However, as $n_t\rightarrow-\infty$ the energy contributions tend to zero, because $1/\|\ionp_{i,j,n}\|\rightarrow 0$. Thus, we can define the finite set $\mathcal{N}$ with $n\in\mathcal{N}$ using symmetric lower and upper bounds for the tuple's components, in order to include the important energy contributions in the summation. For the rest of the article, a summation over $n$ is a summation with all $n\in \mathcal{N}$. There are different ways for defining these bounds in literature, but we propose a simple geometric way called Inflated Cell Truncation, which is presented in Section~\ref{ICT}. 

The convergence of the sum over all pairwise interactions is strongly affected by the selection of interacting images included in the set of $n$ triplets $\mathcal{N}$, but also the order of the summation in Equations~(\ref{eq:ECoul}),(\ref{eq:EBuck})~\cite{1980SimulationConditions}. In this form, $\energy_{Coul}$ is conditionally convergent; the order of the terms in the summation determines whether it will finally converge. Moreover, the need for finite interaction terms imposes the use of methods which abruptly terminate the summation up until some designated distance away from the ions of the central unit cell. This term exclusion causes precision loss problems, discontinuity of limits, such as the derivatives, and others ~\cite{Feller1996EffectWater,Hunenberger1998AlternativeWater,Steinbach1994NewSimulation}. For this reason, we instead expand $\energy_{Coul}$ using the Ewald summation.

\subsection{Long Range Term}
\label{longrangeterm}

As already discussed, the Coulomb potential in its original form is evaluated directly in real space and is conditionally convergent, the conditions being dependent on the order of the summation. When employing the Ewald method, our target is to split the summation into short and long range contributions, so as to treat each part differently and arrive to two rapidly and absolutely convergent summation parts. The part of the summation that is responsible for the long range interactions will be evaluated in reciprocal space, hence the term reciprocal part of the summation. Many reports have presented different versions of the Coulomb energy potential expanded with Ewald summation~\cite{Pickard2018RealBackground,Wang2019EwaldDipoles,Stamm2018ASelf-terms,Lee2009EwaldSupercell,Toukmaji1996EwaldSurvey,Wang2019EwaldDipoles}. There are also some reports presenting the dispersion energy of interatomic potentials like Buckingham~\cite{InTVeld2007ApplicationForces}. For the rest of this paper, let $\cell_n : (\lattice, \allionp)$ denote a unit cell image described by real lattice vectors $\lattice$ and ion positions $\allionp$. Furthermore, let $\cell$ be the set of all such images and $\langle \cell, \mathcal{N}, \mathcal{M} \rangle$ denote a crystal structure whose energy is calculated using the real vectors $\lattice, n \in \mathcal{N}$ and the reciprocal vectors $\rlattice, m \in \mathcal{M}$.

\begin{restatable}{prop}{Cshort}
\label{prop:real}
Let $\energy_{Coul}:\mathbb{R}^{3N+9}\rightarrow\mathbb{R}$ be the Coulomb potential function describing the energy of crystal structure $\langle \cell, \mathcal{N}, \mathcal{M} \rangle$. By expansion with the Ewald summation, the short range interactions $\energy^S_{Coul}$ are calculated as in the following
\begin{equation}
    \energy^S_{Coul}(\allionp, L) = \frac{k_e}{2}\sum_n\sum_{i=1}^N\sum_{j=1}^{N^{'}}\frac{q_iq_j}{\|\ionp_{i,j,n}\|}\erfc{\alpha\|\ionp_{i,j,n}\|}.
\end{equation}
\end{restatable}
The proof of Proposition~\ref{prop:real} can be found in Appendix~\ref{app1}. The $\alpha$ parameter controls the balance between real and reciprocal space terms, meaning that it defines the boundary that separates short and long range interactions. Its value can differ depending on the empirical experimental results. For our case, we use Catlow's definition~\cite{Jackson1988ComputerStructure} so that
\begin{equation}
\label{eq:alpha}
    \alpha = \frac{N^{1/6}\sqrt{\pi}}{V^{1/3}}.
\end{equation}

\begin{restatable}{prop}{Clong}
\label{prop:longrange}
Let $\energy_{Coul}$ and $\langle \cell, \mathcal{N}, \mathcal{M} \rangle$ be defined as in Proposition~\ref{prop:real}. Then, the long range interactions $\energy^L_{Coul}$ of the Coulomb potential summation are calculated as
    \begin{equation}
    \energy_{Coul}^L(\allionp,L) = \frac{k_e}{2}\sum_n\sum_{i=1}^N\sum_{j=1}^N \frac{q_iq_j}{\|\ionp_{i,j,n}\|}\erf{\frac{\|\ionp_{i,j,n}\|}{\sqrt{2}\sigma}}.
\end{equation}
\end{restatable}

Proof for Proposition~\ref{prop:longrange} is in Appendix~\ref{app1}. The sum of $\energy^S_{Coul}+\energy^L_{Coul}$ gives the total value of electrostatic energy of the crystal structure. By using a Gaussian distribution for modelling each charge concentration, $\energy^L_{Coul}$ achieves fast and absolute convergence~\cite{Pickard2018RealBackground,Wang2019EwaldDipoles,Stamm2018ASelf-terms,Lee2009EwaldSupercell,Toukmaji1996EwaldSurvey,Wang2019EwaldDipoles}.

\begin{restatable}{prop}{Crecipexp}
\label{prop:reciprocalexp}
Let $\energy_{Coul}$ and $\langle \cell, \mathcal{N}, \mathcal{M} \rangle$ be defined as in Proposition~\ref{prop:real}. Then, the long range interactions $\energy^L_{Coul}$ can be calculated in reciprocal space as follows
    \begin{equation}
        \energy_{Coul}^L(\allionp,L) = \sum_{m}\sum_{i=1}^N\sum_{j=1}^{N}\frac{4\pi}{V\|\rlattice\|^2}\exp{\left(-\frac{\|\rlattice\|^2}{4\alpha^2}\right)}\exp{(i\rlattice\ionp_{i,j})}.
    \end{equation}
\end{restatable}
This form of $\energy_{Coul}^{L}$ achieves fast and absolute convergence. The subscript $m$ of vector $\rlattice$ stands for a tuple of three integers, so that  and $m = (m_1, m_2, m_3)$, $m_1, m_2, m_3\in \mathbb{Z}$, equivalently to real space $\lattice$. The tuple $m$ corresponds to the translated image $\cell_{m}$ of the central unit cell $\cell_{(0,0,0)}$ in reciprocal space and $\ionp_{i,j} = \ionp_i - \ionp_j$ is the separation vector of two ions in the same unit cell. For the rest of the document, all summations over $m$ include all tuples $m\in \mathcal{M}$ except for $m = (0,0,0)$. $\mathcal{M}$ is populated with our Inflated Cell Truncation method in the same manner that $\mathcal{N}$ is populated, but in reciprocal space. The proof of Proposition~\ref{prop:reciprocalexp} can be found in Appendix~\ref{app1}.

We also employ techniques, such as a Fourier series expansion and symmetry conventions, to dispose of imaginary terms. The result can be seen in Proposition~\ref{prop:reciprocal}.
\begin{restatable}{prop}{Crecip}
\label{prop:reciprocal}
Let $\energy_{Coul}$ and $\langle \cell, \mathcal{N}, \mathcal{M} \rangle$ be defined as in Proposition~\ref{prop:real}. Then the fast convergent $\energy^L_{Coul}$ is given by the function
\begin{equation}
    \label{eq:potentialLong}
    \energy^L_{Coul}(\allionp, L) = \frac{2\pi k_e}{V\|\rlattice\|^2}\sum_{m}\sum_{i=1}^N\sum_{j=1}^{N} q_iq_j\exp{\left(-\frac{\|\rlattice\|^2}{4\alpha^2}\right)}\cos{(\rlattice\ionp_{i,j})}.
\end{equation}
\end{restatable}
The proof of Proposition~\ref{prop:reciprocal} can be found in Appendix~\ref{app1}.

\begin{restatable}{prop}{Coulall}
\label{prop:ECoul}
Let $\energy_{Coul}:\mathbb{R}^{3N+9}\rightarrow\mathbb{R}$ be the Coulomb potential function describing the energy of crystal structure $\langle \cell, \mathcal{N}, \mathcal{M} \rangle$. The fast convergent form of $\energy_{Coul}$ is formed by the following equations:
\begin{equation}
\label{eq:ECoulEwald} 
\begin{aligned}
& \energy^S_{Coul}(\allionp, L) & = & \quad k_e\sum_{i,j}^{N'}\sum_n q_iq_j \frac{erfc(\alpha \|\ionp_{i,j,n}\|)}{2\|\ionp_{i,j,n}\|} \\
& \energy^L_{Coul}(\allionp, L) & = & \quad k_e\sum_{i,j}^N \sum_{m} q_iq_j\frac{2\pi}{V\|G_m\|^2}\cdot\exp{\left(-\frac{\|G_m\|^2}{4\alpha ^2}\right)}\cdot\cos{(\rlattice\cdot \ionp_{i,j})} \\
& \energy^{self}_{Coul}(\allionp, L) & = & -k_e\sum_{i=1}^N q_i^2\frac{\alpha}{\sqrt{\pi}} \\
\end{aligned}
\end{equation}
so that 
\begin{gather}
    \label{eq:ECoulall} 
    \energy_{Coul} = \energy^S_{Coul} + \energy^L_{Coul} + \energy^{self}_{Coul}.
\end{gather}
\end{restatable}
 
The proof of Proposition~\ref{prop:ECoul} can be found in Appendix~\ref{app1}. It is important to notice that for the long range term $\energy^L_{Coul}$ the summation includes the pair $i=j$ for $n=0$ and a new set of triplets $m\in\mathcal{M}$ that correspond to the lattice vectors of unit cell images in reciprocal space, so that $G=(\rvect_1,\rvect_2,\rvect_3)$ are the reciprocal vectors and the sum is over $\rlattice$. The summation cost is owed to the number $N$ of ions, the set of real $\mathcal{N}$ and reciprocal $\mathcal{M}$ vectors, as well as the $\alpha = 1\sqrt{2}\sigma$ parameter, and reaches a complexity of $O(N^\frac{3}{2})$~\cite{Lee2009EwaldSupercell}. 

Our empirical results showed that, when the number of summands of the Buckingham potential is small, the abrupt exclusion of Buckingham terms can cause incontinuities of the gradient and steep slopes on the PES. For this reason, we revisit the problematic dispersion term $\|\ionp_{i,j,n}\|^{-6}$, which is the prevailing term of a series with larger powers $\|\ionp_{i,j,n}\|^{-7}, \|\ionp_{i,j,n}\|^{-8}$ and so on~\cite{Buckingham1965TheoryForces}. The dispersion term is expanded using the Ewald method in the same fashion as before, so that the final form of the Buckingham potential that we use is the presented in Proposition~\ref{prop:EBuck}.

\begin{prop}
\label{prop:EBuck}
Let $\energy_{Buck}:\mathbb{R}^{3N+9}\rightarrow\mathbb{R}$ be the Buckingham potential function describing the energy of crystal structure $\langle \cell, \mathcal{N}, \mathcal{M} \rangle$. If the
dispersion term is expanded using Ewald summation, $\energy_{Buck}$ is formed by  $\energy_{Buck}^S$,  $\energy_{Buck}^L$ and  $\energy_{Buck}^{self}$ in the following way:
\begin{equation}\label{eq:EBuckEwald} 
\begin{aligned}
& \energy^S_{Buck}(\allionp, L) & = & \quad \frac{1}{2}\quad \sum_{i,j}^{N'}\sum_{n} \left[ A_{ij}\exp{(-\frac{\|\ionp_{i,j,n}\|}{\rho})} - \right. \\ 
& & & \left.\frac{C_{ij}}{\|\ionp_{i,j,n}\|^6}\left( 1+\alpha^2\|\ionp_{i,j,n}\|^2+\frac{\alpha^4\|\ionp_{i,j,n}\|^4}{2} \right)\exp{(-\alpha^2\|\ionp_{i,j,n}\|^2)}\right] \\
& \energy^L_{Buck}(\allionp, L) & = & - \frac{1}{2}\sum_{i,j}^N C_{ij}\frac{\pi^{3/2}}{12V}\sum_{m}\left[\sqrt{\pi}\cdot\erfc{\frac{\|\rlattice\|}{2\alpha}} + \right. \\ 
& & & \left. \qquad \left( \frac{4\alpha^3}{\|\rlattice\|^3}-\frac{2\alpha}{\|\rlattice\|}\right)\exp{\left(-\frac{\|\rlattice\|^2}{4\alpha^2}\right)}\right]\cdot\cos{(\rlattice\ionp_{i,j})}\|\rlattice\|^3\\ 
& \energy^{self}_{Buck}(\allionp, L) & = & \quad \frac{1}{2}\sum_{i,j}^N\left(-\frac{C_{ij}}{3V}\pi^{3/2}\alpha^3 \right)+\sum_{i=1}^N \frac{C_{ii}\alpha^6}{12}\\
\end{aligned}
\end{equation}
so that
\begin{gather}
    \label{eq:EBuckall} 
    \energy_{Buck} = \energy^S_{Buck} + \energy^L_{Buck} + \energy^{self}_{Buck}.
\end{gather}
\end{prop}
The proof follows the same pattern as the proof of Proposition~\ref{prop:ECoul}.

\subsection{Inflated Cell Truncation}
\label{ICT}

There are theoretical and practical issues that dictate the exclusion of energy terms which stand for far too long range interactions. On the one hand, there is a balance to maintain between short range and distant energy contributions, as the interaction terms' number increases along with the distance from a certain ion. This is mitigated owing to the slow decay of the functions that emerge from the Ewald expansion. On the other hand, implementation issues call for finite sums, hence, one of the decisions to be made regarding the energy model concerns the possible values of each of the integers $n_1, n_2, n_3, m_1, m_2, m_3$. These numbers declare which neighbouring copies of the unit cell are to be taken into account in the energy and derivatives evaluation. In other words, they define how many terms are to be added into the summation of the energy and forces, thus, affect greatly the convergence of these numbers. The starting point to construct the set of these numbers is the cutoff~\cite{Jackson1988ComputerStructure} $r_{off}$, which reflects the maximum distance $\|\ionp_{i,j,n}\|$ in Angstroms that should separate two ions $i,j$. We set the cutoff for real $r^S_{off}$ and reciprocal space $r^L_{off}$ according Catlow's formulae~\cite{Jackson1988ComputerStructure}, which aim to minimise the number of summands to achieve an accuracy $A$:
\begin{equation} \label{eq:cut}
    r^S_{off} = \frac{\sqrt{-\log{(A)}}}{\alpha}, \quad r^L_{off} = 2\alpha\sqrt{-\log{(A)}},
\end{equation}
where $V$ is the volume of one unit cell. Another thing that affects the result of the summation is the symmetry of the terms of interaction that would be included. More specifically, it would not be realistic to select neighbouring cells in a way that ultimately approaches the shape of a tube.

Traditionally, the simplest method used in order to keep the energy summation finite is the truncation method~\cite{Linsew1986TruncationLiquids}. According to this method, all summands are multiplied by a function $\psi(\|\ionp_{i,j,n}\|)$
\begin{equation}
    \psi(\|\ionp_{i,j,n}\|) = \left\{
\begin{array}{ll}
      1 ,& \|\ionp_{i,j,n}\|<r_{off} \\
      0 ,& \|\ionp_{i,j,n}\|\geq r_{off}
\end{array} 
\right. 
\end{equation}
that allows only ions $j$ within the range of some ion $i$ to be studied. Following this scheme, various ways of creating neighbour lists that store the neighbours within range
for each atom in the central unit cell have been developed. One of the most widely used algorithms include the Verlet Neighbour list~\cite{Frenkel1996UnderstandingApplications}, which maintains the array of the neighbours of the atom within a fixed cutoff, an d the Linked-Cell method~\cite{Brooks1989ComputerLiquids}, which creates a list for each tile in a divided supercell. However, the first method requires a potentially large memory allocation and both methods suffer from a time-consuming bookkeeping, especially for our constant changing setting of the lattice. Other more modern approaches include that of Mason~\cite{Mason2005FasterRepresentation}, who focuses on the storing power of a bitmap, so that the neighbours are represented in memory in a way that enables locality information storage, and that of Zhang et al.~\cite{Zhang2018ASimulation}, whose work revolves around the consideration of potential neighbours near the margin of the ion's sphere of interaction. These techniques, however, are designed to facilitate Molecular Dynamics processes, whereas in our case we can fully exploit the symmetry of the ionic crystal.
 
For this work, we provided a geometric solution that makes use of the properties of the unit cell called InflatedCellTruncation. Let $S$ be a sphere $(O,r)$ where $O$ is the centre of gravity of the central unit cell $\cell_{(0,0,0)}$ and $r$ is a radius with length equal to the cutoff value $r=r_{off}$. 
\begin{assumption}
We assume that 
\begin{equation}
    r_{off}\geq max\{\|\vect_1\|,\|\vect_2\|,\|\vect_3\|\}
\end{equation}
\end{assumption}
meaning that there is at least one whole unit cell in $S$. We, then, identify each plane $P$ to which a face of the unit cell parallelepiped belongs and we assume the corresponding translations $P'$ of the planes so that $P'$ is parallel to $P$ and tangent to $S$. Because of symmetry, we only need to move the 3 adjacent planes defined by $\vect_1, \vect_2, \vect_3$ and apply the opposite movements to the rest. Then, we compute the length of the translation vector $t$ which performs $P'=P+t$ and enumerate all images of unit cells that are encased in the plane translations $P'$.

\begin{restatable}{thm}{translation}\label{th:translation}
    Let $\cell_n$ be the parallelepiped of a unit cell and $A$ a face of $\cell_n$ with $\upsilon$ the height that corresponds to $A$. The vector $t$ that translates the plane $P$ of $\cell_n$ to a parallel plane $P'$ tangent to $(O,r_{off})$ has length $\|t\|=r_{off}-\frac{\upsilon}{2}$ and is parallel to the normal $N_{P'}$ of plane $P'$.
\end{restatable}
The proof of Theorem~\ref{th:translation} can be found in Appendix~\ref{app2}.

\begin{restatable}{corollary}{imgs}
\label{cor:imgs}
    The distance between $P$ and $P'$ can fit $\frac{\|t\|}{\|\upsilon\|} + 1/2$ many unit cell images.
\end{restatable}

The proof of Corollary~\ref{cor:imgs} can be found in Appendix~\ref{app2}.The previous results reveal that, in order to include a reasonable number of neighbours for the ions of the central unit cell $\cell_{(0,0,0)}$, we can find the translation vector $t$ for each face and separate it into a number of segments equal in length to half the height of the parallelepiped parallel to the normal of that face. With the pseudocode of Algorithm~\ref{alg:ICT} we present the procedure of InflatedCellTruncation, which utilises these results to enumerate the images of unit cells that we include in our energy summation. For our setting, we have included whole unit cell images instead of excluding all terms outside the cutoff sphere $(O,r_{off})$.

\begin{algorithm}[H]
  \caption{Algorithm for the procedure that computes the triplets $n$ and $m$ for energy related calculations}
  \label{alg:ICT}
  \textbf{Input} the transpose of the matrix of lattice vectors $L^T$, the cutoff $r_{off}$ \\
  \textbf{Output} an array of triplets $n_1, n_2, ..., n_{c-1}$ with $n_i=(n_{i,1},n_{i,2},n_{i,3})$
  \begin{algorithmic}[1]
  \Procedure{InflatedCellTruncation}{$L^T, r_{off}$}
    \State $V \gets \text{det}(L^T)$
    \State $O \gets \langle [0.5, 0.5, 0.5],L^T\rangle$
    \State $v_0, v_1, v_2 \gets \text{Normals}(\vect_1, \vect_2, \vect_3)$ \label{opNormals}
    \For{$i=0,1,2$}
    \State $\upsilon \gets V/\|v_i\|$
    \State $t_{(i+2)\%3} \gets \ceil*{(r_{off}-\frac{\upsilon}{2})/\upsilon}$
    \EndFor
    \State $c \gets (2t_0+1) \cdot (2t_1+1) \cdot (2t_2+1) -1$
    \For{$(s_0,s_1,s_2) \gets \text{enumerate}(2t_0+1, 2t_1+1, 2t_2+1)$}
        \If{$(s_0,s_1,s_2) \neq (t_0, t_1, t_2)$}
        \State $n_i \gets (s_0,s_1,s_2) - t$
        \EndIf
    \EndFor
    \Return {$n_1, ..., n_{c-1}$}
    \EndProcedure
  \end{algorithmic}
\end{algorithm}

In Algorithm~\ref{alg:ICT} the function Normals in line~\ref{opNormals} calculates the respective normal vectors for each of the 3 faces defined by $\vect_1, \vect_2, \vect_3$. We only need 3 normals as the symmetrical operations are performed in the next lines for the rest of the faces. This algorithm is used both for the real cutoff with triplets $n$, but also the cutoff in reciprocal space with triplets $m$.

\label{methodology}

During structural optimization, there are two main aspects of the crystal that need to be examined in terms of changes and how these changes affect the structure energy. These are the internal and external coordinates of the crystal. The internal coordinates describe relative positioning of the ions, whilst external coordinates describe the lattice formation. In order to perform energy minimisation, we use the derivatives of the energy function with respect to both kinds of parameters. Here, the internal coordinates' degrees of freedom are $3N$ and correspond to the ion position Cartesian coordinates, while the external coordinates' degrees of freedom are nine and correspond to the lattice vectors. Accordingly, in order to perform optimization, we need two kinds of derivatives of the energy function. The first kind is straightforward, since the parameters' involvement is clearly defined from the definition of the energy potential function.

\subsection{Internal Coordinates}

In this paragraph we present the formulae related to the internal forces acting in a crystal structure, or, more precisely, the forces that act on the ion position vectors. These arise as the negative of the gradient of the energy function $\energy$, whose components are the derivatives $\frac{\partial\energy}{\partial\ionp_t}$ of $\energy$ with respect to each ion position $\ionp_t$. Each such derivative corresponds to a 3-dimensional vector whose components are the partial derivatives of $\energy$ with respect to Cartesian coordinates of $\ionp_t$. Different forms of the position derivatives can be found in different places of literature~\cite{Holden2019AnalyticElectron}, as the ion interaction forces need to be studied for various applications. For the sake of completeness and verification, we calculate and present the derivatives analytically.

\begin{restatable}{prop}{Cforces}
\label{prop:Cforces}
    Let $\energy_{Coul}$ and $\langle \cell, \mathcal{N}, \mathcal{M} \rangle$ be defined as in Proposition~\ref{prop:ECoul}. The internal electrostatic forces $\mathcal{F}_{Coul}$ of $\langle \cell, \mathcal{N}, \mathcal{M} \rangle$ can be written as
    \begin{equation}
    \begin{aligned}\label{eq:FCoul}
        \mathcal{F}_{Coul} &= -\nabla_\ionp\energy_{Coul} \\
                           &= -\left(\frac{\partial \energy_{Coul}}{\partial \ionp_1},\frac{\partial \energy_{Coul}}{\partial \ionp_2},...,\frac{\partial \energy_{Coul}}{\partial \ionp_N}\right)
    \end{aligned}
    \end{equation}
    where
    \begin{equation}\label{eq:partialCoul}
    \begin{gathered}
        \frac{\partial \energy_{Coul}}{\partial \ionp_t} = \frac{k_e}{2}\sum_n\left[   -\sum_{j=1}^{N'} q_tq_j\left( \frac{2\alpha}{\sqrt{\pi}}\exp{(-\alpha^2\|\ionp_{t,j,n}\|^2)}+ \frac{\erfc{\alpha \|\ionp_{t,j,n}\|}}{\|\ionp_{t,j,n}\|}\right)\frac{r_{t,j,n}}{\|\ionp_{t,j,n}\|^2} + \right. \\ 
        \left.\sum_{i=1}^{N'} q_iq_t\left( \frac{2\alpha}{\sqrt{\pi}}\exp{(-\alpha^2\|\ionp_{i,t,n}\|^2)}+ \frac{\erfc{\alpha \|\ionp_{i,t,n}\|}}{\|\ionp_{i,t,n}\|}\right)\frac{r_{i,t,n}}{\|\ionp_{i,t,n}\|^2} \right] + \\ 
        \frac{k_e}{2}\sum_{m}\frac{2\pi k_e}{V\|\rlattice\|^2}\exp{\left(-\frac{\|\rlattice\|^2}{4\alpha^2}\right)}\cdot \left[ -\sum_{j=1}^{N} q_tq_j\sin{(\rlattice r_{t,j})} + \sum_{i=1}^{N} q_iq_t\sin{(\rlattice r_{i,t})}\right]\rlattice, \ t\in[N].
    \end{gathered}
    \end{equation}
\end{restatable}
The proof of Proposition~\ref{prop:Cforces} can be found in Appendix~\ref{app3}.

\begin{restatable}{prop}{Bforces}
\label{prop:Bforces}
   Let $\energy_{Buck}$ and $\langle \cell, \mathcal{N}, \mathcal{M} \rangle$ be defined as in Proposition~\ref{prop:EBuck}. The internal Buckingham forces $\mathcal{F}_{Buck}$ can be written as
    \begin{equation}\label{eq:FBuck}
    \begin{aligned}
        \mathcal{F}_{Buck} &= -\nabla_\ionp\energy_{Buck} \\
                           &= -(\frac{\partial \energy_{Buck}}{\partial \ionp_1},\frac{\partial \energy_{Buck}}{\partial \ionp_2},...,\frac{\partial \energy_{Buck}}{\partial \ionp_N})
    \end{aligned}
    \end{equation}
    where
    \begin{equation}\label{eq:partialBuck}
    \begin{gathered}
    \frac{\partial \energy_{Buck}}{\partial \ionp_t} = \\
    \frac{1}{2}\sum_{n}\Bigg\{ \sum_{j=1}^{N'}  \left[ -\frac{A_{tj}}{\rho}\exp{\left(-\frac{\|\ionp_{t,j,n}\|}{\rho}\right)} - C_{tj}\frac{\exp{(-\alpha^2\|\ionp_{t,j,n}\|^2)}}{\|\ionp_{t,j,n}\|^5}\cdot \right. \\
    \left.\left. \left( \frac{6}{\|\ionp_{t,j,n}\|^2}+6\alpha^2+\alpha^6\|\ionp_{t,j,n}\|^4+3\alpha^4\|\ionp_{t,j,n}\|^2 \right) \right] + \right. \\ 
    \sum_{i=1}^{N'}  \left[ \frac{A_{it}}{\rho}\exp{\left(-\frac{\|\ionp_{i,t,n}\|}{\rho}\right)} + \left. C_{it}\frac{\exp{(-\alpha^2\|\ionp_{i,t,n}\|^2)}}{\|\ionp_{i,t,n}\|^5} \right.\right. \\
    \left. \left( \frac{6}{\|\ionp_{i,t,n}\|^2}+6\alpha^2+\alpha^6\|\ionp_{i,t,n}\|^4+3\alpha^4\|\ionp_{i,t,n}\|^2 \right) \right] \Bigg\}\frac{\ionp_{i,t,n}}{\|\ionp_{i,t,n}\|} - \\
    \frac{1}{2}\cdot\frac{\pi^{3/2}}{12V}\sum_{m} \Bigg\{ \sum_{j=1}C_{tj}\left[ \sqrt{\pi}\cdot\erfc{\frac{\|\rlattice\|}{2\alpha}} + \right. \\
    \left.\left. \left( \frac{4\alpha^3}{\|\rlattice\|^3}-\frac{2\alpha}{\|\rlattice\|}\right) \exp{\left(-\frac{\|\rlattice\|^2}{4\alpha^2}\right)} \right]\cdot \sin(\rlattice\ionp_{t,j,n}) + \right. \\ 
    \left. \sum_{i=1}C_{it}\left[ \sqrt{\pi}\cdot\erfc{\frac{\|\rlattice\|}{2\alpha}} + \right.\right. \\
    \left. \left( \frac{4\alpha^3}{\|\rlattice\|^3}-\frac{2\alpha}{\|\rlattice\|}\right) \exp{\left(-\frac{\|\rlattice\|^2}{4\alpha^2}\right)} \right]\cdot\sin(\rlattice\ionp_{t,j,n})\Bigg\}\rlattice, \ t\in[N].
    \end{gathered}
    \end{equation}
\end{restatable}

The proof of Proposition~\ref{prop:Bforces} can be found in Appendix~\ref{app3}.
The sum of Equations~(\ref{eq:partialCoul}),~(\ref{eq:partialBuck}) constitutes the partial derivative of the overall energy potential function $\energy$ with respect to ion position $\ionp_t, \ t\in [N]$. Hence, each such derivative is one of the N components of
\begin{equation}
      \mathcal{F} = -\nabla_\ionp\energy =\mathcal{F}_{Coul} + \mathcal{F}_{Buck}
\end{equation}
that correspond to the internal coordinate forces.

\subsection{External Coordinates}

When performing geometric optimization to a crystal structure, its initial state is assumed to be in a state of agitation. This means that there are two factors preventing the crystal state to reach equilibrium. Since the crystal lattice is in a shape that does not correspond to the equilibrium state, it is deformed. Thus, apart from using the previously listed derivatives with respect to ion positions, we need to change the lattice vectors and express the forces acting on the vectors in such a way, so as to combine position derivatives and lattice forces into one updating step. In order to preserve properties of the material's continuum~\cite{Crandall1960AnSolids} and avoid rigid body movements, we use the symmetrical strain tensor $\epsilon$ to update the lattice vectors. Then at each step of the optimization the lattice vectors are characterised by a state of strain $0,\epsilon_1,...,\epsilon_{n}$.

It is important to highlight here that, since we are changing the lattice cell parameters, the unit cell volume is also affected. Thus derivatives of the unit cell volume $V$ must also be defined. Equation~(\ref{eq:alpha}) suggests that the $\alpha$ parameter \textbf{is a function of the volume}. As a consequence,  have that
\begin{equation}\label{eq:alphadrv}
    \frac{\partial\alpha}{\partial\epsilon_{\lambda\mu}} = \frac{\partial\alpha}{\partial V}\frac{\partial V}{\partial\epsilon_{\lambda\mu}} = \alpha'(V)\delta_{\lambda\mu}V, \quad \alpha'(V) = -\frac{\sqrt[6]{N}\sqrt{\pi}}{3\sqrt[3]{V^4}}.
\end{equation}
We, therefore, stress the fact that some extra terms arise in the strain derivatives owed to the existence of this function in place of a constant parameter. The multiplier $\delta_{\lambda\mu}$ can be used as an indication of which these terms are in Propositions~\ref{prop:Cstress},~\ref{prop:Bstress}. We hereafter refer to $a'(V)$ as $a'$.

\begin{restatable}{prop}{stresstheo}
\label{prop:stresstheo}
    The forces that act on the volume of a unit cell $\cell_n$ and change the shape of the crystal lattice can be expressed with the symmetric stress tensor $\sigma=\{\sigma_{\lambda\mu}\}_{\lambda,\mu\in[3]}$ as a result of an existing strain $\epsilon=\{\epsilon_{\lambda\mu}\}_{\lambda,\mu\in[3]}$.
\end{restatable}
The proof of Proposition~\ref{prop:stresstheo} can be found in Appendix~\ref{app3}.

\begin{restatable}{lemma}{stresscalc}
\label{lem:stresscalc}
    Let $\energy$ be defined as in Equation~(\ref{eq:energyall}). The parameters affected by stress are $\allionp,\lattice,\rlattice$, hence each stress component $\sigma_{\lambda\mu}$ is calculated as
\begin{equation}
\begin{gathered}
    \sigma_{\lambda\mu} = \frac{1}{V} \left( \sum_{i=1}^N\frac{\partial\energy}{\partial\ionp_{i\lambda}}\ionp_{i\mu} + \sum_{n}\frac{\partial\energy}{\partial\lattice_{\lambda}}\lattice_\mu \right. \\
    \pushright{\left. + \sum_{m}\frac{\partial\energy}{\partial\rlattice_{\mu}}\rlattice_\lambda + \frac{\partial\energy}{\partial V}\delta_{\lambda\mu}V \right)}.
\end{gathered}
\end{equation}
\end{restatable}

The proof of Lemma~\ref{lem:stresscalc} can be found in Appendix~\ref{app3}.
The overall stress can be easily calculated by separately evaluating the stress produced by Coulomb and Buckingham stresses so that 
\begin{equation}
    \sigma_{\lambda\mu} = \sigma_{(Coul)\lambda\mu}+\sigma_{(Buck)\lambda\mu}
\end{equation}
Then the lattice vectors $L$ can be updated with the procedure listed in Algorithm~\ref{alg:RLupdate}.
\begin{algorithm}[H]
  \caption{Algorithm of parameter update}
  \label{alg:RLupdate}
  \textbf{Input} step size s, direction vector $d$, transpose of ion positions $\allionp^T$, transpose of lattice vectors $L^T$, strain tensor $E$\\
  \textbf{Output} transpose of ion positions $\allionp^T$, transpose of lattice vectors $L^T$ and strain matrix E
  \begin{algorithmic}[1]
    \Function{Update}{$s,d,\allionp^T, L^T, E$}
    \State $\allionp^T \gets \allionp^T + s\cdot d[1,..,N]$ \Comment{Ion positions' update}
    \State $E[DU] \gets E[DU] + s\cdot d[N+1,...,N+6]$ \Comment{Strain update}
    \State $E[L^T] \gets E[U]$ \Comment{Render strain matrix symmetric}
    \State $\Delta \gets (E-J_3)+I_3$
    \State $L^T \gets L^T\Delta^T$ \Comment{Apply strains to lattice vectors}
    \State $\allionp^T \gets \allionp^T\Delta^T$ \Comment{Apply strains to ion vectors}
    \State \Return {$\allionp^T, L^T, E$}
    \EndFunction
  \end{algorithmic}
\end{algorithm}
In Algorithm~\ref{alg:RLupdate} we denote with $L$,$D$,$U$ the matrix decomposition in lower triangular, diagonal and upper triangular parts, so that matrix $E$ will be populated with the six values $d[N+1,..,N+6]$ according to the Voigt notation. Moreover, $R$, $L$, $J_3$, $I_3$ stand for the following matrices 
\begin{gather*}
    R = \begin{bmatrix}
        \vertbar & \vertbar &        & \vertbar \\
        \ionp_1 & \ionp_2, & ..., & \ionp_N \\
        \vertbar & \vertbar &        & \vertbar 
    \end{bmatrix},
    L = \begin{bmatrix}
        \vertbar & \vertbar & \vertbar \\
        \vect_1 & \vect_2 & \vect_3 \\
        \vertbar & \vertbar & \vertbar 
    \end{bmatrix},
    I_3 = \begin{bmatrix}
        1 & 0 & 0 \\
        0 & 1 & 0 \\
        0 & 0 & 1
    \end{bmatrix}, 
    J_3 = \begin{bmatrix}
        1 & 1 & 1 \\
        1 & 1 & 1 \\
        1 & 1 & 1
    \end{bmatrix}
\end{gather*}

In order to avoid a large percentage of Buckingham catastrophes, a strain reset must be in place. This means that there must be an interval after which the current configuration is assumed to be the initial one, the crystal structure without any stress present. This is achieved by setting the components of strain matrix $E$ to 1, so that $E=I_3$, every $3N+9$ iterations. We selected this interval as a reset point according to the number of parameters of the problem~\ref{eqn:Problem} and the least steps possible to relax the structure, however its duration is yet to be optimised.

By Proposition~\ref{prop:stresstheo} and Lemma~\ref{lem:stresscalc}, Proposition~\ref{prop:Cstress} and Proposition~\ref{prop:Bstress} arise naturally. This can be seen in their proofs located at Appendix~\ref{app3}.
\begin{prop}
\label{prop:Cstress}
    Let $\energy_{Coul}^S$, $\energy_{Coul}^L$, $\energy_{Coul}^{self}$ be defined as in Proposition~\ref{prop:ECoul}. The stress applied on the unit cell volume due to Coulomb forces is calculated as
    \begin{equation}
    \begin{aligned}
         \frac{\partial \energy^S_{Coul}}{\partial \epsilon_{\lambda\mu}} &= \sum_{i,j,n}^{N'}q_iq_j\left[ k_e\frac{-\alpha'V}{\sqrt{\pi}}\exp{(-\alpha^2\|\ionp_{i,j,n}\|^2)\delta_{\lambda\mu}} +  f_S'(\|\ionp_{i,j,n}\|)\frac{\ionp_{i,j,n(\lambda)}}{\|\ionp_{i,j,n}\|}\ionp_{{i,j,n}(\mu)} \right] \\
         \frac{\partial \energy^L_{Coul}}{\partial \epsilon_{\lambda\mu}} &= \frac{2\pi k_e}{V}\sum_{i,j,m}q_iq_j\frac{\exp{\left(-\frac{\|\rlattice\|^2}{4\alpha^2}\right)}}{\|\rlattice\|^2}\cos{(\rlattice\ionp_{i,j})}\cdot \\
         & \pushright{\left[ \left( \frac{1}{2\alpha^2}+\frac{2}{\|\rlattice\|^2}\right)\rlattice_{\mu}\rlattice_{\lambda} - \delta_{\lambda\mu}\left( 1-\frac{\|\rlattice\|^2}{2\alpha^3}\alpha'V \right) \right]} \\
         \frac{\partial \energy^{self}_{Coul}}{\partial \epsilon_{\lambda\mu}} &=   -\frac{\alpha'k_e}{\sqrt{\pi}}\sum_i^Nq_i^2\cdot \delta_{\lambda\mu}V
    \end{aligned}
    \end{equation}
    and ultimately the stress tensor component is
    \begin{equation}
        \sigma_{(Coul)\lambda\mu} = \frac{1}{V}\left( \frac{\partial \energy^S_{Coul}}{\partial \epsilon_{\lambda\mu}} + \frac{\partial \energy^L_{Coul}}{\partial \epsilon_{\lambda\mu}} - \frac{\partial \energy^{self}_{Coul}}{\partial \epsilon_{\lambda\mu}} \right).
    \end{equation}
\end{prop}

\begin{prop}
\label{prop:Bstress}
    Let $\energy_{Buck}^S$, $\energy_{Buck}^L$, $\energy_{Buck}^{self}$ be defined as in Proposition~\ref{prop:EBuck}. The stress applied on the unit cell volume due to Buckingham forces is calculated as
    \begin{equation}
    \begin{aligned}
        \frac{\partial \energy^{S}_{Buck}}{\partial \epsilon_{\lambda\mu}} &= -\frac{1}{2}\sum^{N'}_{i,j,n}\frac{C_{ij}}{\|\ionp_{i,j,n}\|^6}\exp{(-\alpha^2\|\ionp_{i,j,n}\|^2)}\cdot  \\
        & \bigg[\left( \frac{6}{\|\ionp_{i,j,n}\|^2} + 6\alpha^2+\alpha^6\|\ionp_{i,j,n}\|^4+3\alpha^4\|\ionp_{i,j,n}\|^2 \right)r_{{i,j,n}(\lambda)}r_{{i,j,n}(\mu)} + \alpha'\alpha^5V\|\ionp_{i,j,n}\|^6\delta_{\lambda\mu} \bigg] \\
        \frac{\partial \energy^{L}_{Buck}}{\partial \epsilon_{\lambda\mu}} &= \frac{1}{2}\sum_{i,j}^N C_{ij}\frac{\pi^{3/2}}{12V}\sum_{m} \cos{(\rlattice\ionp_{i,j})}\cdot \\
        & \left[ \left(3\sqrt{\pi}\|\rlattice\|\erfc{\frac{\|\rlattice\|}{2\alpha}}-6\alpha\exp{\left(-\frac{\|\rlattice\|^2}{4\alpha^2}\right)}\right)\rlattice_\mu\rlattice_\lambda  \right. - \\ 
        & \left. \delta_{\lambda\mu}\left( -\sqrt{\pi}\erfc{\frac{\|\rlattice\|}{2\alpha}}\|\rlattice\|^3 + (-2\alpha^2+\|\rlattice\|^2+6V\alpha\alpha')2\alpha\exp{\left(-\frac{\|\rlattice\|^2}{4\alpha^2}\right)} \right) \right] \\
        \frac{\partial \energy^{self}_{Buck}}{\partial \epsilon_{\lambda\mu}} &=  \delta_{\lambda\mu} \bigg( -\frac{1}{2}\sum_{i,j}^NC_{ij}\frac{\pi^{3/2}\alpha^2}{3V}\cdot(3\alpha'V-\alpha)+\frac{1}{2}\sum_i^NC_{ii}\alpha^5\alpha'V \bigg)
    \end{aligned}
    \end{equation}
    and ultimately the stress tensor component is
    \begin{equation}
        \sigma_{(Buck)\lambda\mu} = \frac{1}{V}\left( \frac{\partial \energy^S_{Buck}}{\partial \epsilon_{\lambda\mu}} + \frac{\partial \energy^L_{Buck}}{\partial \epsilon_{\lambda\mu}} - \frac{\partial \energy^{self}_{Buck}}{\partial \epsilon_{\lambda\mu}} \right).
    \end{equation}
\end{prop}
This concludes the presentation of the energy function and its first order derivatives, along with the necessary theory that accompanies the derivations. In the next we emphasize on the experimentation process starting with a description of the employed algorithms.

\subsection{First order algorithms}
\label{algorithms}

We investigate the performance of Steepest Descent (aka Gradient Descent) and Conjugate Gradient. The choice of these algorithms is based on the prospect of selecting stable recipes that are sure to converge to the local energy minimum and investigate how parameter tuning can affect the relaxation. More specifically, we are interested in comparing the stability and convergence speed of their respective  updating schemes. \textit{Gradient Descent} is a simple method relying on the function gradient to define the search direction of the optimization. Let $\mathcal{F}:\mathbb{R}^n\rightarrow\mathbb{R}$ a multivariate differentiable function.  Then the negative of the gradient vector $\nabla\mathcal{F}$ determines the direction with maximum decrease. For the second part of our experiments we use \textit{Conjugate Gradient}. This is a conjugate directions method for nonlinear problems, meaning that every produced direction is targeted to be conjugate to all previous directions. Its nonlinear version is made possible using the Gram-Schmidt orthonormalising process to create the search directions. We use the Polak–Ribière method~\cite{Polak1969NoteConjuguees} to update the direction vector, which is proven to have good performance in various similar problems. It is an inherently restarting method~\cite{Powell1977RestartMethod} that avoids repeatedly small steps when the direction vector is almost orthogonal to the residual of the function.

For our experiments, the input variables are the positions of the ions and the strains of the lattice vectors. The strain tensor is the infinitesimal version of the change seen in solids under acting forces, as found in~\cite{Crandall1960AnSolids}. The tensor is used in place of the actual lattice vectors in order to maintain the symmetry of the unit cell and restrict rigid body movements.

\section{Experimental Setting} \label{eexperiments}

In this section we will introduce the experimentation process. We first describe the purpose of our experiments and provide an outline of what is presented in the results. Afterwards, we include a technical description of the input and software used for their execution.

\subsection{Description of experiments}\label{experimentsdescr}
Little work has been done on studying the structural relaxation as a Computer Science algorithmic application and formally exhibiting its properties and requirements. Gradient Descent is usually the to-go algorithm for any non-convex local optimization problem, as it is intuitively easy to understand and it can be robust enough to eventually lead to the minimiser using only the first, and easiest to compute, derivatives. However, our results suggest that Conjugate Gradient is a valuable alternative optimization algorithm which is not only more trustworthy than Gradient Descent, but also generally quicker to converge. Thus, we test the following hypotheses:
\begin{itemize}
    \item Conjugate Gradient is more robust than Gradient Descent.
    \item Careful step size selection improves the optimization's performance.
     \item There is a trade-off between convergence speed and robustness.
\end{itemize}

We test each of the above hypotheses experimentally and we provide an enhanced algorithmic recipe reflecting the conclusions drawn from our results. We trial the convergence and speed of Gradient Descent and Conjugate Gradient under 4 step size adaptations: a) constant step size, b) exponential scheduled step size (\textit{es}), c) constant scheduled bisection (\textit{bisect}) for the step size and c) gradient-norm-related scheduled (\textit{gbisect}) step size. We examine the impact of these adaptations under two conditions, namely, the number of successful experiments and the runtime, in terms of number of iterations.  We, then, analyse the behaviour of all step size recipes under the same conditions and select the optimal per method. The results of the experiments will give rise to a \textbf{trade-off}, whereby, step size recipes that reduce the runtime result into more failed experiments. This becomes particularly prominent with constant step sizes, as, the larger the step size gets, the more the number of experiments that do not finish before a set deadline increases. Therefore, a function will be proposed to measure the utility of each algorithm with each step size adaptation.

\subsection{Technical Information} \label{data}
For our experiments we have used a set of 200 crystal structures 
produced with a stable Strontium Titanate (\ch{Sr_3Ti_3O_9}) as a reference point and the introduction of randomness to  the unit cell.  More specifically, after defining the length of each lattice vector from a set of values of 4, 6, 8, 10, and 12 {\AA}, an orthorhombic unit cell is formed and 15 ions -- 3 strontium, 3 titanium and 9 oxygen ions -- are placed in a random manner on grid points defined by a 1 {\AA}  grid spacing. The placement is such that the negative ions are placed on grid points with even indices, and positive ions and are placed on grid points with odd indices. This construction method provides input elements that are not likely to lie near the PES minimum, thus allowing to test if an algorithm can find the minimum despite it being far away. The 200 structures were divided into 5 groups of 40 randomly to provide an unbiased basis for statistical analysis of the results.

Our implementation (found in https://github.com/lrcfmd/veltiCRYS) offers the aforementioned energy and forces calculations as in Equations~(\ref{eq:potentialLong})-(\ref{eq:EBuckEwald}) using the \textit{ICT} method. Input can be defined using ASE's\cite{Hjorth_Larsen_2017} \textit{Atoms} class or read as a \textit{CIF} file\cite{Hall:es0165}. Other ASE tools for geometry and input-output tools have also been used. The implementation's backbone is written in Cython 0.29.30 and parts such as input handling are written in Python 3.10.4. The output can be configured to extract PNG images and CIF files for each produced crystal structure configuration every requested number of iterations. The experiments were run on a Intel Xeon Gold Skylake processor with 9.6 GB of memory per core with Linux.

\section{Results and discussion}
\label{results}

Every structure, needing several structural modifications to approach equilibrium, underwent a procedure with which each minimization iteration corresponds to two parameter updates. Firstly, a displacement of all ions in the unit cell and, secondly, a length and angle adaptation of the unit cell vectors; in other words, a structural relaxation iteration. In simple terms, with each parameter update we move the ions $\allionp$ in the unit cell, then we stretch or shrink the lattice vectors $L$ and change the three angles in between the lattice vectors. The relaxation stops when one of the following occurs:
\begin{enumerate}
    \item \textbf{The gradient norm $g$ has fallen below the tolerance value $tol = \num{0.001}$}, and the resulting potential energy $\energy$ is less than the energy $\energy_0$ of the initial configuration $\allionp_0,L_0$. We refer to the following as the gradient norm
    \begin{equation} \label{eq:gnorm}
        g(\allionp,L) = \frac{
        \sqrt{
            \sum_{i=1}^{N}\left((\frac{\partial \energy}{\partial \ionp_{ix}})^2 + (\frac{\partial \energy}{\partial \ionp_{iy}})^2 + (\frac{\partial \energy}{\partial \ionp_{iz}})^2\right) + 
            \sum_{i=1}^6\frac{\partial \energy}{\partial \epsilon_{i}}^2
        }
        }{3N+6}
    \end{equation}
    and we announce a successful relaxation when 
    \begin{equation}
    \label{eq:termination}
    g < \num{0.001} \text{ and } \energy < \energy_0.
    \end{equation}
    Since our experiments' input is constructed in a way such that $\allionp_0,L_0$ is not very ``close" to the stable \ch{Sr_3Ti_3O_9}, we allow for a margin of error $\epsilon$ which enables our methods to converge to the \textit{true} local minimum of $\allionp_0,L_0$.

    \item \textbf{Buckingham catastrophe} happens, which results into a constantly increasing gradient norm and decreasing energy value. Further description can be found in~\ref{buckcat}.

    \item The conditions~(\ref{eq:termination}) are not fulfilled and the \textbf{iteration number $i$ has reached 50000 iterations}, which we describe as \textit{overtime}. This result indicates that the step size magnitude is large enough to prevent convergence below the selected tolerance in a reasonable time margin and cannot guarantee a finite sequence of iterations. We observed that when the step size is too big the gradient and energy started to oscillate and their values could not decrease below some threshold in a sensible amount of time. 

    \item This case arose only for the category of experiments on the largest constant step size and concerns only a minority of the dataset, while the rest of the structures of this category of experiments came to failure due to the aforementioned reasons. This is when, in spite of the iteration number $i<I$, \textbf{the experiment has not completed in the time interval of 3 days after its commencement}. In this situation we observe that some of the lattice vectors increase abnormally in size, causing a lot of calculations with large numbers. This slows down the experiment's progress, while moving away from the expected local minimum and hence eventually not creating the conditions for convergence.   
\end{enumerate}

\subsection{Trivial case benchmark}
\label{experimentsI}

\begin{figure}
\centering
    \begin{subfigure}{0.49\linewidth}
    \includegraphics[width=0.9\textwidth]{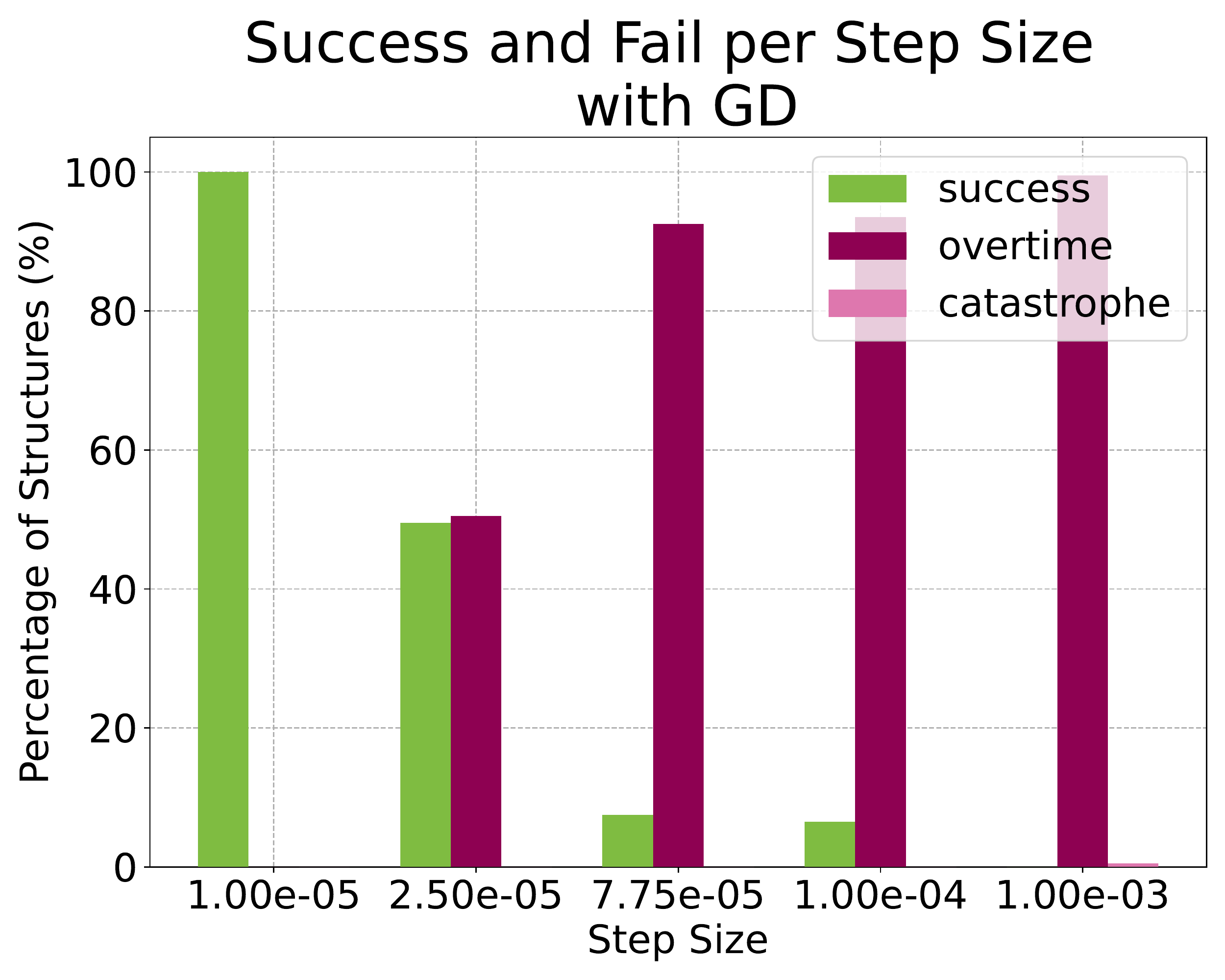}
    \caption{Gradient Descent success rate across step sizes}
    \label{fig:stepsuccgd}
    \end{subfigure}
    \begin{subfigure}{0.49\linewidth}
    \includegraphics[width=0.9\textwidth]{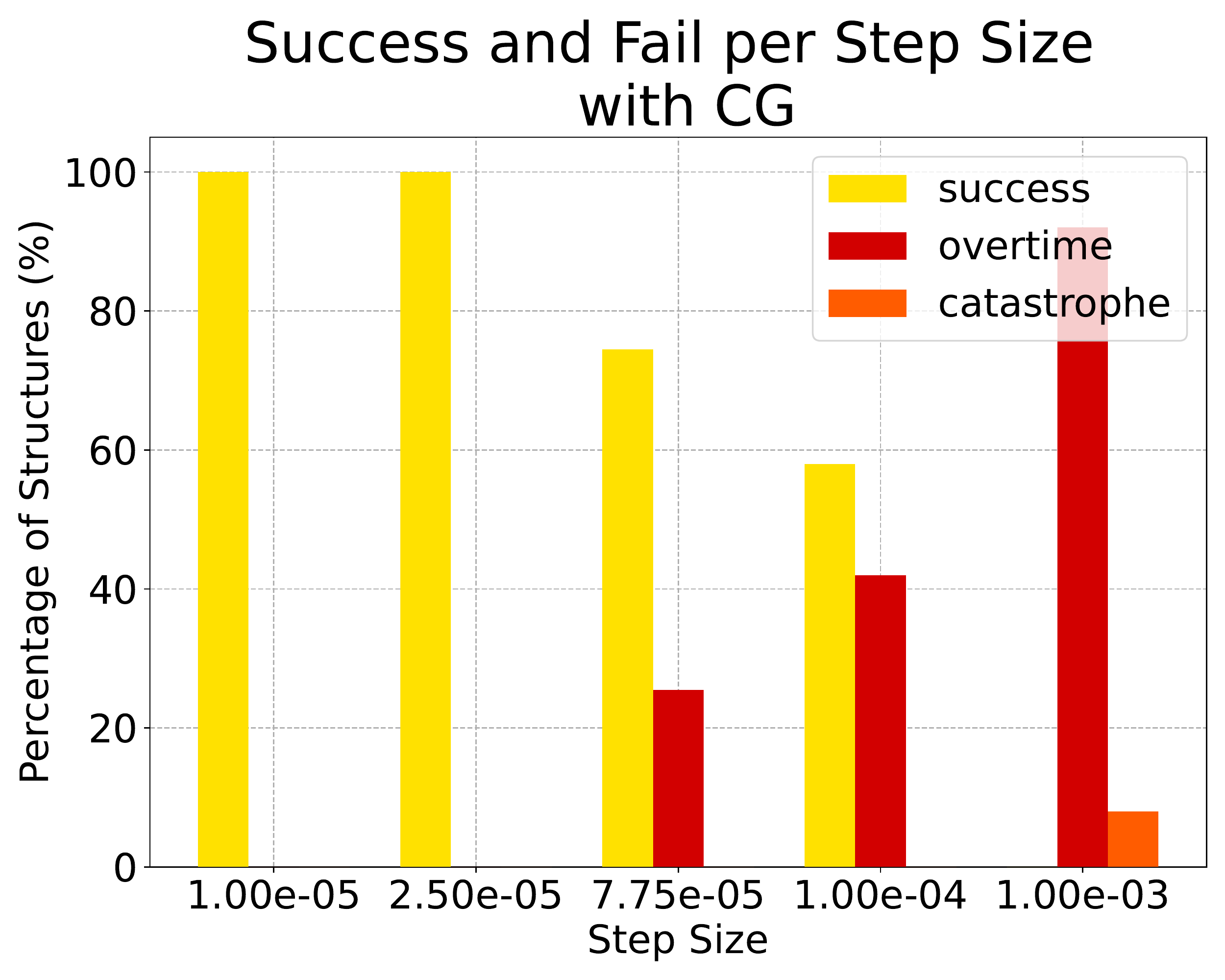}
    \caption{Conjugate Gradient success rate across step sizes}
    \label{fig:stepsucccg}
    \end{subfigure}
    \caption[Percentage of success for relaxations with constant step sizes]{\textbf{Percentage of success over all 200 structures for relaxations with constant step sizes of $\mathcal{S}$.} The first bar on the left per step size shows the percentage of 200 structures that were successfully relaxed and the next two bars per step size show the ones that failed due to overtime running -- reached the budget of iterations -- or Buckingham catastrophe.}
    \label{fig:stepsucc}
\end{figure}

We test five values for the constant step size, namely $\num{0.00001}$ (small), $\num{0.000025}$, $\num{0.0000775}$, $\num{0.0001}$ (medium) and $\num{0.001}$ (large). Our choice of the lower step size value was determined empirically and stemmed out of the hypothesis that there is at least one local minimiser in the feasible neighbourhood around $\allionp_0, L_0$. This is a justified hypothesis taking into account that we already know a stable structure with similar configuration to each $(\allionp_0, L_0)$, the stable \ch{Sr_3Ti_3O_9}. After several experiments we achieved a successful relaxation for \textbf{100\%} of the dataset using a constant $s = \num{0.00001}$ for both examined methods, as seen in Figure~\ref{fig:stepsucc}. The results from these runs showed that a successful relaxation would terminate in fewer than 50000 cycles for all structures, consequently, we set our limit for iterations to this number. Figures~\ref{fig:stepsuccgd},~\ref{fig:stepsucccg} illustrate that, both for Gradient Descent and Conjugate gradient, 0\% of instances converged with $s = \num{0.001}$. Thus, $\num{0.001}$ will be the highest value of the step size that we test.  In order to investigate the algorithms' behaviour for step sizes in the range $(\num{0.00001},\num{0.001})$, we complete a set $\mathcal{S}$ of step size values to test with  $\num{0.000025}$, $\num{0.0000775}$ and $\num{0.0001}$ so that
\begin{equation}
    \label{eq:setS}
    \mathcal{S} = \{\num{0.00001}, \num{0.000025}, \num{0.0000775}, \num{0.0001}, \num{0.001}\}.
\end{equation}

\subsubsection{Small step size}
When the step size is constantly very small with unchanged value, Gradient Descent and Conjugate Gradient have similar performance. The smallest step size $s = \num{0.00001}$ cancels the conjugacy between directions produced from Conjugate Gradient, as expected in cases of employment of the Polak–Ribière update~\cite{Grippo1997AMethod}. This happens because when the step size value approaches zero $s\rightarrow 0$ then the difference between consecutive gradients decreases $\gamma_i - \gamma_{i-1} \rightarrow 0$ and the direction vector approaches the negative gradient $d_i\rightarrow -\gamma_i$. As a result, Conjugate Gradient acts as Gradient Descent and the two methods appear to have almost \textit{no differences} in consecutive steps.

\subsubsection{Medium step size}
For this range of step sizes, a large percentage of Gradient Descent experiments fail, but only up to half structures fail with Conjugate Gradient and the relaxations are accelerated as the step size increases. When $\num{0.000025}$, almost half of the Gradient Descent experiments fail, but all experiments of Conjugate Gradient are successfully completed and finish in fewer than half the number of iterations from before. When $s=\num{0.0000775}$ and $s=\num{0.0001}$,  Gradient Descent fails to converge for 98\% of the structures and Conjugate Gradient failure numbers can reach almost $50\%$ of the structures. Nonetheless, the \textit{successful} relaxations are accelerated to a mean $\approx 3700$ iterations per optimization run, resulting in a 7 times speedup from the initial experiments. Ultimately, a trade-off between success rate and convergence speed is evident, but the high failure rate of Gradient Descent renders it useless for step size values $>\num{0.000025}$.  

\subsubsection{Large step size}
For this step size value \textit{all} structures fail to relax with a budget of 50000 iterations for either of the two methods. Following our previous observations, we find that almost all failures are caused due to the large step size that fails to accurately approach the small neighbourhood of the local minimum in $I=50000$ iterations. Here we must report that Conjugate Gradient with this constant step leads several structures to unit cells with very large pairwise distances, so this confirms that the \textit{true} local minimum was missed and another extremum is being followed. All structures that were run with Gradient Descent, except for one that fell into Buckingham catastrophe, failed due to exhausting the iteration budget I. We come across with a frequent \textbf{Buckingham catastrophe} when using Conjugate Gradient with $s=\num{0.001}$.  In this case at least one experiment out of 40 falls into Buckingham catastrophe creating a 2.5\% to 5\% possibility for a structure to fail because of it. Therefore, we see that large displacements push ion components together into an energy regime that does not allow the procedure to recover from pacing towards $-\infty$. After all, it has been shown in literature~\cite{Dai2010ConvergenceStepsizes} that large values in constant step size can cause problems to Conjugate Gradient's convergence and lead to uphill directions. 

\subsubsection{Some further results}
\begin{figure}
    \begin{subfigure}[t]{0.49\linewidth}
    \includegraphics[width=0.95\linewidth]{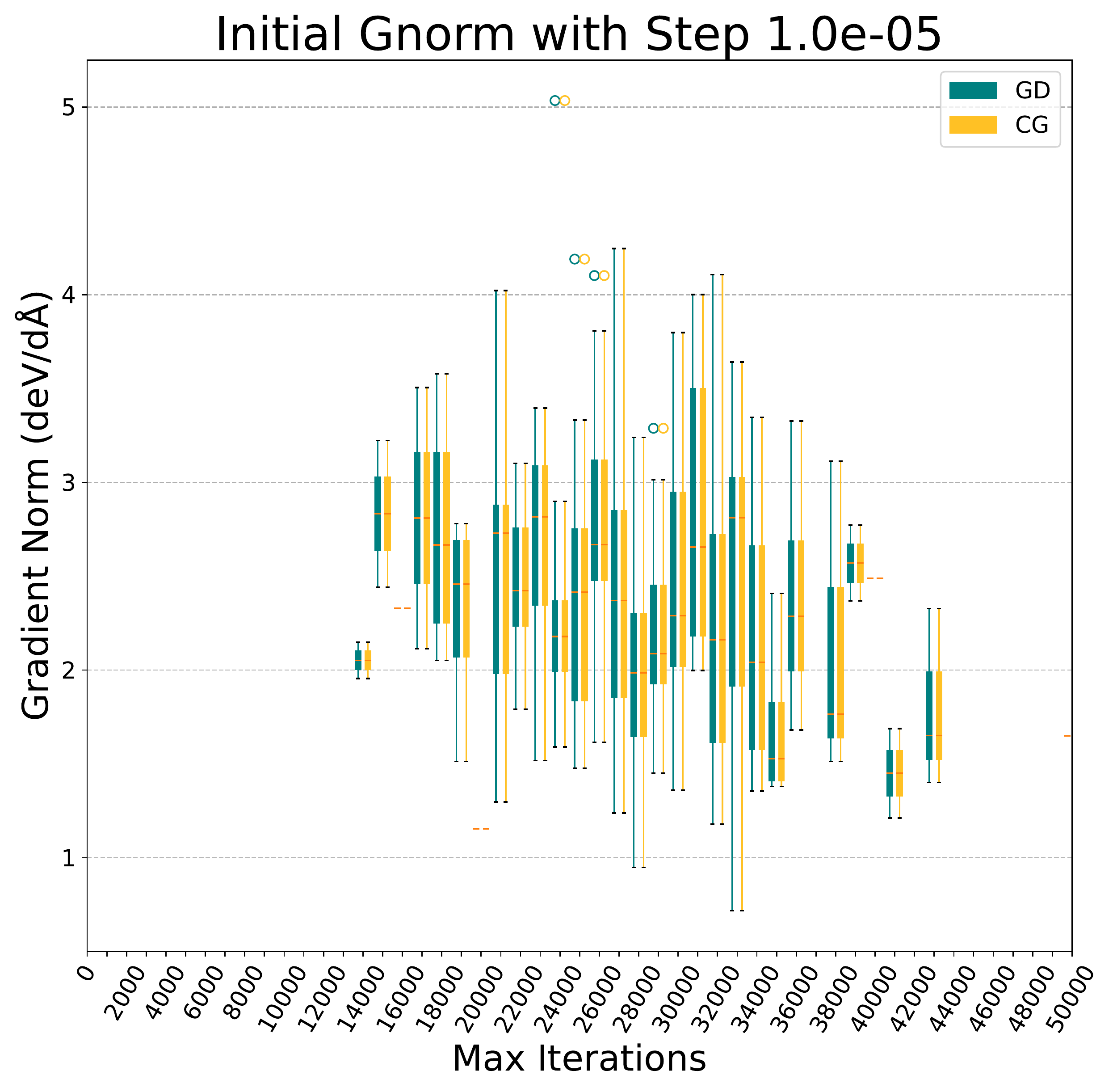}
    \caption{Initial gradient norm versus (rounded to 1000) total iteration number for successful relaxations with constant step $s=\num{0.00001}$.}
    \label{fig:initgnorm1e5}
    \end{subfigure}
    \begin{subfigure}[t]{0.49\linewidth}
    \includegraphics[width=0.95\linewidth]{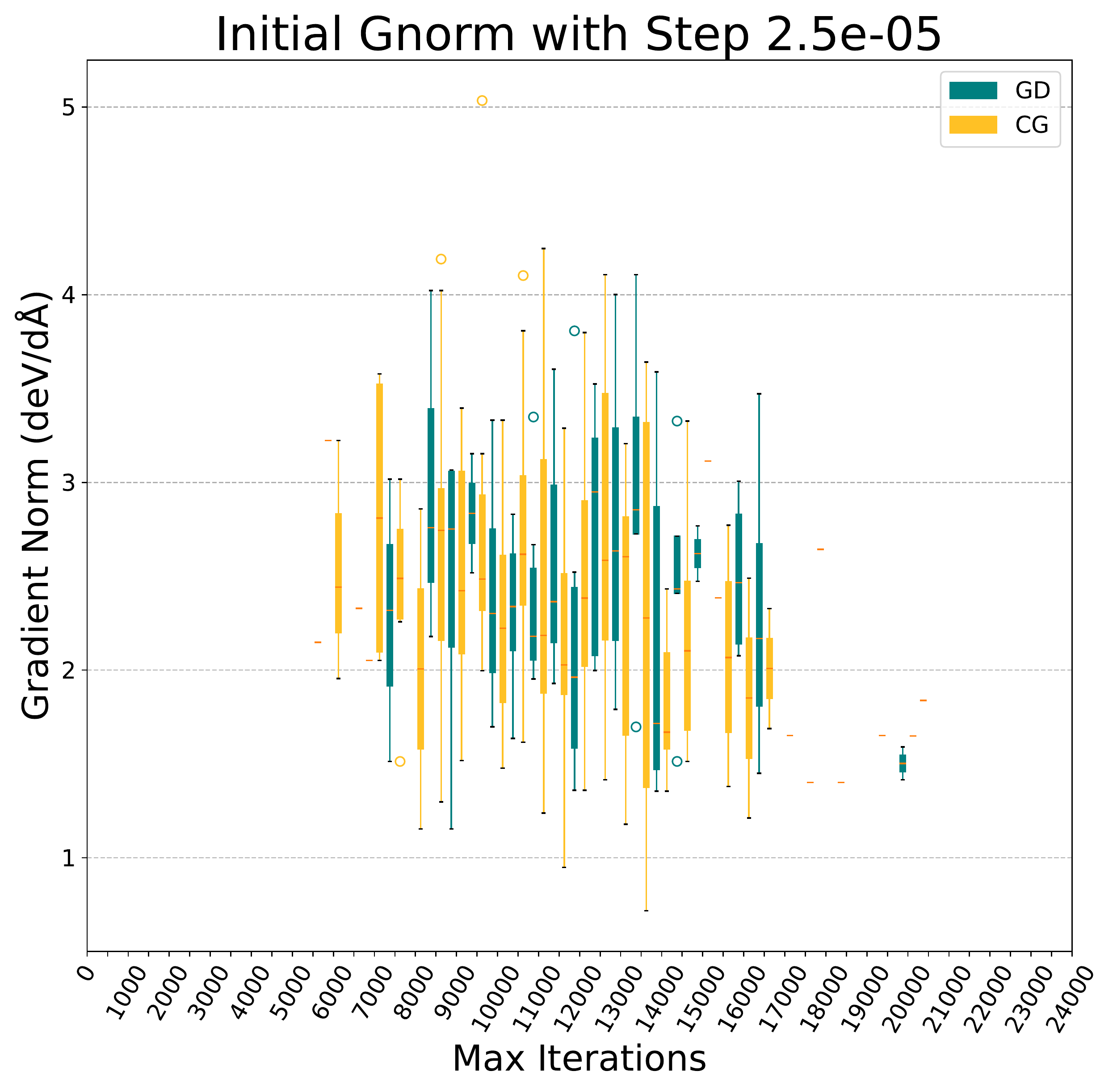}
    \caption{Initial gradient norm versus (rounded to 1000) total iteration number for successful relaxations with constant step $s=\num{0.000025}$.}
    \label{fig:initgnorm2e5}
    \end{subfigure}
    \medskip
    \begin{subfigure}[t]{0.49\linewidth}
    \includegraphics[width=0.95\linewidth]{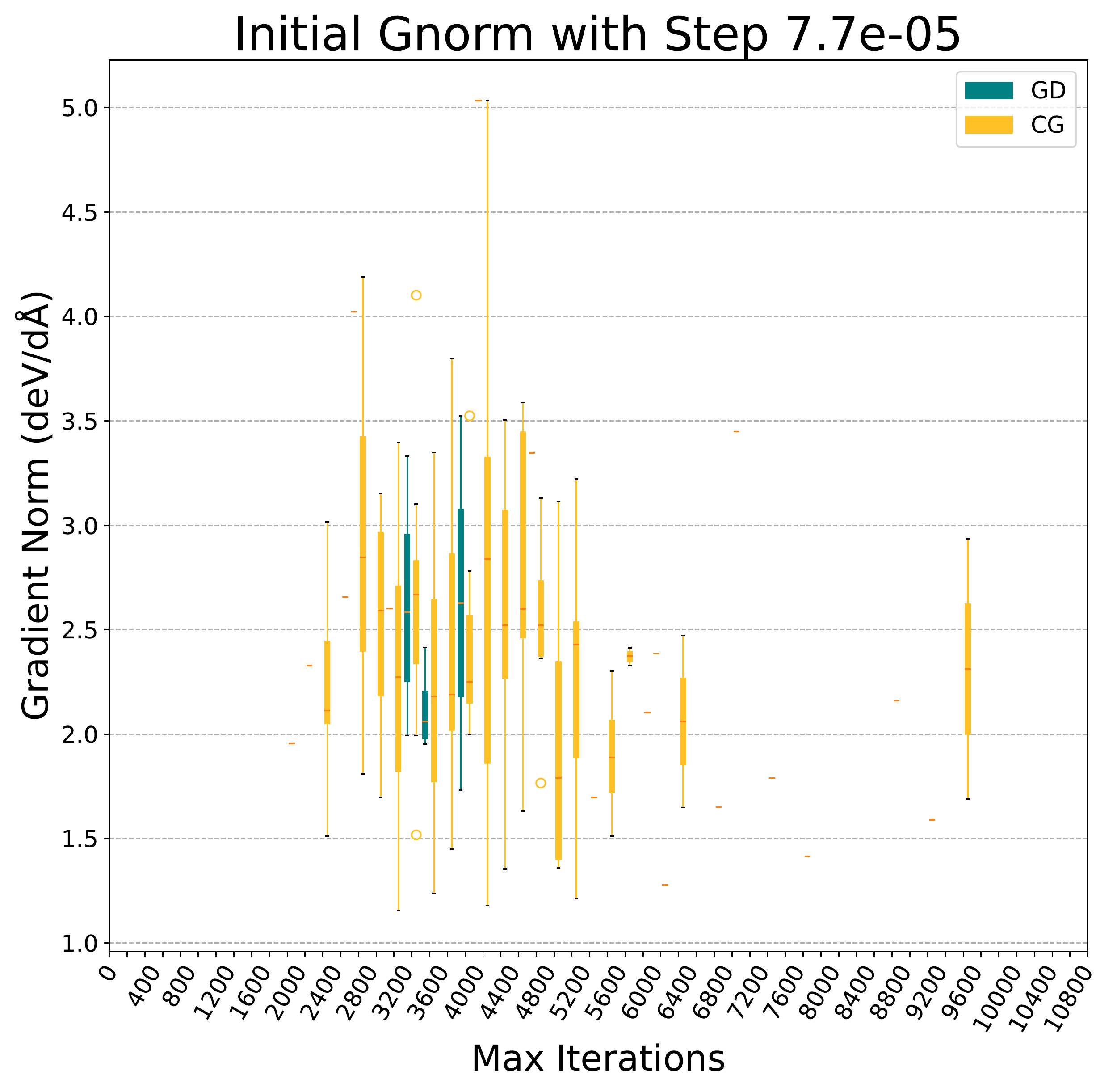}
    \caption{Initial gradient norm versus (rounded to 1000) total iteration number for successful relaxations with constant step $s=\num{0.0000775}$.}
    \label{fig:initgnorm7e5}
    \end{subfigure}
     \hspace{0.2cm}
    \begin{subfigure}[t]{0.49\linewidth}
    \includegraphics[width=0.95\linewidth]{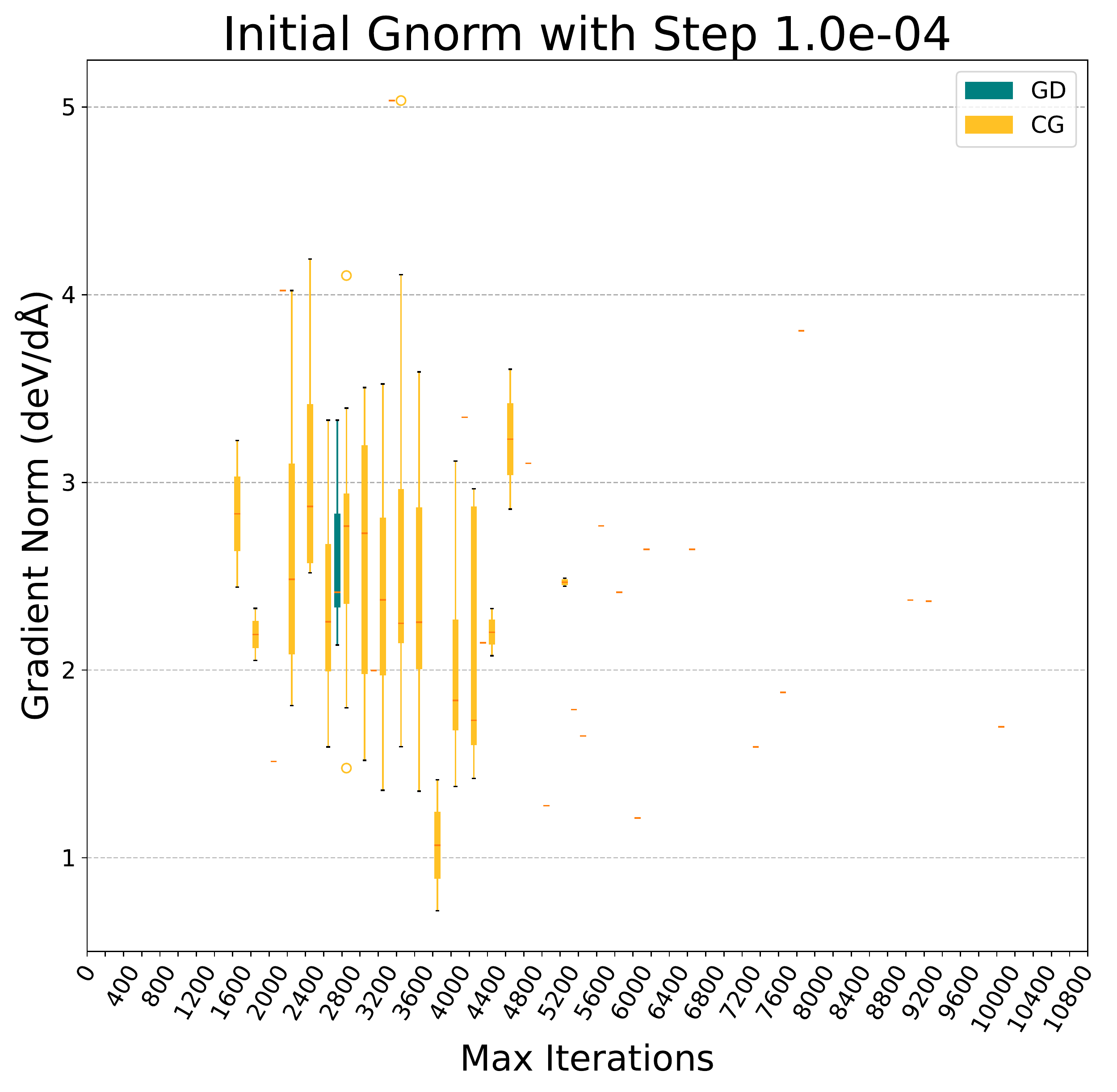}
    \caption{Initial gradient norm versus (rounded to 1000) total iteration number for successful relaxations with constant step $s=\num{0.0001}$.}
    \label{fig:initgnorm1e4}
    \end{subfigure}
    \caption[Range of the initial gradient norm with respect to the number of iterations to success]{\textbf{Range of the initial gradient norm with respect to the number of iterations to success.} Figures (a),(b),(c) and (d) show the distribution of the initial gradient norm among total iteration number. The initial gradient norm is the norm of the gradient of a structure that has not undergone any relaxation yet. Each box matches a range of gradient norm values from the y-axis to a rounded total iterations from the x-axis. The purpose of these plots is to show the relation of the initial gradient norm with the number of total iterations that the successful structures underwent. Green boxes correspond to Gradient Descent and red boxes correspond to Conjugate Gradient.}
    \label{fig:initgnorm}
\end{figure}

\begin{figure}
    \begin{subfigure}[t]{0.49\linewidth}
    \includegraphics[width=0.95\linewidth]{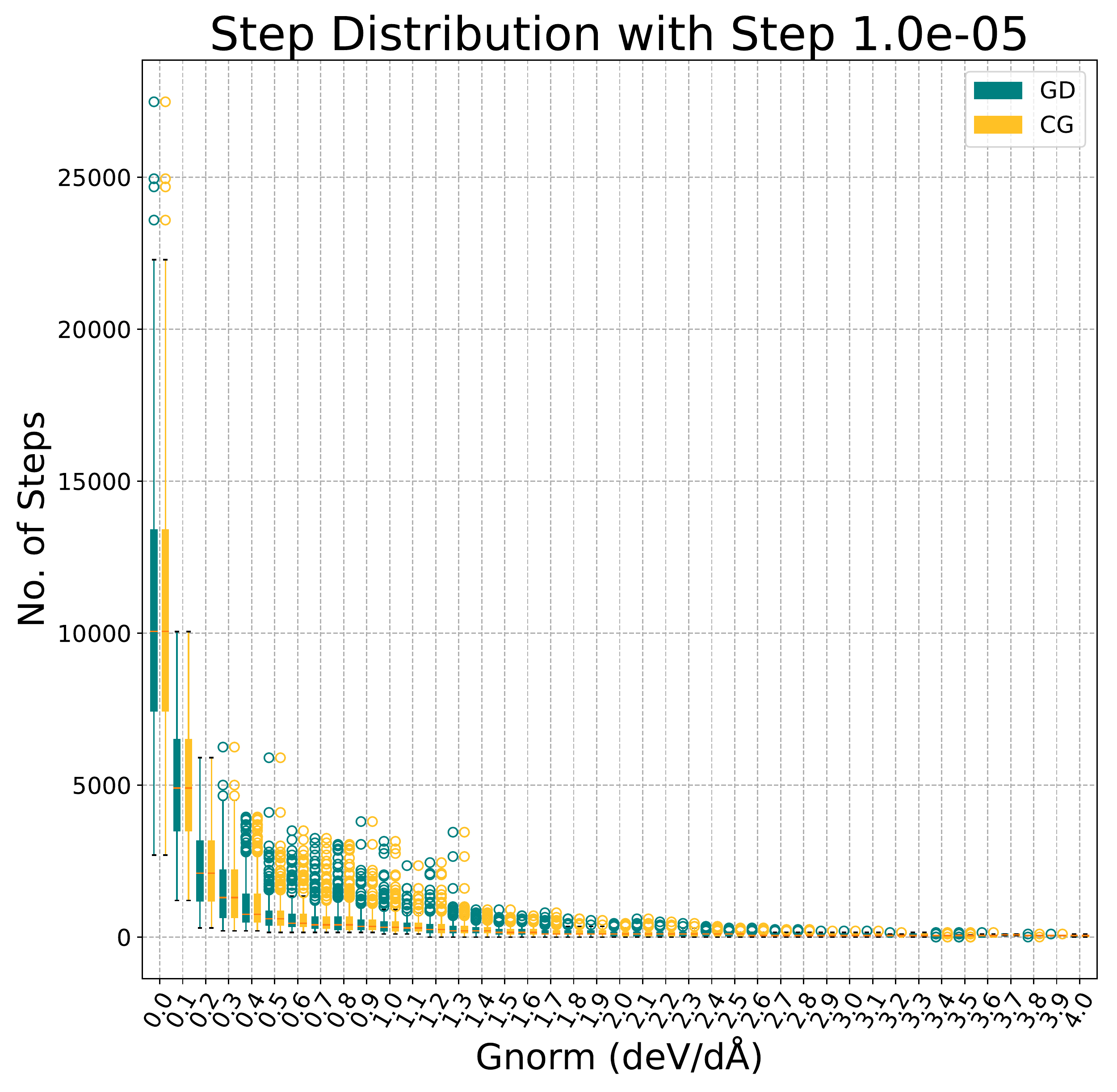}
    \caption{Distribution of iterations to reach values of the gradient norm rounded to the first decimal for experiments run with constant step size $s=\num{0.00001}$. }
    \label{fig:stepdistro1e5}
    \end{subfigure}
    \begin{subfigure}[t]{0.49\linewidth}
    \includegraphics[width=0.95\linewidth]{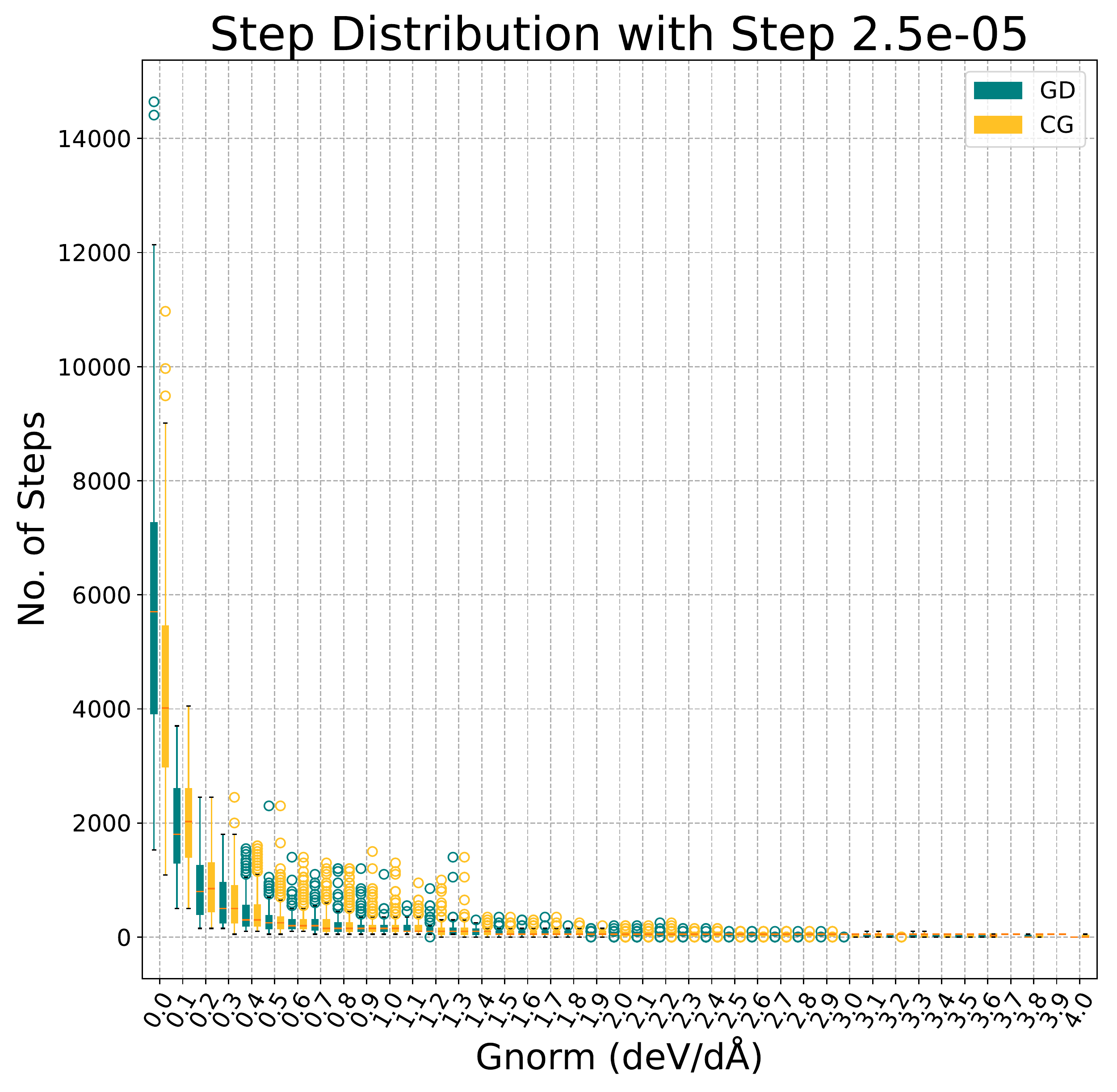}
    \caption{Distribution of iterations to reach values of the gradient norm rounded to the first decimal for experiments run with constant step size $s=\num{0.000025}$.}
    \label{fig:stepdistro2e5}
    \end{subfigure}
    \medskip
    \begin{subfigure}[t]{0.49\linewidth}
    \includegraphics[width=0.95\linewidth]{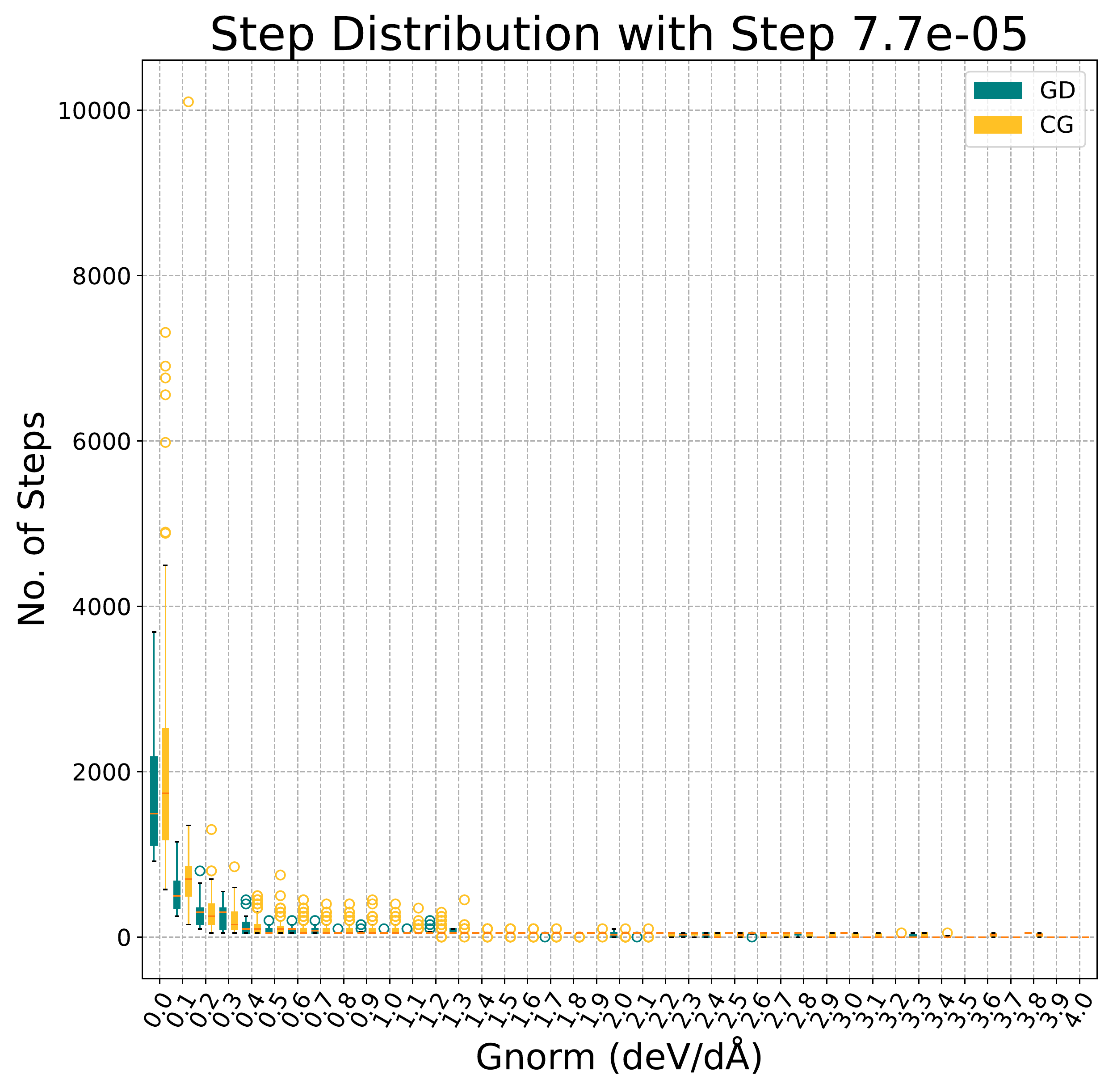}
    \caption{Distribution of iterations to reach values of the gradient norm rounded to the first decimal for experiments run with constant step size $s=\num{0.0000775}$. }
    \label{stepdistro7e5}
    \end{subfigure}
     \hspace{0.2cm}
    \begin{subfigure}[t]{0.49\linewidth}
    \includegraphics[width=0.95\linewidth]{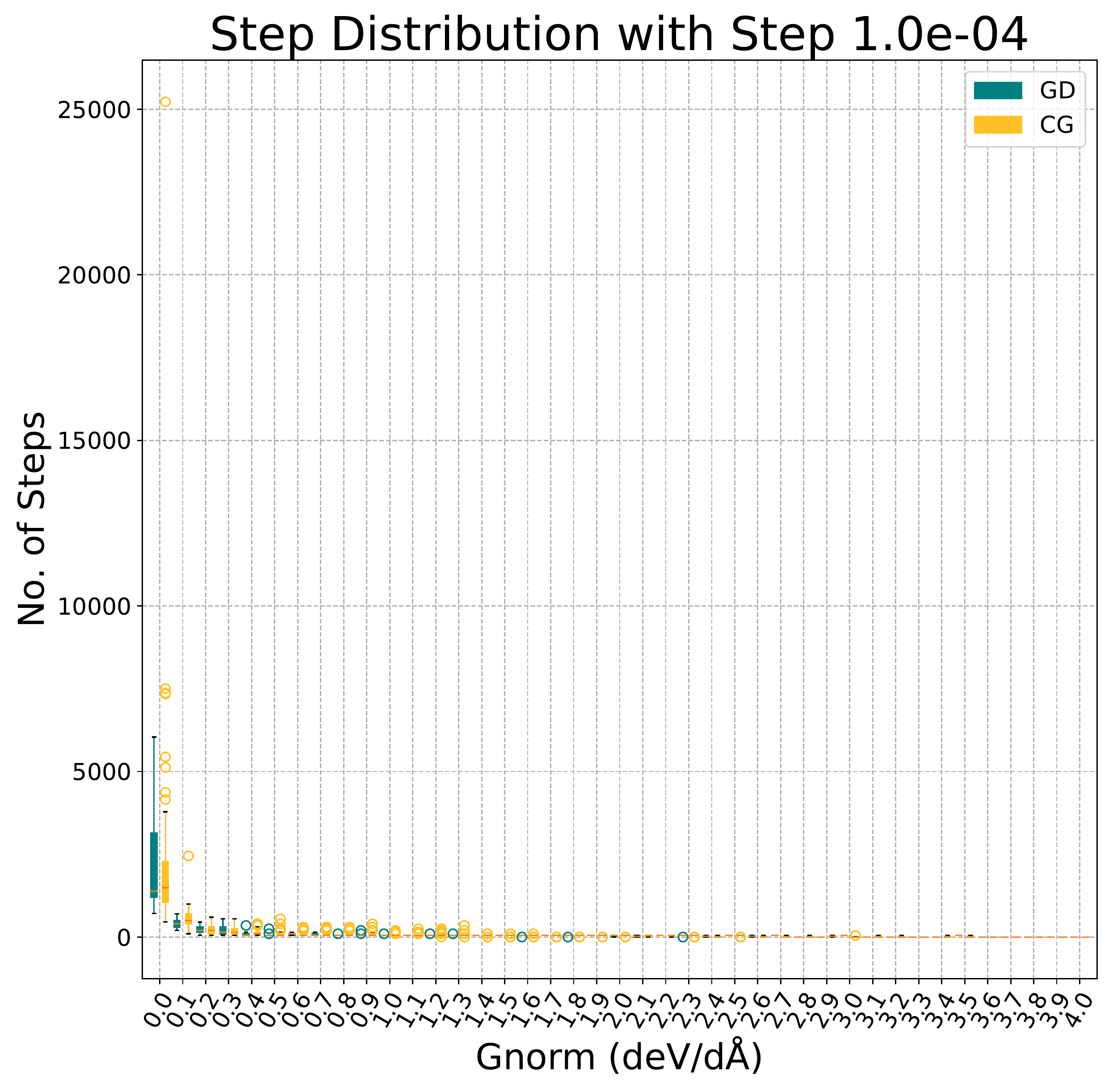}
    \caption{Distribution of iterations to reach values of the gradient norm rounded to the first decimal for experiments run with constant step size $s=\num{0.0001}$.}
    \label{stepdistro1e4}
    \end{subfigure}
    \caption[Distribution of relaxation iterations with respect to gradient norm values]{\textbf{Distribution of relaxation iterations with respect to gradient norm values.} In the above plots the y-axis corresponds to number of steps/iterations and the x-axis corresponds to values of gradient norm. Each box shows the range of number of iterations that was needed so that the norm would decrease by 0.1 deV/dÅ. The ranges include the relaxations that were successfully completed using Gradient Descent (green) and Conjugate Gradient (red) with an input of 200 instances.}
    \label{fig:stepdistro}
\end{figure}

\subparagraph*{Initial gradient norm vs total iterations}
Unfortunately, the gradient norm of the initial configuration does not provide any indication for the total iteration number of successful relaxations. Figure~\ref{fig:initgnorm} illustrates this argument by associating the range of initial gradient norm values among the 200 structures with the number of steps to success. We observe that the norm varies with undefined probability in relation to the number of iterations a structure goes through until successful completion. Hence, no predictions can be made with respect to the running time using this information. This is because, when the optimisation starts, we notice a large drop both in energy and gradient values, that eventually stabilises. Figure~\ref{fig:stepdistro}, which depicts the range of number of steps taken to go from a gradient norm value to the immediately smaller, it is apparent that the majority of iterations is realised for gradient norm values smaller than 0.3. This strengthens our argument that in the first few iterations the decrease in gradient and energy is fast and exponential to a degree that their values are not associated with the final iterations.

\begin{figure}
    \begin{subfigure}[t]{0.49\linewidth}
    \includegraphics[width=0.95\linewidth]{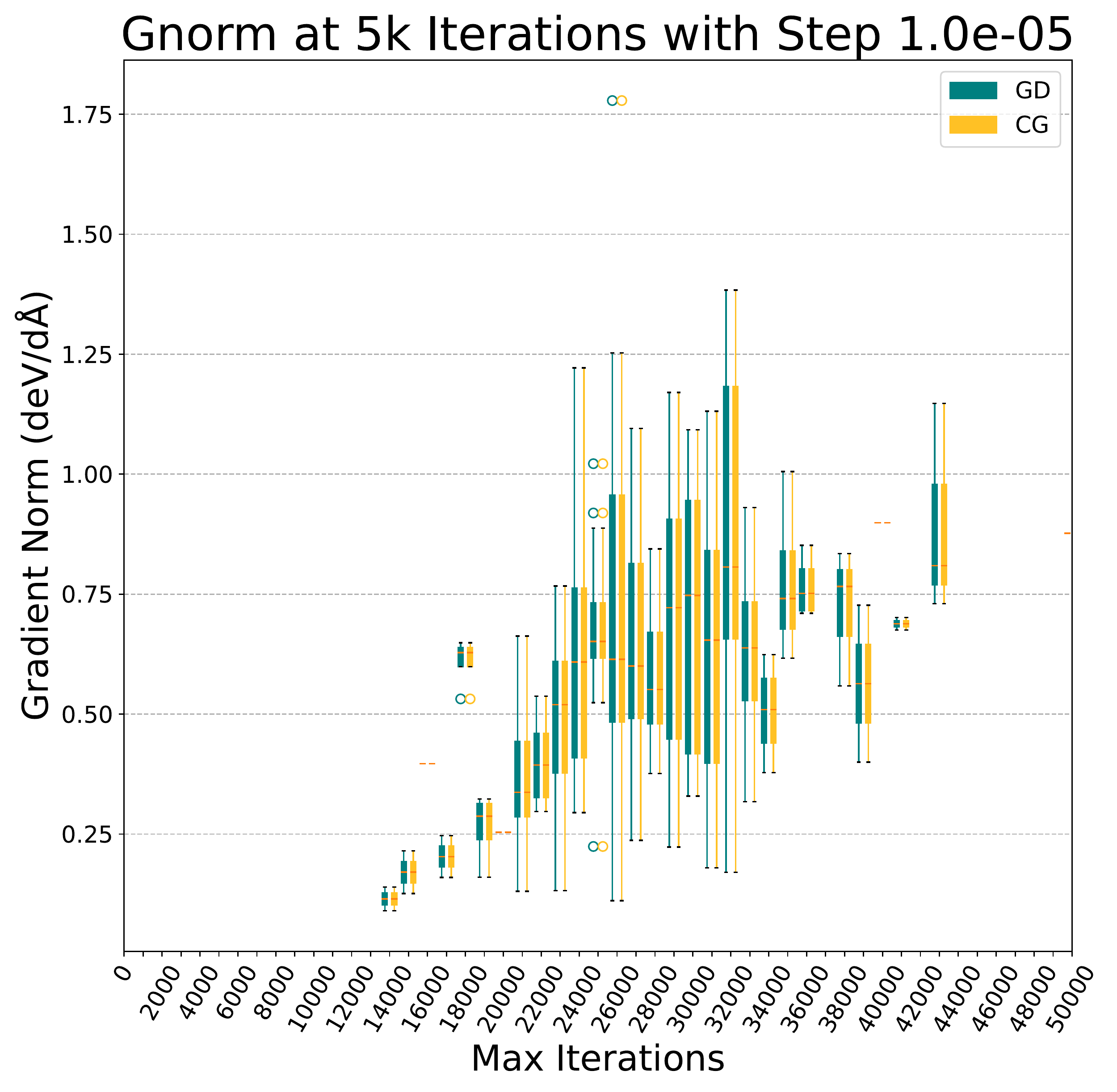}
    \caption{Gradient norm value range at iteration 5000 when the step size is a constant $s=\num{0.00001}$. The x-axis shows the iterations' number to success rounded to 1000.}
    \label{fig:ingnorm1e5}
    \end{subfigure}
    \begin{subfigure}[t]{0.49\linewidth}
    \includegraphics[width=0.95\linewidth]{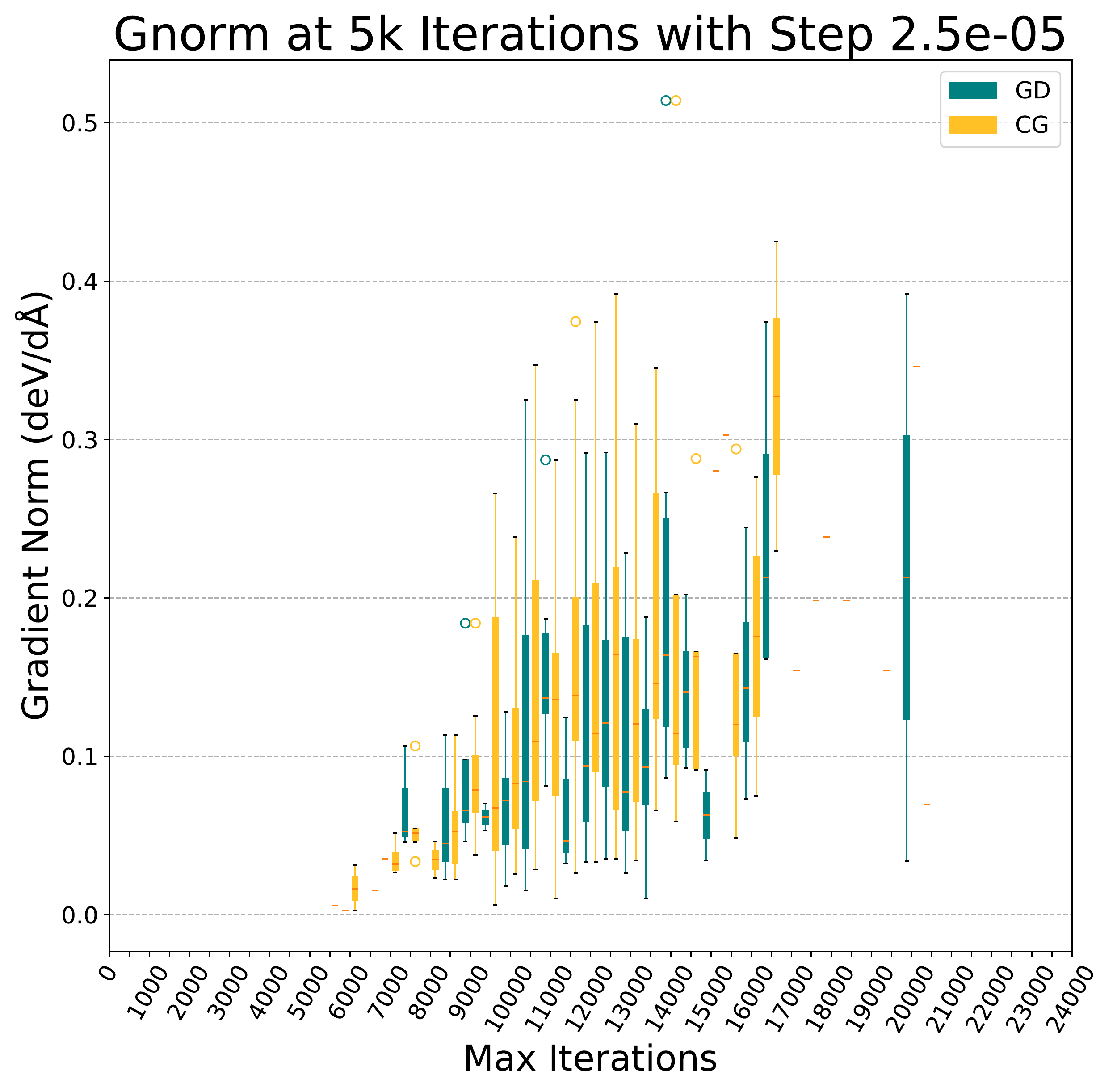}
    \caption{Gradient norm value range at iteration 5000 when the step size is a constant $s=\num{0.000025}$. The x-axis shows the iterations' number to success rounded to 500.}
    \label{fig:ingnorm2e5}
    \end{subfigure}
    \medskip
    \begin{subfigure}[t]{0.49\linewidth}
    \includegraphics[width=0.95\linewidth]{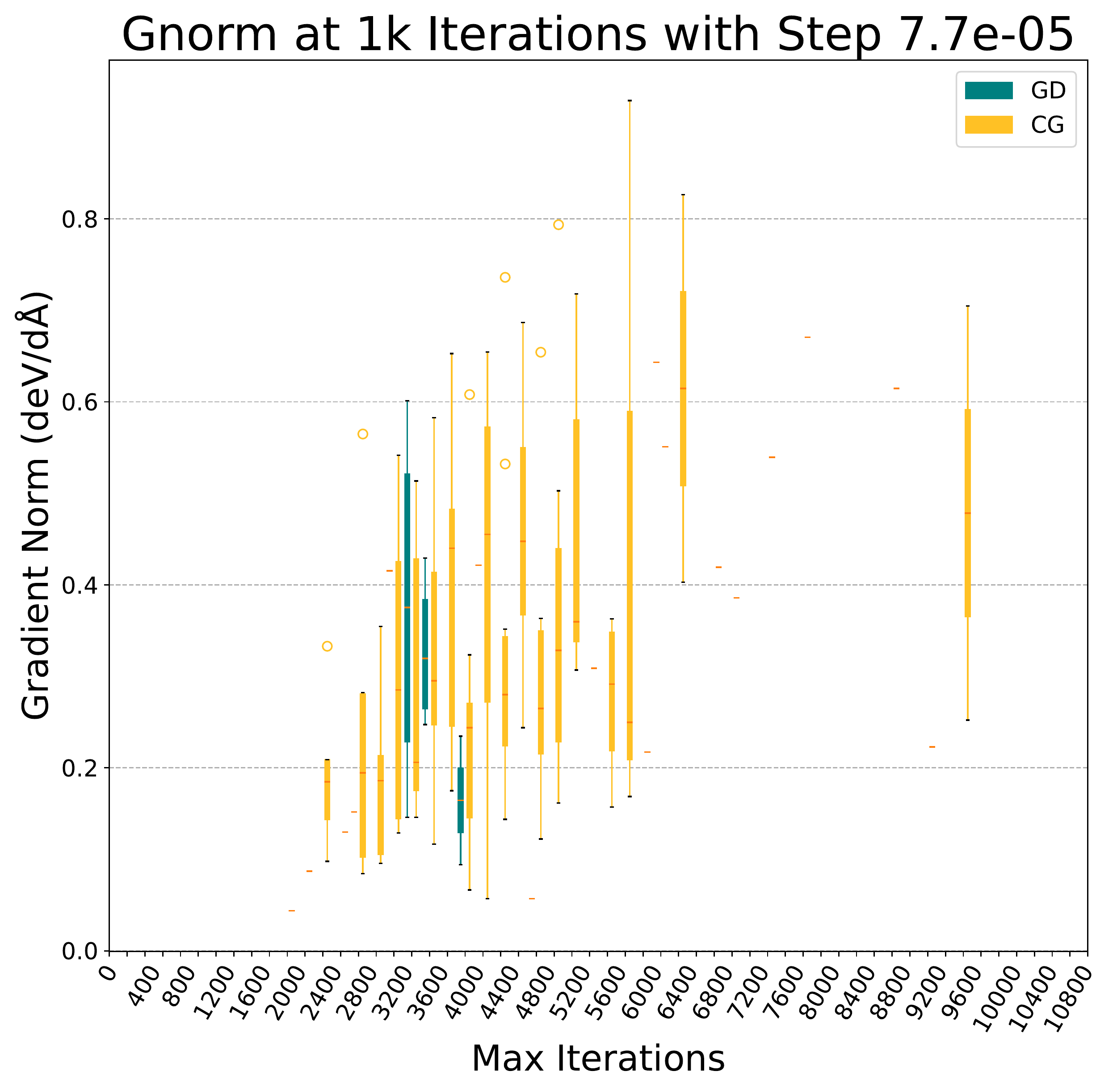}
    \caption{Gradient norm value range at iteration 1000 when the step size is a constant $s=\num{0.0000775}$. The x-axis shows the iterations' number to success rounded to 200.}
    \label{fig:ingnorm7e5}
    \end{subfigure}
    \hspace{0.2cm}
    \begin{subfigure}[t]{0.49\linewidth}
    \includegraphics[width=0.95\linewidth]{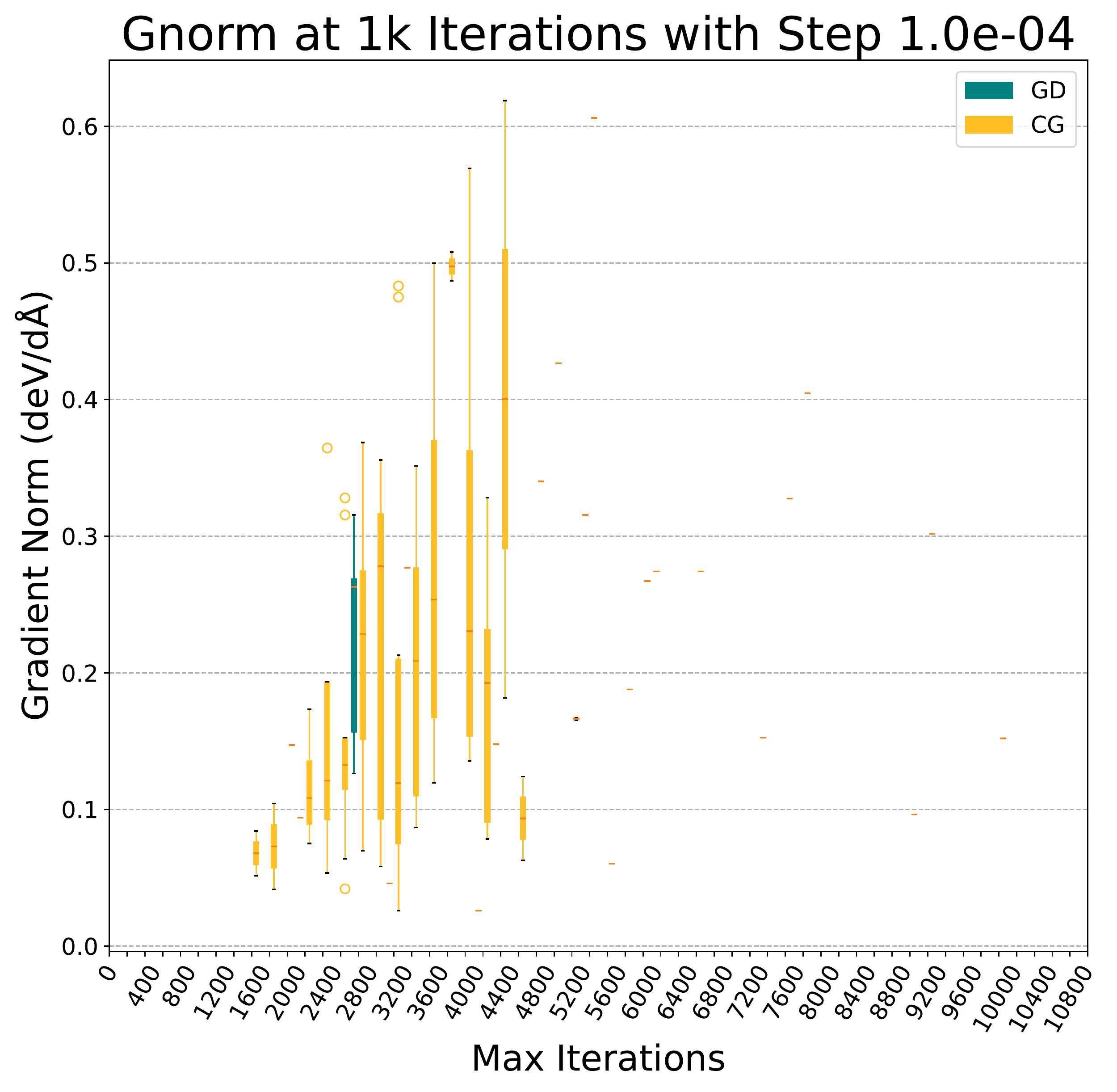}
    \caption{Gradient norm value range at iteration 1000 when the step size is a constant The x-axis shows the iterations' number to success rounded to 500.}
    \label{fig:ingnorm1e4}
    \end{subfigure}
    \caption[Distribution of the gradient norm values with respect to the number of iterations to success]{\textbf{The range of the gradient norm values with respect to the number of iterations to success.} Each box matches a range of values of the gradient norm of successful experiments from the y-axis to a number of total iterations on the x-axis. The iteration number has been rounded to different levels per constant step size, in accordance with the overall experimentation process length for the particular step size. The methods used on all 200 structures are Gradient Descent (green) and Conjugate Gradient (red).
    }
    \label{fig:ingnorm}
\end{figure}

\subparagraph*{Late gradient norm vs total iterations}
This is not true when we examine the gradient norm of a configuration later in the relaxation process. Considering the norm of the gradient around iteration 5000 for $s = \num{0.00001}$ we can match the increase of its magnitude to the increase in iteration number, according to Figure~\ref{fig:ingnorm}. The almost linear connection of gradient norm and total iteration number becomes even more clear for Conjugate Gradient and step size $s=\num{0.000025}$. This proportionality is to be expected, since a gradient with larger values implies more steps to be taken in order for it to decrease below a certain tolerance value $tol$. This relation cannot be observed when $s = \num{0.0001}$, for which the gradient norm is not proportional to the iteration number increase. For Gradient Descent with $s=\num{0.000025}$, on the other hand, same total iteration numbers as when $s=\num{0.00001}$ are associated with smaller gradient norm values. This means that some experiments studied at the same stage of the relaxation -- at 5000 iterations -- but with different step sizes, $s=\num{0.00001}$ and $s=\num{0.000025}$,  had reached different energy values with the latter being closer to the minimum. Even so, the relaxation procedure lasted for the same number of iterations, thus more steps were taken close to the minimum for the second case. 

\begin{figure}
    \begin{subfigure}[c]{0.49\linewidth}
    \includegraphics[width=0.95\linewidth]{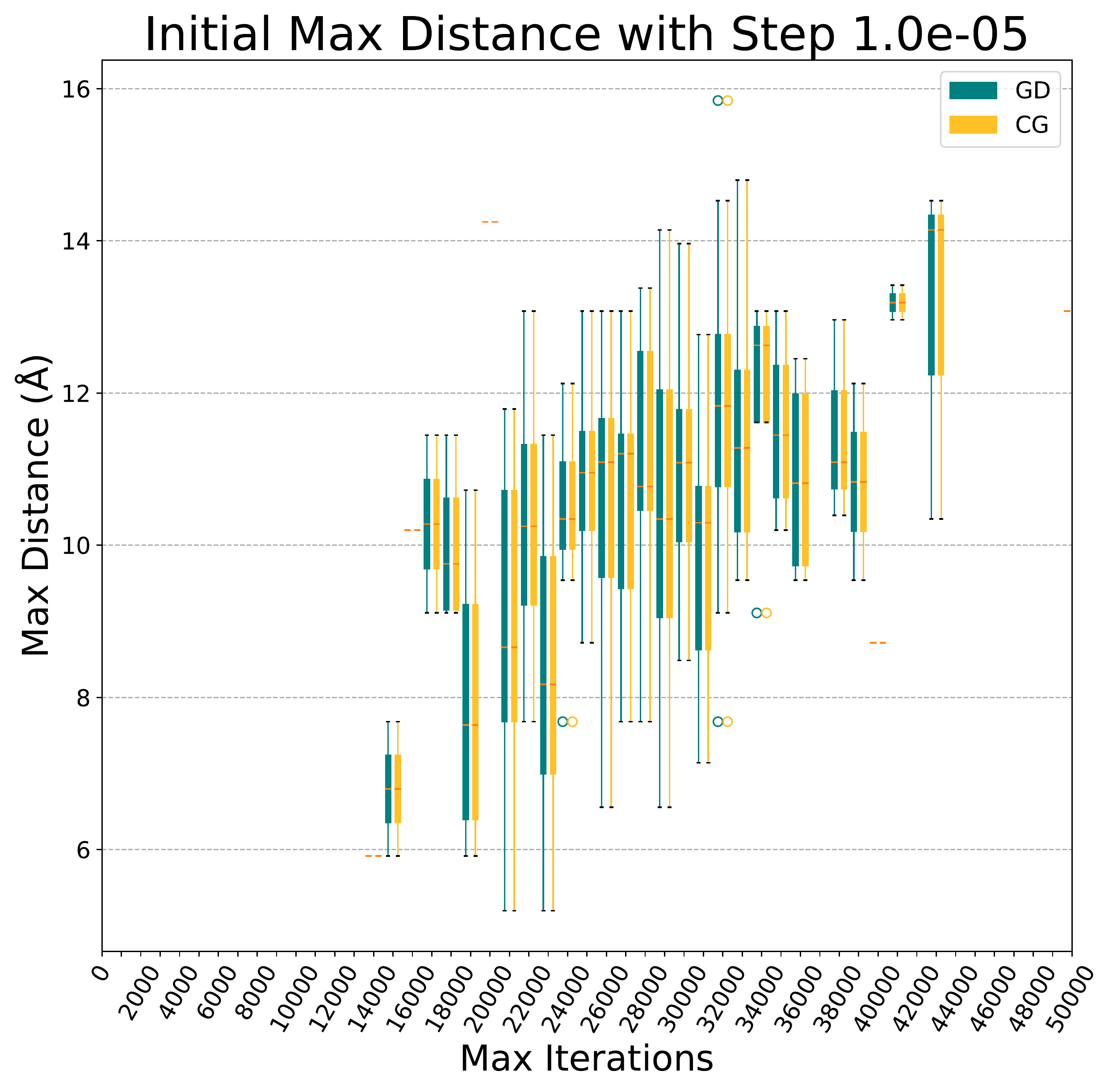}
    \caption{Max pairwise distances in the initial unit cell with respect to the  number of iterations to success rounded to 1000. The step size used is a constant $s=\num{0.00001}$.}
    \label{fig:initdist1e5}
    \end{subfigure}
    \begin{subfigure}[c]{0.49\linewidth}
    \includegraphics[width=0.95\linewidth]{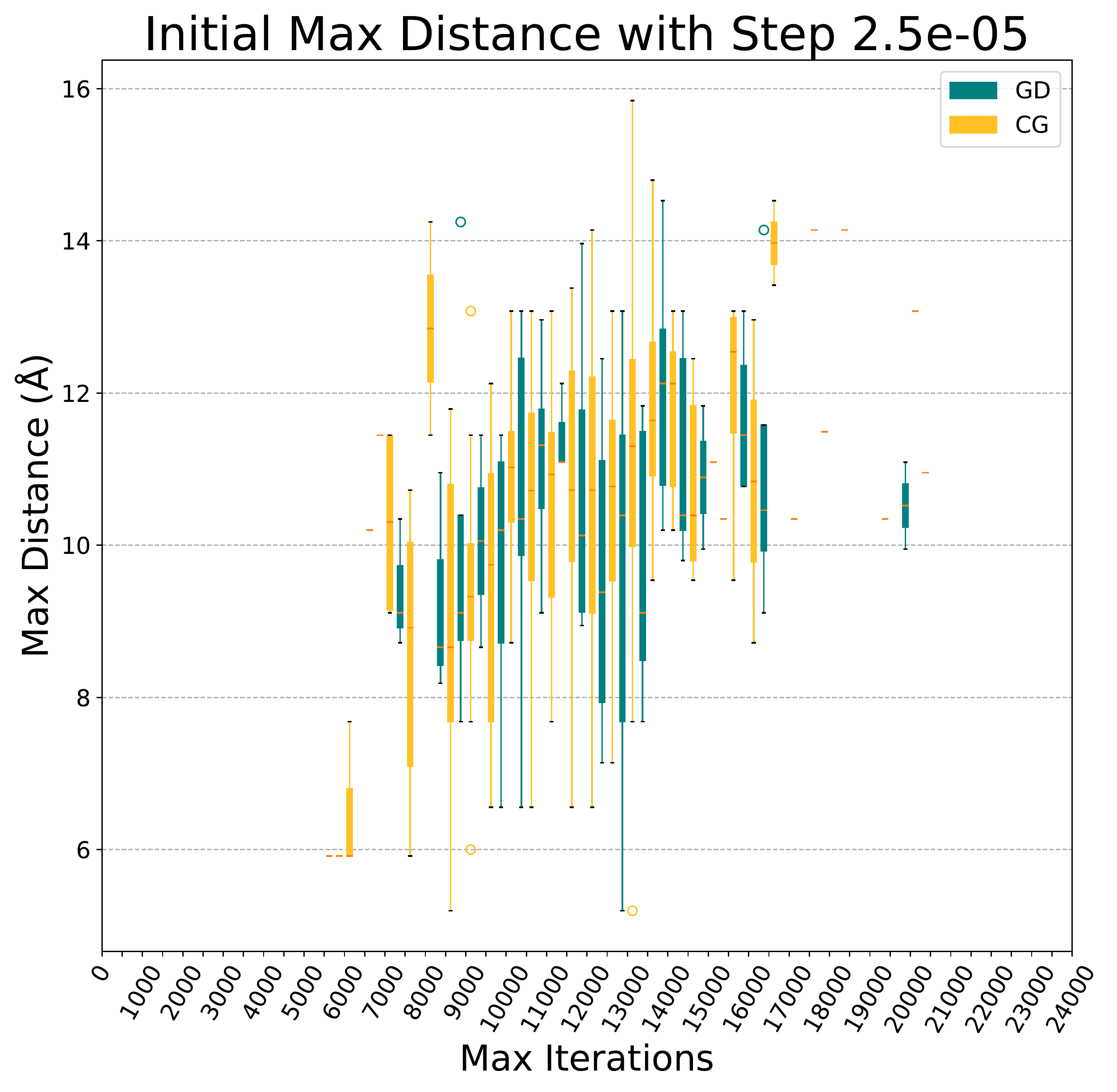}
    \caption{Max pairwise distances in the initial unit cell with respect to the  number of iterations to success rounded to 500.  The step size used is a constant $s=\num{0.000025}$.}
    \label{fig:initdist2e5}
    \end{subfigure}
    \medskip
    \begin{subfigure}[c]{0.49\linewidth}
    \includegraphics[width=0.95\linewidth]{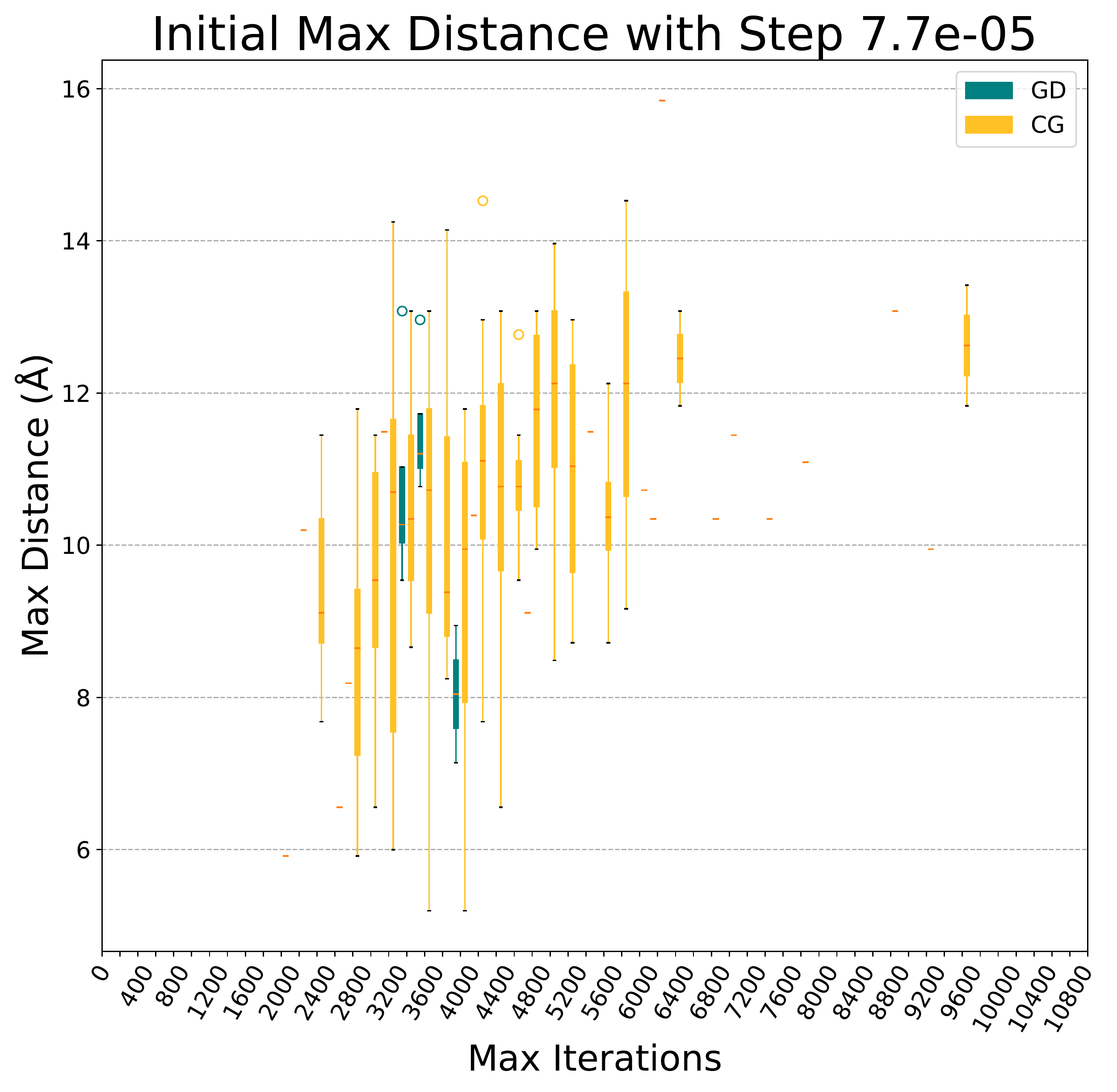}
    \caption{Max pairwise distances in the initial unit cell with respect to the  number of iterations to success rounded to 200. The step size used is a constant $s=\num{0.0000775}$.}
    \label{fig:initdist7e5}
    \end{subfigure}
    \hspace{0.2cm}
    \begin{subfigure}[c]{0.49\linewidth}
    \includegraphics[width=0.95\linewidth]{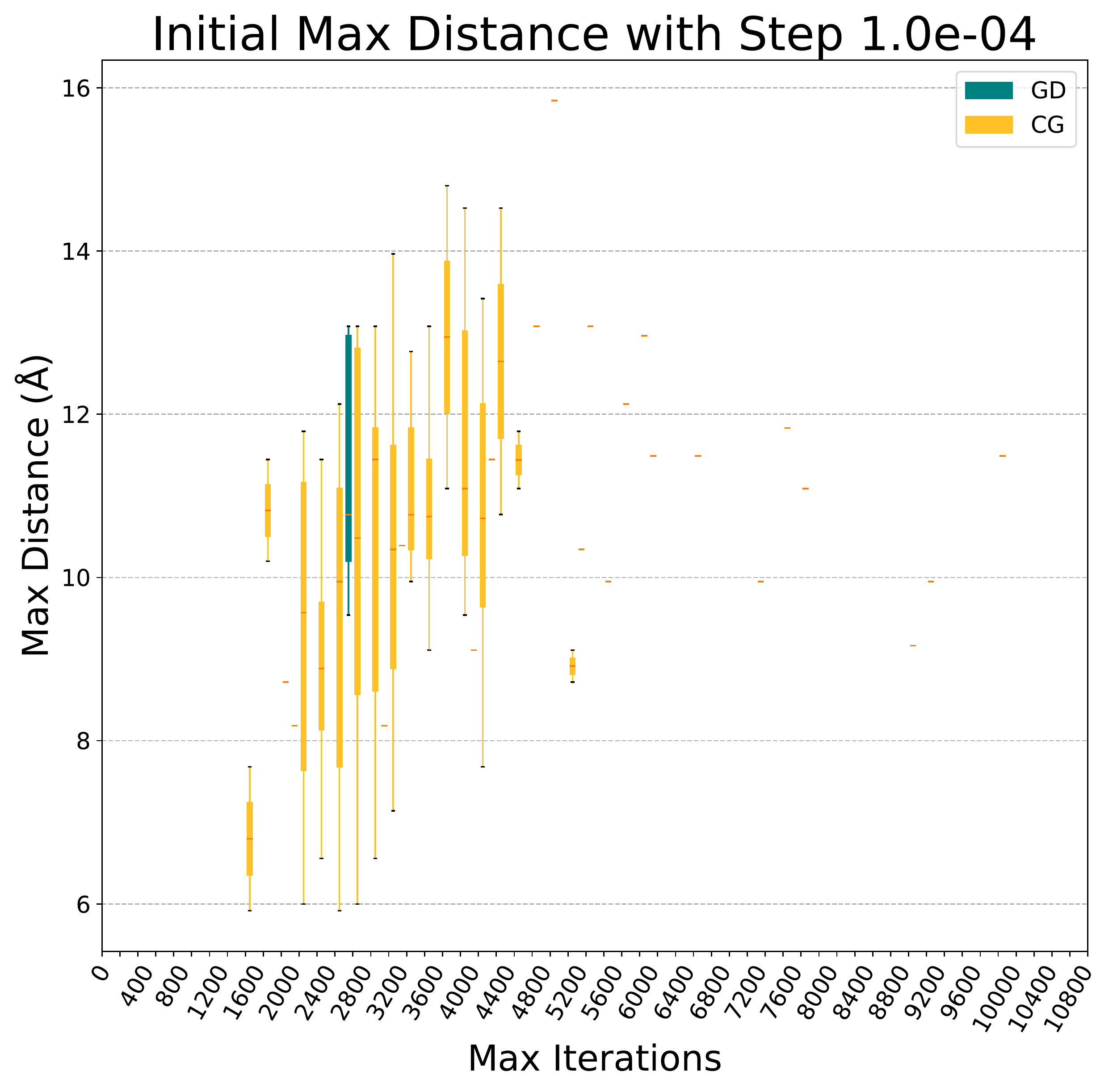}
    \caption{Max pairwise distances in the initial unit cell with respect to the  number of iterations to success rounded to 500. The step size used is a constant $s=\num{0.0001}$.}
    \label{fig:initdist1e4}
    \end{subfigure}
    \caption[Distribution of the maximum pairwise distances in the unit cell with respect to the number of steps to success]{\textbf{Range of the maximum pairwise distances in the unit cell with respect to the number of steps to success.} The distances in question are the pairwise distances of the ions before the relaxation starts. The iteration number has been rounded to different levels per constant step size, in accordance with the overall experimentation process length for the particular step size. The methods used on all 200 structures are Gradient Descent (green) and Conjugate Gradient (red).}
    \label{fig:initmaxdist}
\end{figure}

\subparagraph*{Initial max pairwise distance vs total iterations}
In contrast to the initial gradient norm, we can predict the number of iterations to success using the maximum pairwise distance of ions in the initial configuration. We anticipate that longer pairwise distances in the cell imply more optimisation steps, and this is confirmed by the following. The range of max pairwise distances in the initial unit cell versus iterations to success is depicted in Figure~\ref{fig:initmaxdist}. Pairwise distances in the initial unit cell require accordingly long relaxations to arrive to completion when $s = \num{0.00001}$ for either Gradient Descent or Conjugate Gradient (Figure~\ref{fig:initdist1e5}). This is also the case for $s = \num{0.000025}$ and Conjugate Gradient, as seen in Figure~\ref{fig:initdist2e5}, yet Gradient Descent starts to show different behaviour with this step size for structures that finish in more than 15.000 iterations. We observe that for the same initial pairwise distances the number of iterations to success has increased for some cases of $s=\num{0.000025}$ compared to $s=\num{0.00001}$. This happens because Gradient Descent produces large updating steps that cannot lead directly to the minimiser and more iterations are needed to redirect the procedure back to it. However, Conjugate Gradient maintains a seemingly linear relation between total number of iterations and maximum initial pairwise distance for step size values up to $s=\num{0.0001}$. When $s\geq \num{0.0001}$ the step size is large enough to break this pattern and we can thus observe that structures with small interatomic distances took longer to converge compared to these with larger pairwise distances.

\begin{figure}
    \centering
    \begin{subfigure}[c]{0.49\linewidth}
        \includegraphics[width=0.95\linewidth]{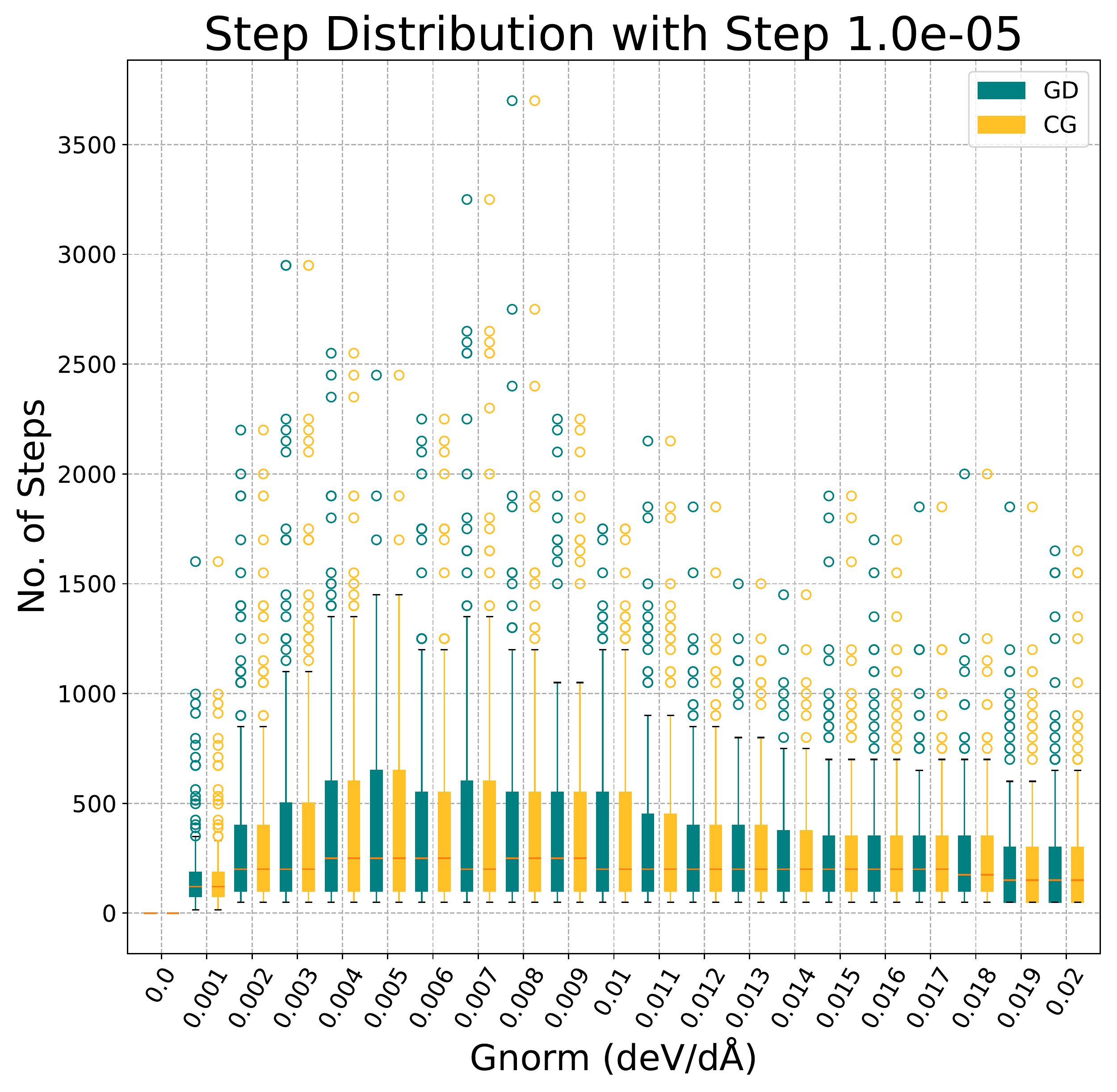}
        \caption{Distribution of steps with respect to a 0.001 Å gradient norm decrease. The constant step size for this plot was $s=\num{0.00001}$.}
        \label{fig:stepdistrozoom1e5}
    \end{subfigure}
    \begin{subfigure}[c]{0.49\linewidth}
        \includegraphics[width=0.95\linewidth]{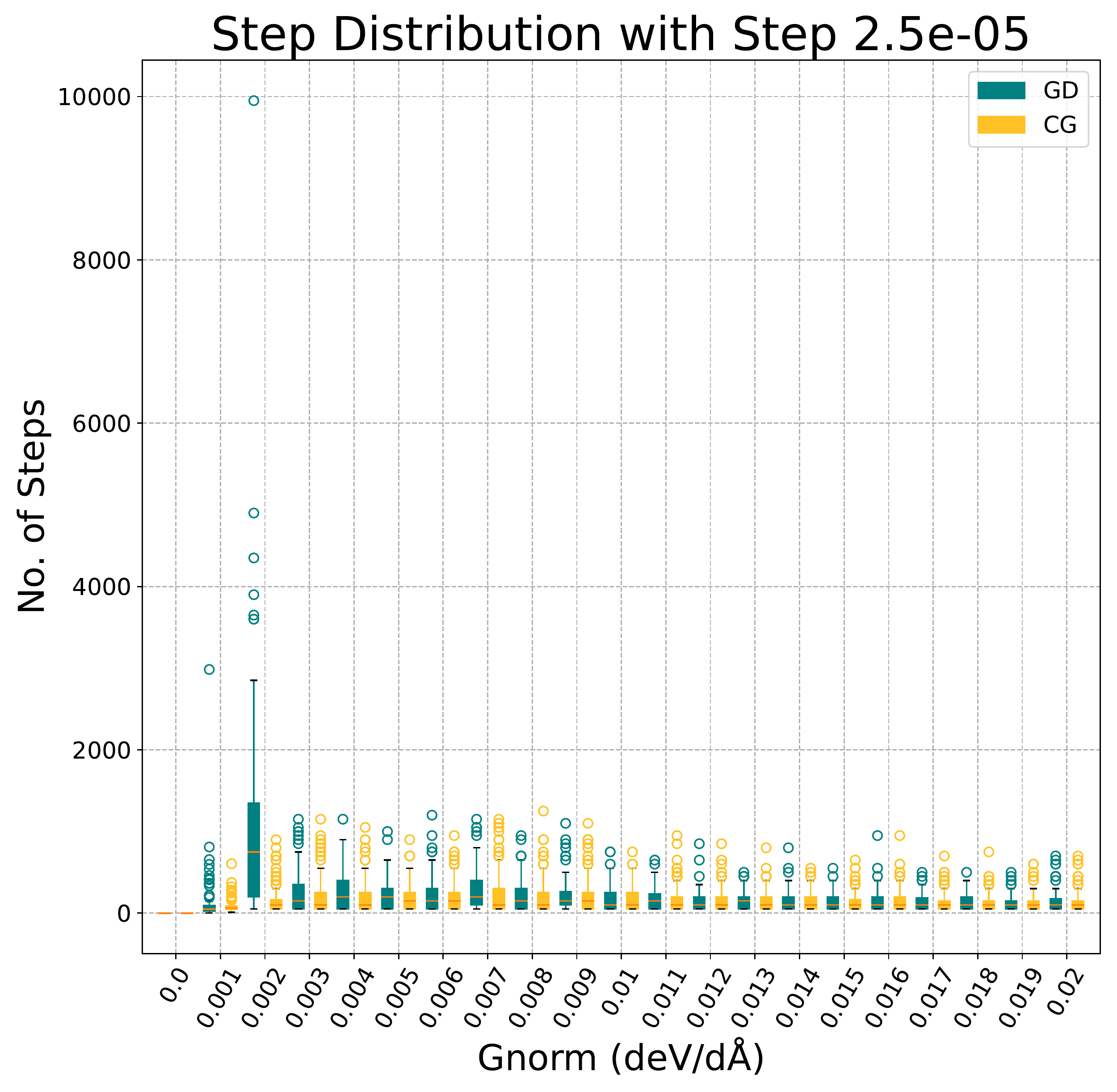}
        \caption{Distribution of steps with respect to a 0.001 Å gradient norm decrease. The constant step size for this plot was $s=\num{0.000025}$.}
        \label{fig:stepdistrozoom2e5}
    \end{subfigure}
    \medskip
    \begin{subfigure}[c]{0.49\linewidth}
        \includegraphics[width=0.95\linewidth]{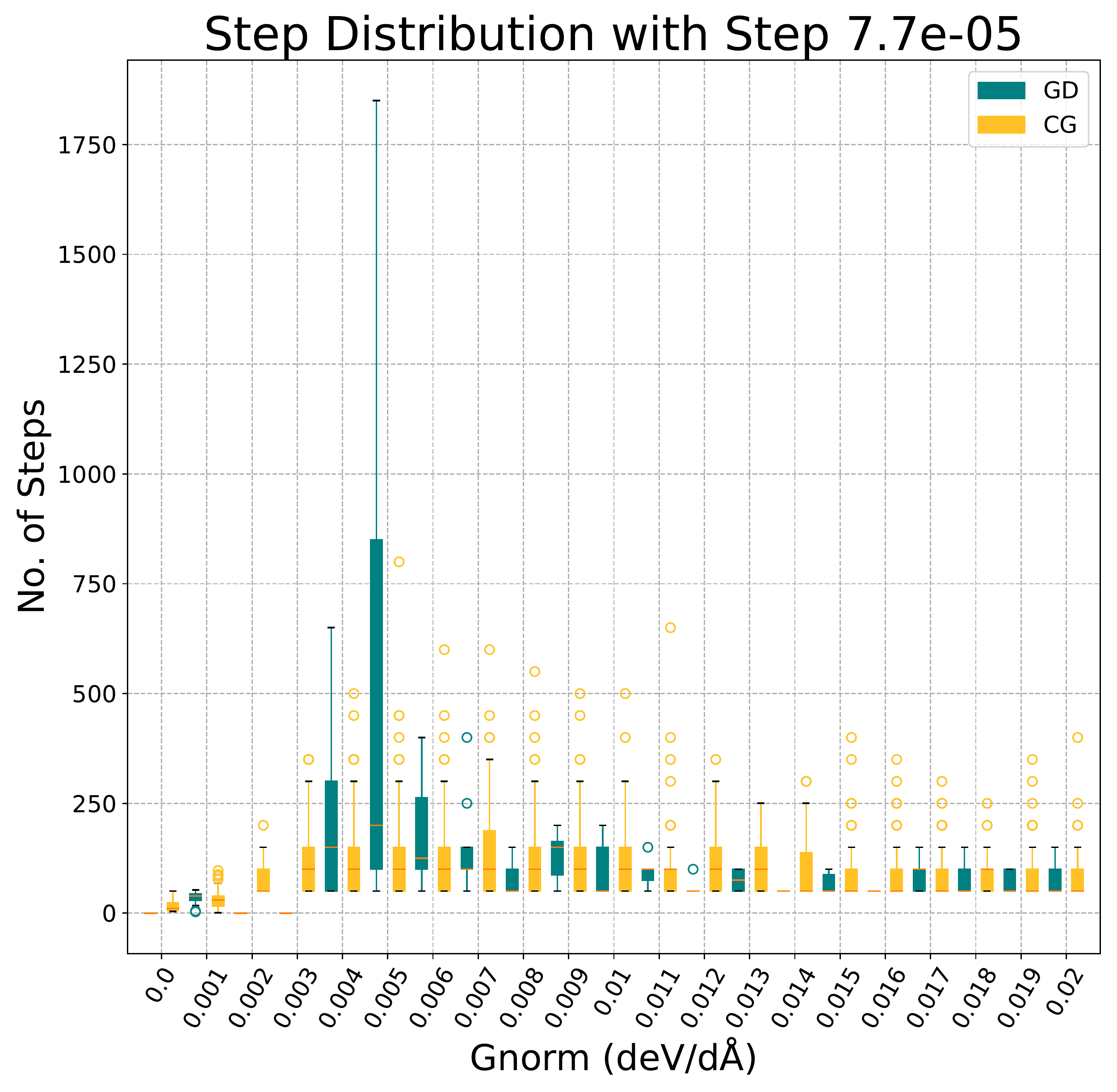}
        \caption{Distribution of steps with respect to a 0.001 Å gradient norm decrease. The constant step size for this plot was $s=\num{0.0000775}$.}
        \label{fig:stepdistrozoom7e5}
    \end{subfigure}
    \begin{subfigure}[c]{0.49\linewidth}
        \includegraphics[width=0.95\linewidth]{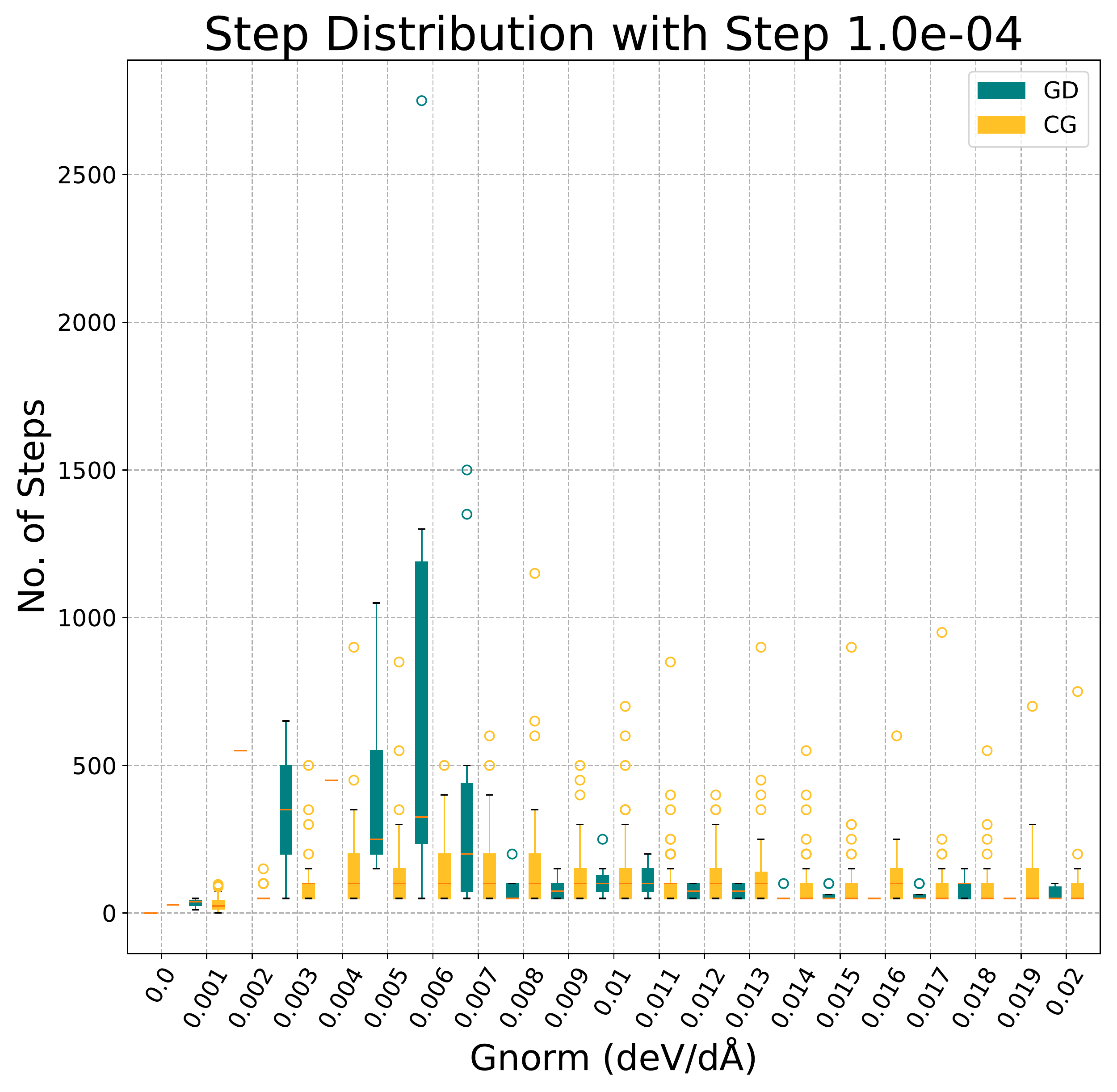}
        \caption{Distribution of steps with respect to a 0.001 Å gradient norm decrease. The constant step size for this plot was $s=\num{0.0001}$.}
        \label{fig:stepdistrozoom1e4}
    \end{subfigure}
    \caption[Distribution of steps with respect to gradient norm close to zero]{\textbf{Step number distribution with respect to gradient norm close to zero.} The boxes in the above figures show the ranges of number of steps that were needed to decrease from one gradient norm value to the smaller by 0.001 deV/dÅ value. These correspond to the successful experiments out of 200 input structures per method for Gradient Descent (green) and Conjugate Gradient (red) with constant step size. Each plot shows the results for one value of constant step size. The results are focused around a small neighbourhood of the minimum where the gradient norm was less than $4\cdot 10^{-2}$ deV/dÅ.}
    \label{fig:stepdistrozoom}
\end{figure}

\subparagraph{}
The previously found difference in robustness between Gradient Descent and Conjugate Gradient is owed to the nature of the two algorithms. It is apparent that the first resorts to many iterations of redirection until the path to the minimiser is retrieved, even for as small step sizes as $s=\num{0.000025}$. On the other hand, Conjugate Gradient succeeds to steadily approach the minimiser with gradual progress even for $s=\num{0.0001}$. As already mentioned, most of the iterations happen when the gradient norm is below 0.3 and Figure~\ref{fig:stepdistrozoom} reveals that, for $s>\num{0.00001}$, a vast number of Gradient Descent iterations is taken towards the end. Nonetheless, with Conjugate Gradient we observe an evenly distributed number of iterations. As a consequence, we can conclude that Conjugate Gradient steadily follows a smooth path constructed to lead towards the minimum, whereas Gradient Descent follows small steps that need constant readjustment in order to eventually point to the correct direction. 

To summarise, we have seen that Conjugate Gradient maintains a success rate of at least 50\% in 200 experiments for up to 10 times the smallest constant step size value $s=\num{0.00001}$. On the other hand, Gradient Descent's success rate falls by roughly 50\% when the step size is a little more than doubled from $s=\num{0.00001}$ to $s=\num{0.000025}$. Conjugate Gradient can perform well with a wider range of constant step sizes revealing an increased utility both in terms of speed and success rate. This is due to its ability to adapt well the direction of relaxation and make steady progress. Ultimately, its versatility can accelerate a process of multiple experiments without having to resort to failed experiments. We observe, nonetheless, that there is a threshold ($s=\num{0.0001}$) above which step size increase affects the convergence speed of Conjugate Gradient negatively, so one cannot increase the step size value further and expect quicker convergence. For our energy model, a large uneducated increase to the step size can lead to Buckingham catastrophe, which sentences the optimisation procedure to failure. Hence, we conclude that a small step size value that can slowly and steadily lead to the minimum using an update with only the first derivatives, but a larger step size in combination with Conjugate Gradient can accelerate the process.

\subsection{Utility comparison}
\label{experimentsII}

We now introduce a utility function to evaluate and compare the performance of Gradient Descent and Conjugate Gradient with different step size arrangements according to Graham D. et al~\cite{Graham2022FormalizingDistributions}. We set two preferences: success rate and iteration number. The goal is to maximize the utility function $u_{\text{FP}}\in \mathbb{R}$ defined as follows
\begin{equation}\label{eqn:utility}
    u_{\text{FP}} = (1-\lambda)\cdot\frac{I-i_f}{I} + \lambda\tau
\end{equation}
The function is evaluated on the result that each structure produces and then the mean value is used for the overall result per method. The success rate $\tau$ is defined as the percentage of structures that was successfully relaxed from the batch that this structure belonged to. The iteration number $i_f$ is the the total number of iterations per experiment. Given the iterations' upper bound $I$, we consider $i_f$ to be capped by $I$. By introducing function $p(i_f,I)=\frac{I-i_f}{I}$ wherein $I$ is a known captime, we construct $u_{\text{FP}}$ as in Equation~(\ref{eqn:utility}). Following the notation of Graham D. et al, we define the two constants $c_0, c_1$ as $0 < c_1=1-\lambda < 1$ and $c_0 = \lambda\tau$. The $\lambda$ parameter designates the side -- success rate or speed in iterations -- to which we place the most preference. According to this preference, we can select the algorithmic recipe that would mostly correspond to our needs.

\subsubsection{Constant step size utility}
\begin{figure}
    \centering
    \begin{subfigure}[c]{0.49\linewidth}
        \includegraphics[width=\textwidth]{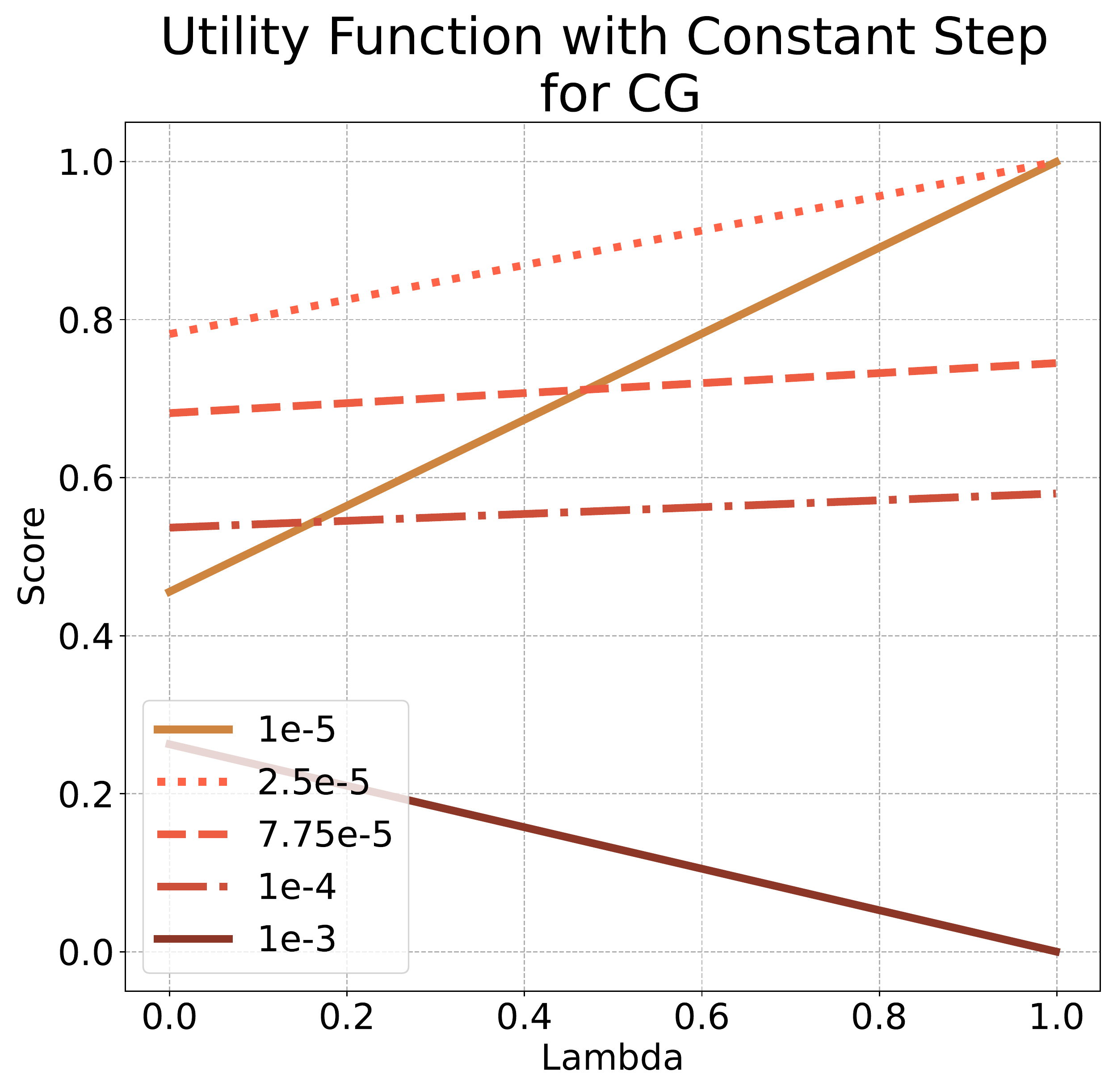}
    \end{subfigure}
    \begin{subfigure}[c]{0.49\linewidth}
        \includegraphics[width=\textwidth]{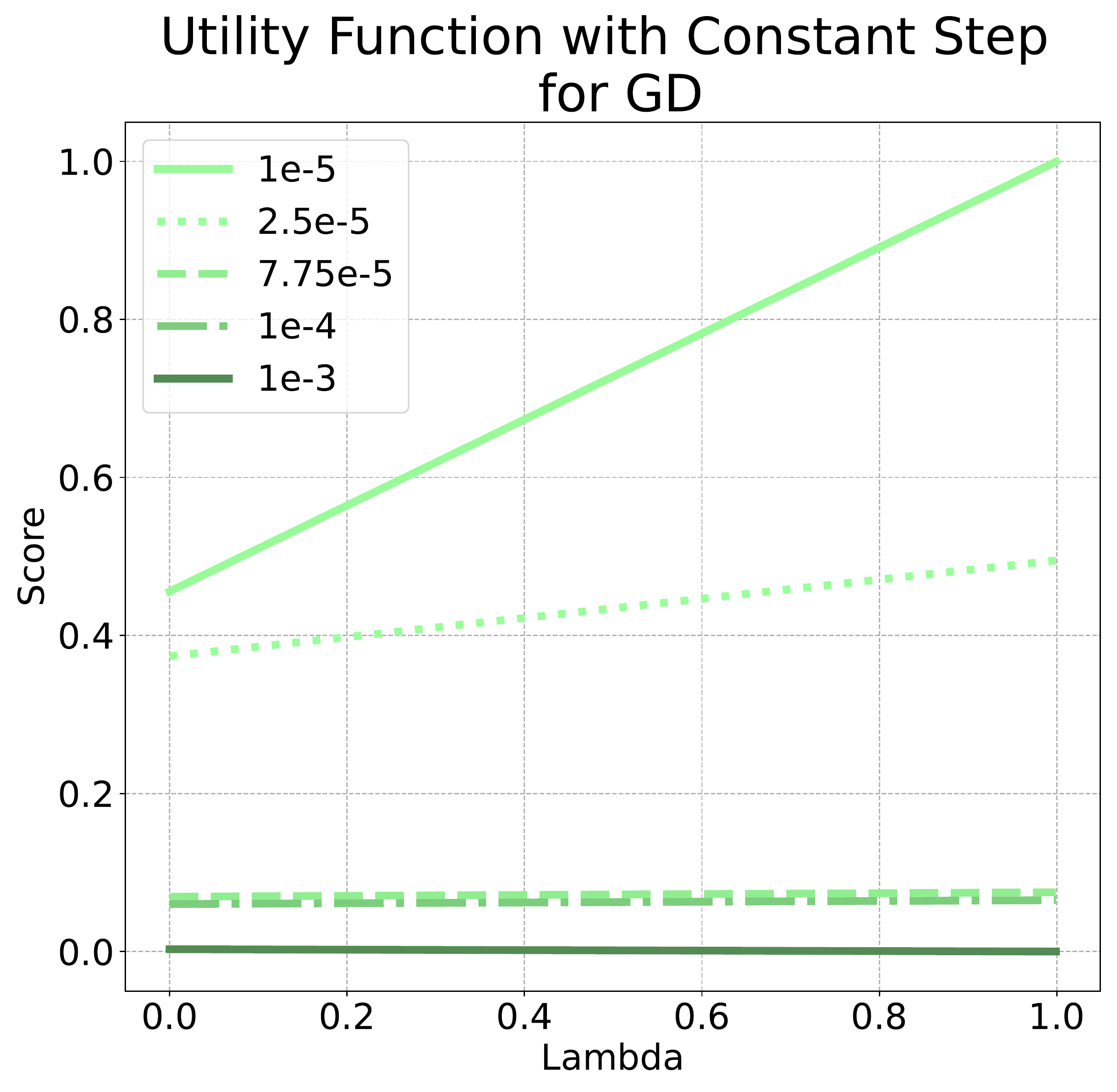}
    \end{subfigure}
    \caption[Utility scores of Gradient Descent and Conjugate Gradient with constant step]{\textbf{Utility scores of Gradient Descent and Conjugate Gradient with constant step.} Each line represents results over 200 structures. \textit{Left}, Scores of Gradient Descent run with constant steps of set $S$ with respect to the $\lambda$ parameter of the utility function $u_{\text{FP}}$; \textit{right}, scores of Conjugate Gradient run with constant steps of set $\mathcal{S}$~\ref{eq:setS} with respect to the $\lambda$ parameter of the utility function $u_{\text{FP}}$.}
    \label{fig:uconst}
\end{figure}

As depicted in Figure~\ref{fig:uconst}, Gradient Descent and Conjugate Gradient display different levels of utility for different constant step sizes. The utility of Gradient Descent decreases as the step size increases, appointing $s=\num{0.00001}$ as the definitive best choice for this algorithm. While Conjugate Gradient features the highest utility with a constant $s=\num{0.000025}$ for $\lambda\rightarrow 0$, it manages to relax all structures with both $\num{0.00001}$ and $\num{0.000025}$, thus the two compete for the highest utility score when $\lambda\rightarrow 1$. However, $s=\num{0.00001}$ becomes a good choice only for $\lambda \geq 0.6$, when high success rate needs to be ensured. When speed is at least as much important, $s\geq\num{0.0000775}$ is a better match for Conjugate Gradient. Figure~\ref{fig:uconst} confirms once more that Gradient Descent and Conjugate Gradient have similar performance for a small $s=\num{0.00001}$, thus, the aforementioned comparison of Conjugate Gradient with $s=\num{0.000025},\num{0.0000775},\num{0.0001}$ and $s=\num{0.00001}$ can be directly applied to Gradient Descent with $s=\num{0.00001}$.

\subsubsection{Scheduled step size utility}

We will first describe various methods of scheduling step size. Let $\underline{s}=\num{0.00001}$ and $\overline{s}=\num{0.001}$, the smallest and largest step size values from the experiments with constant step size. We set the initial value ($s_0$) of a scheduled step size to be $\overline{s}$, our lower bound to be $\underline{s}$ and test the scheduling rules over groups of 40 structures. The scheduled step size methods that we test are the following:

\begin{itemize}
    \item \textit{Bisection (bisect)}. The first scheduling rule we employed is a simple bisection. The initial value of the step size and a lower bound for it are provided. After initialisation, the step size is updated as the mean of its current value and the lower bound every 100 iterations.
    \item \textit{Gradient norm-Scheduled Bisection (gbisect)}. The step size is initialised with a given upper value. A lower bound is also provided. Once the gradient norm is decreased by some order of magnitude $\beta$, the next step size value becomes the mean of the current step size and the lower bound.
    \item \textit{Exponential Scheduled (expo)}.  With this rule an upper bound and a lower bound for the value of the step size are provided. The step size is initialised with the upper bound and is then multiplied by a fixed constant number $0 < \gamma < 1$ at every subsequent iteration until it reaches the lower bound.
\end{itemize}

\begin{figure}
    \begin{subfigure}[c]{0.49\linewidth}
        \includegraphics[width=\textwidth]{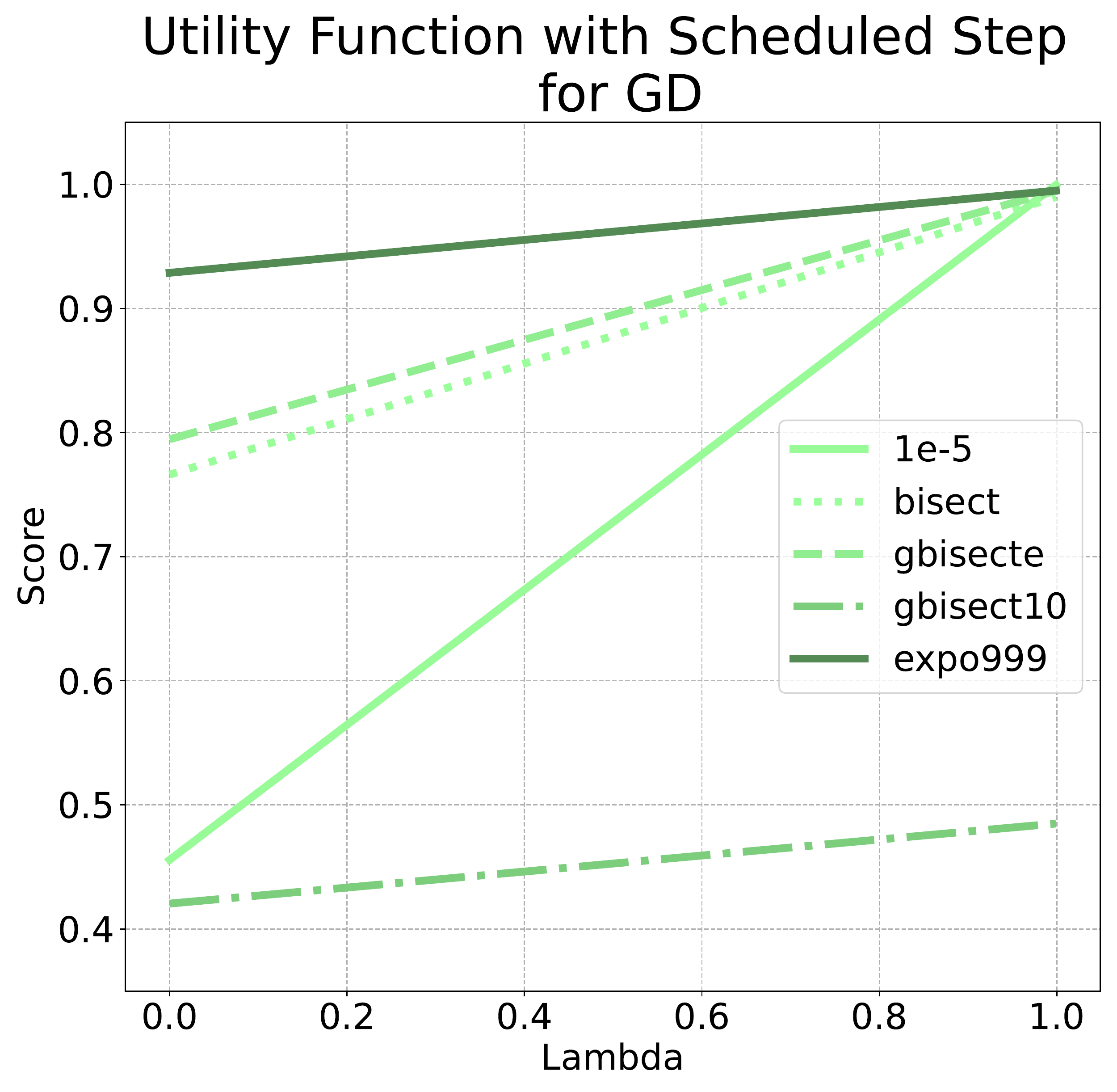}
    \end{subfigure}
    \begin{subfigure}[c]{0.49\linewidth}
        \includegraphics[width=\textwidth]{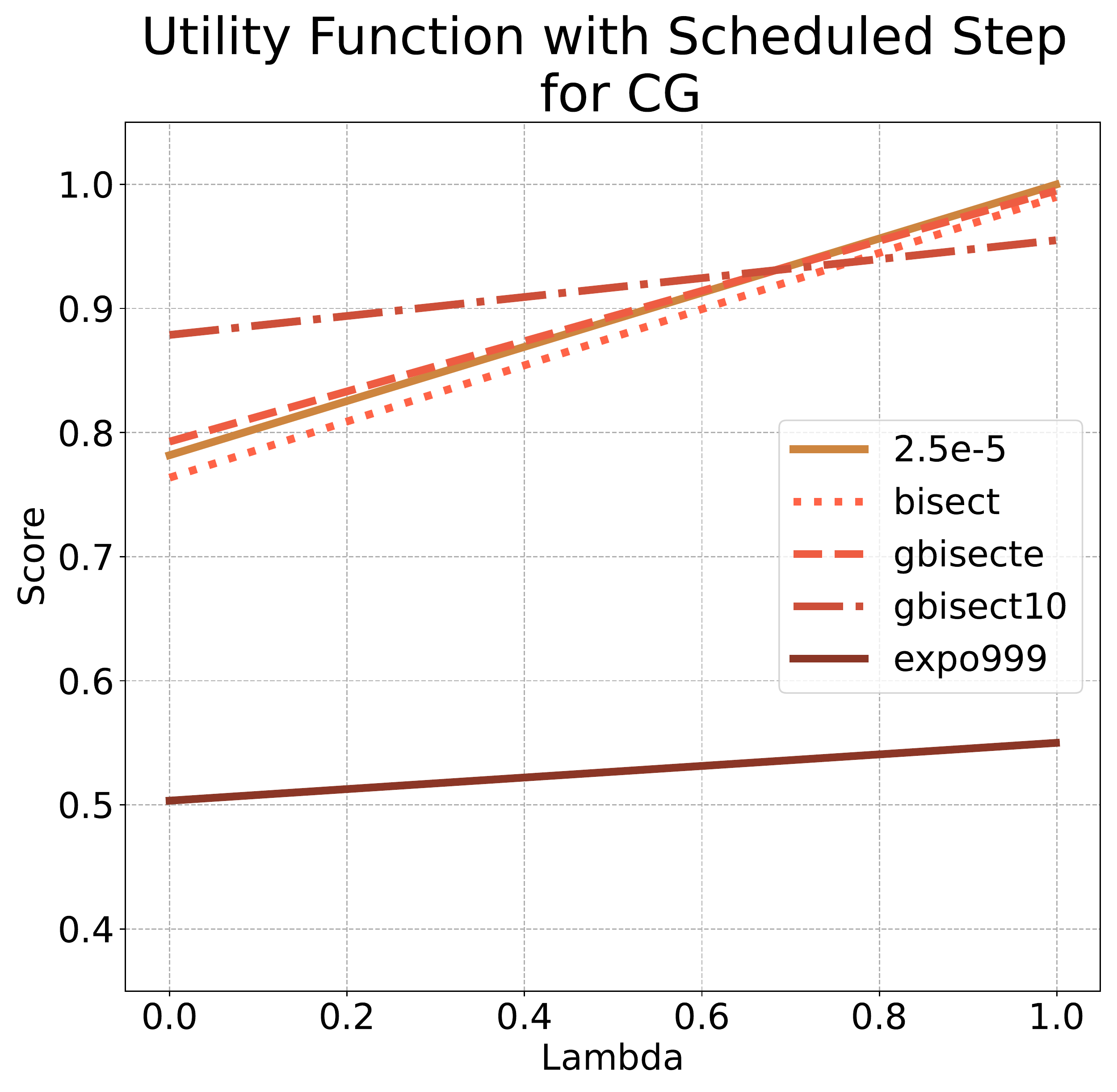}
    \end{subfigure}
    \caption[Utility scores of Gradient Descent and Conjugate Gradient with constant and scheduled step with respect to $\lambda$.]{\textbf{Utility scores of Gradient Descent and Conjugate Gradient with constant and scheduled step with respect to $\lambda$.} The scores correspond to results of relaxations of 200 structures. For all rules except GD with \num{0.00001} and CG with \num{0.000025} the initial step size is $s=\overline{s}$ and its lower bound is $\underline{s}$.}
    \label{fig:uall}
\end{figure}

The results from scheduled step size experiments, much like the previous experiments, show that different step size recipes provide high utilities for different preferences. Various hyperparameter configurations show that when the step size remains in small values then more structures can be relaxed, while larger step sizes can increase convergence speed. It is also shown that each scheduling rule is effective with a different algorithm between Gradient Descent and Conjugate Gradient, as seen in Figure~\ref{fig:uall}. We find that \textit{expo} is more beneficial to Gradient Descent, while \textit{gbisect} with $\beta=10$ is more beneficial to Conjugate Gradient, and both have roughly the same performance with \textit{bisect}. In the case of Gradient Descent we can observe the same monotonicity as with the constant step sizes in Figure~\ref{fig:uall}. No crossing lines exist, meaning that, when a step size scheduling rule is quick to converge, it is also able to relax more structures. However, further results show that there is at least one failed experiment when the step size is scheduled to drop many times during each experiment run, such as when using \textit{bisect} or \textit{expo} with $\gamma=0.999$ and $s_0=\overline{s}$, thus a constant step size is optimal for $\lambda\rightarrow 1$.

Due the great performance of Gradient Descent with \textit{expo999}, we experimented with different values of the parameter $\gamma$. All rules \textit{expo99}, \textit{expo999}, \textit{expo9999} have relaxed all structures in the same group of 40 structures, thus their speed is what distinguishes them. On the contrary, as long as Conjugate Gradient is concerned, different rules show different score order depending on the preferences case. While \textit{gbisect} with $\beta=10$ provides the lowest iteration number, it does not manage to relax a considerable amount of structures. In other words, while for $\lambda\leq 0.6$ \textit{gbisect} with $\beta=10$ is the best choice, for $\lambda > 0.6$ other rules with continuous and small reductions, such as \textit{bisect}, \textit{gbisect} with $\beta = e$, and especially the constant $s=\num{0.000025}$, have the best results.

\begin{figure}
    \centering
    \begin{subfigure}[c]{0.49\linewidth} 
    \includegraphics[width=\textwidth]{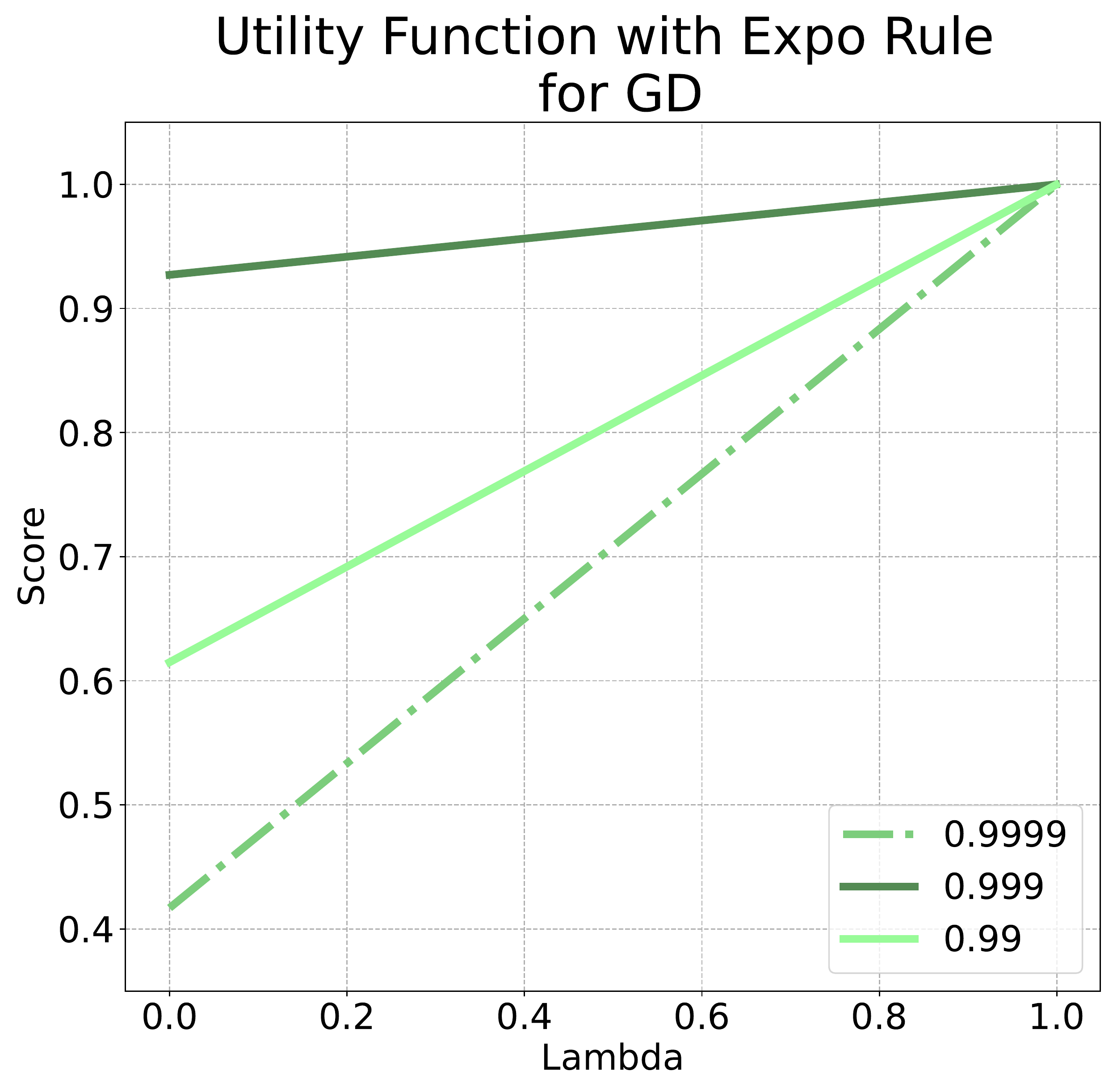}
    \caption{Utility scores of GD with rule \textit{expo} when $\gamma$ is 0.99, 0.999 and 0.9999. The initial step size is $s_1 = \num{0.001}$.}
    \label{fig:ues}
    \end{subfigure}

    \begin{subfigure}[t]{0.49\linewidth}  
    \includegraphics[width=\textwidth]{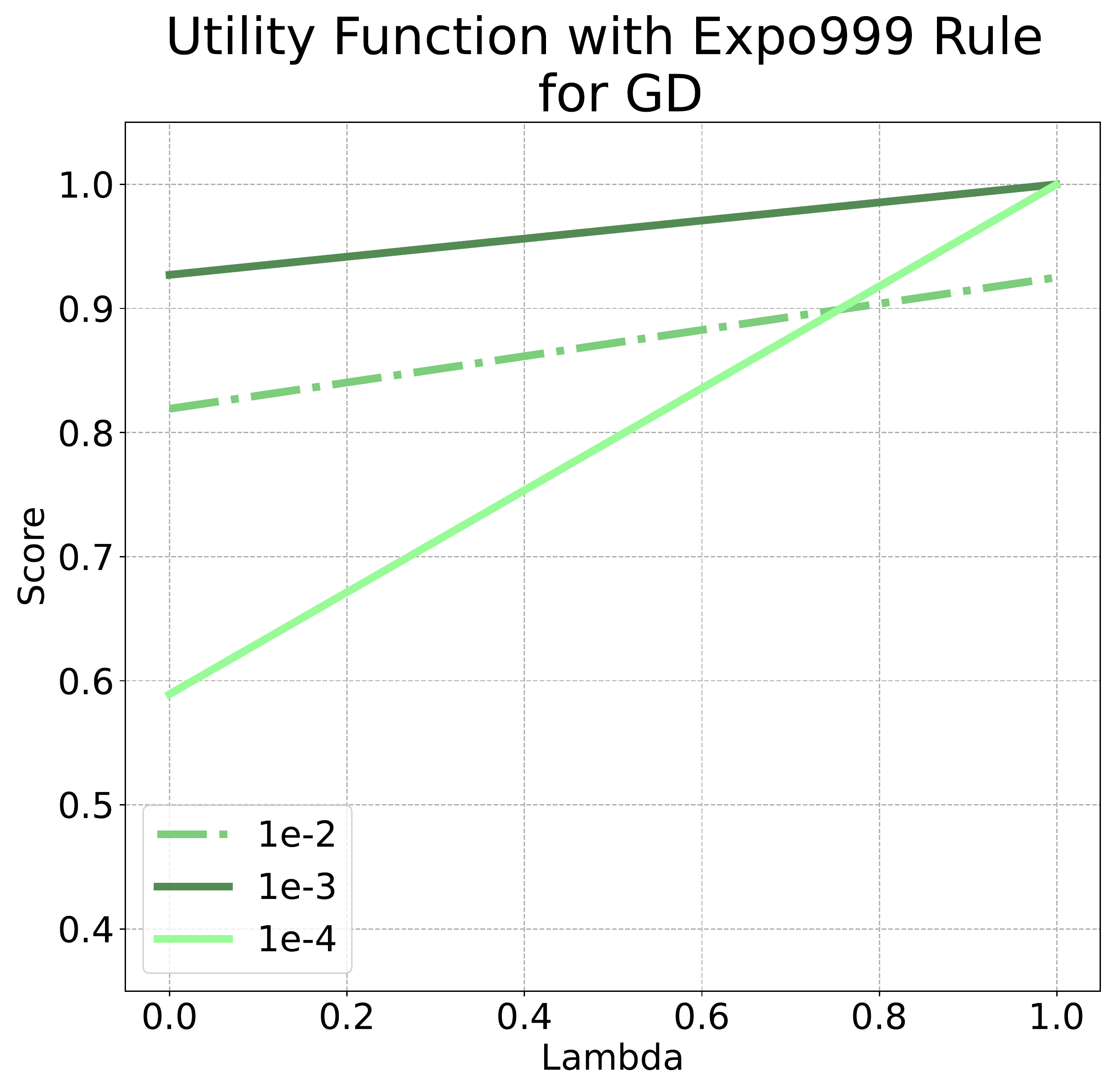}
    \caption[Utility scores of GD with the exponential rule (expo) with different initial step sizes]{Utility scores of GD with \textit{expo} when $\gamma = 0.999$ and the initial step size is \num{0.01}, \num{0.001}, \num{0.0001}.}
    \label{fig:ues999}    
    \end{subfigure}
    \begin{subfigure}[t]{0.49\linewidth}
    \includegraphics[width=\textwidth]{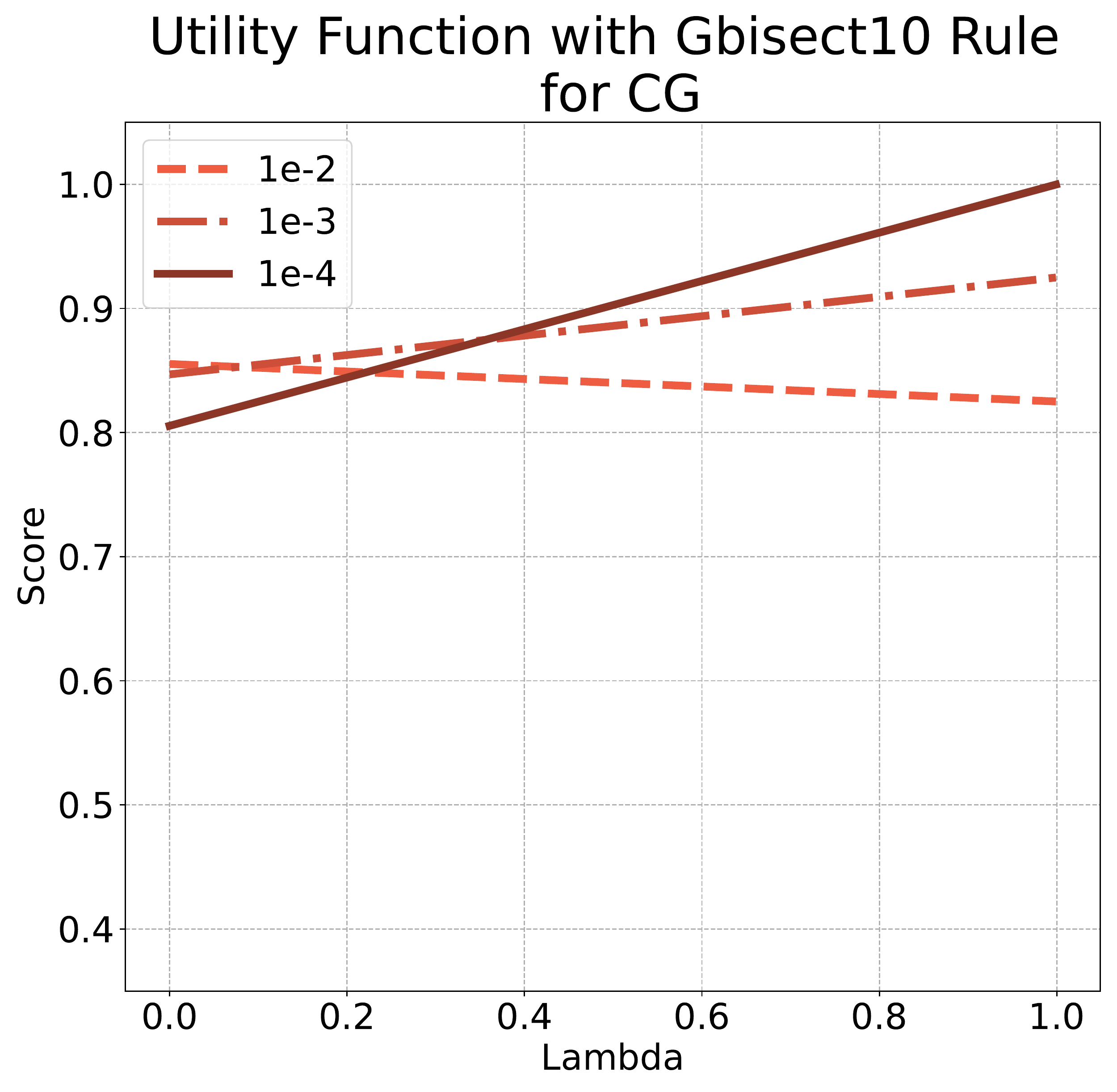}
    \caption[Utility scores of CG  with the gradient norm-scheduled rule with different initial step sizes]{Utility scores of CG with the gradient norm-scheduled rule \textit{gbisect} when $\beta = 10$ and the initial step size is \num{0.01}, \num{0.001}, \num{0.0001}.}
    \label{fig:ugsb10}
    \end{subfigure}

    \caption{\textbf{Utility scores with different values for various parameters.} Gradient Descent \textit{(green)} and Conjugate Gradient \textit{(red)} combined with the exponential and gradient norm scheduling rules respectively. The plots represent a parameter exploration in order to understand how changes affect the result. The scores correspond to experiments with 40 structures per parameter value.}
\end{figure}

Another parameter change of exponential scheduled step size was tested, concerning its initial value. We increased it from $s=\num{0.001}$ to $s=10^{-2}$ and decreased it from $s=\num{0.001}$ to $s=\num{0.0001}$. The analysis of this change through the utility function showed that $\overline{s}$ was the best choice. More specifically, Gradient Descent with \textit{expo999} and first step size $s=\overline{s}$ has the best scores for all $\lambda\in [0,1)$ compared to $s=\num{0.0001}$ and $s=\num{0.01}$, as seen in Figure~\ref{fig:ues999}. A similar initial step size analysis for \textit{gbisect10} is provided, with the same three different values like before and their utility scores, as found in Figure~\ref{fig:ugsb10}. It appears that when a low number of iterations is important, $s=\num{0.01}$ and  $s=\num{0.001}$ compete for the best option. As $\lambda\rightarrow 1$ it is implied that a large number of structures being successfully relaxed is our preference and speed is unimportant, consequently, a small step size is the safest and most reliable option.

\subsubsection{Analysis} 
The utility scores in the first two columns of the Tables~\ref{tab:uvalsGD},~\ref{tab:uvalsCG} with $\lambda\in [0,0.5]$  which correspond to constant step sizes confirm that Conjugate Gradient can greatly reduce the number of iterations compared to Gradient Descent. In conjuction with the last column, it is apparent that Conjugate Gradient's utility is the best in all cases with consistenly high success rate, thus more robust. We believe that the conjugacy of the produced direction vectors, along with the `memory' that the updating scheme carries, can more accurately and quickly traverse the PES, hence the results.

\begin{table}
\centering
\begin{tabular}{lllccc|} 
 & & & \multicolumn{3}{c}{$\lambda$} \\ \cline{4-6} 
 & \multicolumn{1}{l}{} & \multicolumn{1}{l|}{} & \multicolumn{1}{c|}{0} & \multicolumn{1}{c|}{0.5} & \multicolumn{1}{c|}{1} \\ \cline{2-6} 
\multicolumn{1}{c|}{\parbox[t]{2mm}{\multirow{13}{*}{\rotatebox[origin=c]{90}{method}}}} &                            \multicolumn{1}{c|}{$\num{0.00001}$} & \multicolumn{1}{l|}{\multirow{5}{*}{const}} & 0.456 & 0.728 & \textbf{1.000} \\ 
\multicolumn{1}{c|}{} & \multicolumn{1}{c|}{$\num{0.000025}$} & \multicolumn{1}{l|}{} & 0.374 & 0.434 & 0.495 \\ 
\multicolumn{1}{c|}{} & \multicolumn{1}{c|}{$\num{0.0000775}$} & \multicolumn{1}{l|}{} & 0.069 & 0.072 & 0.075 \\ 
\multicolumn{1}{c|}{} & \multicolumn{1}{c|}{$\num{0.0001}$} & \multicolumn{1}{l|}{} & 0.060 & 0.063 & 0.065 \\ 
\multicolumn{1}{c|}{} & \multicolumn{1}{c|}{$\num{0.001}$} & \multicolumn{1}{l|}{} & 0.003 & 0.002 & 0.000 \\ \cline{2-3}
\multicolumn{1}{c|}{} & \multicolumn{1}{c|}{$\num{0.001}$} & \multicolumn{1}{l|}{bisect} &  0.766 & 0.878 & 0.990 \\ \cline{2-3}
\multicolumn{1}{c|}{} & \multicolumn{1}{c|}{$\num{0.001}$} & \multicolumn{1}{l|}{gbisecte} & 0.795 & 0.895 & 0.995 \\ \cline{2-3}
\multicolumn{1}{c|}{} & \multicolumn{1}{c|}{$\num{0.001}$} & \multicolumn{1}{l|}{gbisect10} & 0.420 & 0.453 & 0.485 \\ \cline{2-3}
\multicolumn{1}{c|}{} & \multicolumn{1}{c|}{$\num{0.001}$} & \multicolumn{1}{l|}{expo9999} & 0.417 & 0.708 & \textbf{1.000} \\ \cline{2-3}
\multicolumn{1}{c|}{} & \multicolumn{1}{c|}{$\num{0.0001}$} & \multicolumn{1}{l|}{\multirow{3}{*}{expo999}} & 0.589 & 0.795 & \textbf{1.000} \\ 
\multicolumn{1}{c|}{} & \multicolumn{1}{c|}{$\num{0.001}$} & \multicolumn{1}{l|}{} & \textbf{0.929} & \textbf{0.962} & 0.995 \\ 
\multicolumn{1}{c|}{} & \multicolumn{1}{c|}{$\num{0.01}$} & \multicolumn{1}{l|}{} & 0.819 & 0.872 & 0.925 \\ \cline{2-3}
\multicolumn{1}{c|}{} & \multicolumn{1}{c|}{$\num{0.001}$} & \multicolumn{1}{l|}{expo99} & 0.615 & 0.807 & \textbf{1.000}  \\ \cline{2-6}
\end{tabular}
\caption{\textbf{Utility scores for all scheduling rules with Gradient Descent using function $u_{\text{FP}}$.} The scores of constant steps, bisect, gbisecte, gbisect10 and $\num{0.001}$\_expo999 correspond to experiments with 200 structures. The scores of rules expo9999, expo99, $\num{0.0001}$\_expo999 and $\num{0.01}$\_expo999 correspond to experiments on 40 structures.}
\label{tab:uvalsGD}
\end{table}
\begin{table}
\centering
\begin{tabular}{lllccc|}
 & & & \multicolumn{3}{c}{$\lambda$} \\ \cline{4-6} 
 & \multicolumn{1}{l}{} & \multicolumn{1}{l|}{} & \multicolumn{1}{c|}{0} & \multicolumn{1}{c|}{0.5} & \multicolumn{1}{c|}{1} \\ \cline{2-6} 
\multicolumn{1}{c|}{\parbox[t]{2mm}{\multirow{11}{*}{\rotatebox[origin=c]{90}{method}}}} &                            \multicolumn{1}{c|}{$\num{0.00001}$} & \multicolumn{1}{l|}{\multirow{5}{*}{const}} & 0.456 & 0.728 & \textbf{1.000} \\ 
\multicolumn{1}{c|}{} & \multicolumn{1}{c|}{$\num{0.000025}$} & \multicolumn{1}{l|}{} & 0.782 & 0.891 & \textbf{1.000} \\  
\multicolumn{1}{c|}{} & \multicolumn{1}{c|}{$\num{0.0000775}$} & \multicolumn{1}{l|}{} & 0.682 & 0.713 & 0.745 \\ 
\multicolumn{1}{c|}{} & \multicolumn{1}{c|}{$\num{0.0001}$} & \multicolumn{1}{l|}{} & 0.537 & 0.558 & 0.580 \\ 
\multicolumn{1}{c|}{} & \multicolumn{1}{c|}{$\num{0.001}$} & \multicolumn{1}{l|}{} & 0.263 & 0.131 & 0.000 \\ \cline{2-3}
\multicolumn{1}{c|}{} & \multicolumn{1}{c|}{$\num{0.001}$} & \multicolumn{1}{l|}{bisect} & 0.764 & 0.877 & 0.99 \\ \cline{2-3}
\multicolumn{1}{c|}{} & \multicolumn{1}{c|}{$\num{0.001}$} & \multicolumn{1}{l|}{gbisecte} & 0.793 & 0.894 & 0.995 \\ \cline{2-3}
\multicolumn{1}{c|}{} & \multicolumn{1}{c|}{$\num{0.0001}$} & \multicolumn{1}{l|}{\multirow{3}{*}{gbisect10}} & 0.805 & 0.903 & \textbf{1.000} \\ 
\multicolumn{1}{c|}{} & \multicolumn{1}{c|}{$\num{0.001}$} & \multicolumn{1}{l|}{} & \textbf{0.879} & \textbf{0.917} & 0.955 \\ 
\multicolumn{1}{c|}{} & \multicolumn{1}{c|}{$\num{0.01}$} & \multicolumn{1}{l|}{} & 0.855 & 0.840 & 0.825 \\ \cline{2-3}
\multicolumn{1}{c|}{} & \multicolumn{1}{c|}{$\num{0.001}$} & \multicolumn{1}{l|}{expo999} & 0.503 & 0.527 & 0.550 \\ \cline{2-6}
\end{tabular}
\caption{\textbf{Utility scores for all scheduling rules with Conjugate Gradient using function $u_{\text{FP}}$.} The scores of constant steps, bisect, gbisecte, gbisect10 and $\num{0.001}$\_expo999 correspond to experiments with 200 structures. The scores of rules $\num{0.0001}$\_gbisect10 and $\num{0.01}$\_gbisect10 correspond to experiments on 40 structures.}
\label{tab:uvalsCG}
\end{table}

Powell~\cite{Powell1984NonconvexMethod} has showed that, even with an exact line search, Conjugate Gradient must be combined with a changing step size that tends to zero. In our work we confirm that a decreasing step size can increase the utility of Conjugate Gradient in terms of speed, in other words, it can reduce the number of its iterations. This is shown by the large increase in its utility score for $\lambda\rightarrow 0$ in rows of Table~\ref{tab:uvalsCG} concerning \textit{gbisecte} and \textit{gbisect10} compared to the previous rows, concerning constant step. However, for the rest of the scheduling rules there is small or no increase. Scheduling the step size to be reduced according to an observed large drop in the gradient norm, like with \textit{gbisect10}, improves the algorithm's convergence speed and can in most cases relax a crystal structure successfully. Step size scheduling rules that frequently decrease the step size value, like \textit{bisect} and \textit{gbisecte}, improve the success rate of relaxations, so the utility with $\lambda\rightarrow 1$ is improved, but they do not perform well in terms of speed, thus the utility scores $\lambda\rightarrow 0$ are low.  We have already argued that small steps lead Conjugate Gradient updates to imitate Gradient Descent updates. This also means that its ability to keep information from previous iterations vanishes. Consequently, for Conjugate Gradient to be effective, the step size must be reduced at low pace, so that new directions can benefit from past directions visited. This way, the new directions will tend to be more orthogonal to the direction where the minimum lies and Conjugate Gradient can approach it in fewer steps~\cite{Blair2006ProblemYudin}. Interestingly, rules like \textit{bisect} and \textit{gbisecte} have slightly better utility scores with Gradient Descent compared to Conjugate Gradient. We can observe that Gradient Descent with a scheduling rule that reduces the step size frequently achieves better results. For a constant small decrease in step size, like when using \textit{expo999}, Gradient Descent can achieve almost optimal results. Conjugate Gradient's utility enhancement with large step sizes also justifies why rule \textit{gbisect10} with this algorithm is not as effective as \textit{expo999} with Gradient Descent; towards the last iterations the step size is reduced to almost $\num{0.00001}$, ergo Conjugate Gradient behaves like Gradient Descent and loses its previous convergence properties.

\section{Conclusion}

In this paper we have provided the derivation of a set of crucial equations used in crystal structure prediction. We tested Gradient Descent and Conjugate Gradient with a constant step size and other step size adaptive methods, thus providing the foundations for a direct comparison of the two in geometric optimisation of crystals. With these benchmarks in place, we intend to further investigate the algorithms' performance in relaxation when combined with proper line search, which will be designed to avoid skipping the closest \textit{true} local minimum of the function. The outcome of the experimentation process was that Conjugate Gradient is more trustworthy and efficient, due to its ability to adapt to the PES and mark a steady progress. We confirmed that a scheduled decrease to the step size reduces the iteration number, but there is a trade-off between quick convergence and the number of successful experiments. What is more, we showed that Gradient Descent and Conjugate Gradient benefit from different kinds of step size scheduling rules. We concluded that the values extracted from our constant step experiments provide the best schemes, with Gradient Descent and rule \textit{expo999} displaying optimal results for speed preferences and small constant step sizes displaying optimal results for success rate preferences. In the future, we will include second order methods to our study and extend our analysis to accommodate them accordingly.

\section{Acknowledgements}
This work has been funded by the Leverhulme Research Centre for Functional Materials Design.

\Urlmuskip=0mu plus 1mu
\bibliographystyle{plain}
\bibliography{references}

\appendix
\section{Proofs related to the energy function potential} \label{app1}
\subsection{Proof of Propositions~\ref{prop:real} and~\ref{prop:longrange} (See page~\pageref{prop:real})}
The following derivations involving the handling of the ions' potential field have been heavily inspired by the work of H.Lee and W.Cai~\cite{Lee2009EwaldSupercell}.
\begin{proof}
    For the proof of Propositions~\ref{prop:real} and ~\ref{prop:longrange}, we will be referring to the electric potential field ($\phi$) and the charge density distribution ($\rho$) of a point charge ($q_j$) at position $\ionp$:

    \begin{align*}
        \phi_j(\ionp) &= \frac{1}{4\pi\epsilon_0}\frac{q_j}{\|\ionp-\ionp_j\|}, \ \ionp,\ionp_j \in \mathbb{R}^3 \\
        \rho_j(\ionp) &= q_j\delta(\ionp-\ionp_j), \ \ionp,\ionp_j \in \mathbb{R}^3
    \end{align*}
    where $\epsilon_0$ is the vacuum permittivity and $\delta$ is the Dirac delta function. The potential field at position $\ionp = \ionp_i+\lattice$ generated by all ions with positions $\ionp_j, \ j\in[N], \ j\neq i$  is
    \begin{equation}
        \label{eq:potfield_i}
        \phi_j(\ionp_i) = \frac{1}{4\pi\epsilon_0}\sum_n\sum_{j=1}^{N'}\frac{q_j}{\|\ionp_{i,j,n}\|}, \ \ionp_i,\ionp_j \in \mathbb{R}^3 
    \end{equation}
    which gives the conditionally convergent energy potential
    \begin{equation*}
        \energy_{Coul}(\allionp,L) = \frac{1}{4\pi\epsilon_0}\sum_{i=1}^N\phi_j(\ionp_i)
    \end{equation*}
    The following techniques for the mentioned derivations are heavily based on the work of Lee and Cai~\cite{Lee2009EwaldSupercell}. 
    We can consider that the charge distribution of an ion extends in space as a Gaussian distribution
    \begin{equation}
        G_\sigma(\ionp) = \frac{1}{(2\pi\sigma^2 )^{3\slash2}}\exp\left(-\frac{r^2}{2\sigma^2}\right), \ \ionp\in \mathbb{R}^3
    \end{equation}
    Moreover, it is convenient to acknowledge that the delta function is actually the limit of the Gaussian distribution $G_\sigma$ with the standard deviation approaching zero $\sigma \rightarrow 0$
    \begin{equation*}
        \lim_{\sigma \rightarrow 0} G_\sigma(\ionp) = \delta(\ionp)
    \end{equation*}
    Then, we can discriminate between interactions close to charge $q_j$ and interactions in distance from it by adding and subtracting $G_\sigma(\ionp-\ionp_i)$ to charge density distribution
    \begin{equation}
    \begin{aligned}
    \label{eq:split}
        \rho_j(\ionp) &= \ \rho^S_j(\ionp) + \rho^L_j(\ionp) \\
        \rho^S_j(\ionp) &= \ q_j\delta(\ionp - \ionp_j) - q_jG_\sigma(\ionp - \ionp_j) \\
        \rho^L_j(\ionp) &= \ q_jG_\sigma(\ionp - \ionp_j)
    \end{aligned}
    \end{equation}
    The charge density distribution is connected to the potential field $\phi_i$ through Poisson's equation. When it takes the form of a Gaussian distribution, we have
    \begin{equation}
    \label{eq:Poisson}
        \nabla^2\phi_j(\ionp) = -\frac{\rho_j(\ionp)}{\epsilon_0} = -\frac{q_jG_\sigma(\ionp-\ionp_j)}{\epsilon_0}
    \end{equation}
    Here, the employment of some Gaussian charge distribution with standard deviation $\sigma$ serves as the means to embed convergence factors that ultimately convert the Coulombic energy potential function to an absolutely convergent summation. Now, since $G_\sigma$ is a function with only one independent variable $r$, we can express Poisson's equation in spherical coordinates. To be more specific, the Laplacian operator in terms of spherical coordinates appears as
    \begin{equation}
    \label{eq:LaplAll}
        \nabla^2 = \frac{1}{r^2}\frac{\partial}{\partial r}\left(r^2\frac{\partial}{\partial r}\right)+\frac{1}{r^2\sin{\theta}}\frac{\partial}{\partial\theta}\left(\sin{\theta}\frac{\partial}{\partial\theta}\right)+\frac{1}{r^2\sin^2{\theta}}\frac{\partial^2}{\partial\phi^2}
    \end{equation}
    We define as $\overline{\phi}_j(r)$, $\overline{G}_\sigma(r)$ the functions that correspond to $\phi_j(\ionp)$,$G_\sigma(\ionp)$ with the difference that the independent variable becomes the distance $r = \|\ionp-\ionp_j\|$ since it is the only quantity affected from the vector positions $\ionp$. Because of spherical symmetry and after some operations
    \begin{equation*}
    \begin{aligned}
        & \frac{1}{r^2}\frac{\partial}{\partial r}\left(r^2\frac{\partial}{\partial r}\overline{\phi}_j(r)\right) &=& \ \frac{1}{r^2}\left( r^2\frac{\partial^2\overline{\phi}_j}{\partial r^2}(r) + 2r\frac{\partial\overline{\phi}_j}{\partial r}(r)\right)\\
        & &=& \ \frac{\partial^2\overline{\phi}_j}{\partial r^2}(r) + \frac{2}{r}\frac{\partial\overline{\phi}_j}{\partial r}(r) \\
        & &=& \ \frac{1}{r}\frac{\partial^2}{\partial r^2}\left(r\overline{\phi}_j(r)\right)
    \end{aligned}
    \end{equation*}
    we get the following simple equation from Equation~(\ref{eq:Poisson})
    \begin{equation*}
        \frac{1}{r}\frac{\partial^2}{\partial r^2}\left(r\overline{\phi}_j(r)\right) = -\frac{\overline{G}_\sigma(r)}{\epsilon_0}
    \end{equation*}
    By integration we arrive to
    \begin{equation}
    \begin{aligned}
        & & r\overline{\phi}_j(r) =& \ \frac{\sigma}{\epsilon_0}\int_0^r \overline{G}_\sigma(r)dr \\
        & & =& \ \frac{\sigma}{\epsilon_0}\frac{1}{(2\pi\sigma^2 )^{3\slash2}}\sqrt{\frac{\pi}{2}}\sigma \erf{\frac{r}{\sqrt{2}\sigma}} \\
        & \Rightarrow & \overline{\phi}_j(r) =& \ \frac{1}{4\pi\epsilon_0r}\erf{\frac{r}{\sqrt{2}\sigma}}
    \end{aligned}
    \end{equation}
    in which $\erf{z} = \frac{2}{\sqrt{\pi}}\int_0^z\exp{(-t^2)}dt$ is the error function. After the split of charge density distribution in Equation~(\ref{eq:split}), the potential field caused by ion $j$ is correspondingly split into
    \begin{align}
    \label{eq:split1}
        \phi_j(\ionp) &= \ \phi_j^S(\ionp) + \phi_j^L(\ionp) & \\
        \phi_j^S(\ionp) &= \  \frac{1}{4\pi\epsilon_0}\frac{q_j}{\|\ionp-\ionp_j\|}\left[1-\erf{\frac{\|\ionp-\ionp_j\|}{\sqrt{2}\sigma}}\right] = \frac{1}{4\pi\epsilon_0}\frac{q_j}{\|\ionp-\ionp_j\|}\erfc{\frac{\|\ionp-\ionp_j\|}{\sqrt{2}\sigma}} & \label{eq:potfieldS} \\
        \phi_j^L(\ionp) &= \ \frac{1}{4\pi\epsilon_0}\frac{q_j}{\|\ionp-\ionp_j\|}\erf{\frac{\|\ionp-\ionp_j\|}{\sqrt{2}\sigma}} \label{eq:potfieldL} 
    \end{align}
    Finally, we can calculate the electrostatic energy caused by short ranged interactions using Equation~(\ref{eq:potfieldS}). This range is determined by the term $\erfc{}$, which truncates the summation for long distances and converges absolutely
    \begin{equation}
    \begin{aligned}
        \label{eq:potentialShort}
        \energy^S_{Coul}(\allionp, L) &= \frac{k_e}{2}\sum_{i=1}^N\phi_j^S(\ionp) \\
        &= \frac{k_e}{2}\sum_n\sum_{i=1}^N\sum_{j=1}^{N^{'}}\frac{q_iq_j}{\|\ionp_{i,j,n}\|}\erfc{\alpha\|\ionp_{i,j,n}\|}
    \end{aligned}
    \end{equation}
    with $\alpha = \frac{1}{\sqrt{2}\sigma}$ and $k_e = \frac{1}{4\pi\epsilon_0}$.
    In the same fashion, for long ranged interactions we have
    \begin{equation}
    \begin{aligned}
        \energy_{Coul}^L(\allionp,L) &= \frac{k_e}{2}\sum_n\sum_{i=1}^N\sum_{j=1}^N\phi_j^L(\ionp) \\
        &= \frac{k_e}{2}\sum_n\sum_{i=1}^N\sum_{j=1}^N \frac{q_iq_j}{\|\ionp_{i,j,n}\|}\erf{\frac{\|\ionp_{i,j,n}\|}{\sqrt{2}\sigma}}
    \end{aligned} 
    \end{equation}

\end{proof}

\subsection{Proof of Proposition~\ref{prop:reciprocalexp} (See page~\pageref{prop:reciprocalexp})}
The proof of this proposition has been inspired by the techniques presented by D. Wang et al.~\cite{Wang2019EwaldDipoles}.
\begin{proof}
Let us define a function
\begin{equation}
    \label{eq:definePeriodic}
    f(\ionp_i-\ionp_j) = f(\ionp) = \sum_n\frac{1}{\|\ionp+\lattice\|}\erf{\alpha\|\ionp + \lattice\|} = f(\ionp + \lattice)
\end{equation}
according to Equation~(\ref{eq:potfieldL}). The function $f$ is periodic for intervals equal to the unit cell lengths defined by the three linearly independent vectors $\lattice = n_1\vect_1+n_2\vect_2+n_3\vect_3$. We can expand $f$ as a Fourier series of exponential functions
\begin{equation}
\label{eq:expSeries}
    f(\ionp) = \sum_{m\in\mathbb{Z}^3}h(m)\exp{(i\rlattice\ionp)}
\end{equation}
in which we have used the reciprocal vectors' property $\rvect_i\cdot\vect_j = 2\pi\delta_{ij},\ i,j\in\{1,2,3\}$ where $\delta_{ij}$ is the Kronecker delta, and $h(m)$ are Fourier coefficients. It follows from Equation~(\ref{eq:definePeriodic}) that the multipliers of the series of Equation~(\ref{eq:expSeries}) are
\begin{align*}
    h(m) &= \frac{1}{|V|}\iiint\limits_V f(\ionp)\exp{(-i\rlattice\ionp)}d^3\ionp \\
    &= \frac{1}{|V|}\iiint\limits_V\sum_n\frac{1}{\|\ionp+\lattice\|}\erf{\frac{\|\ionp+\lattice\|}{\sqrt{2}\sigma}}\exp{(-i\rlattice\ionp)}d^3\ionp \, , 
\end{align*} where $V$ is the unit cell, and 
$|V|$ is its the volume. Following the derivation of Wang at al.~\cite{Wang2019EwaldDipoles} we change $\ionp$ into spherical coordinates $\rho,\theta,\phi$, which describe vector's $\ionp$ length, angle from x-axis in the xy-plane and angle from z-axis respectively. We can assume that the z-axis is parallel to $\rlattice$ as follows
\begin{align*}
    h(m) &= \frac{1}{|V|}\iiint\limits_V f(\ionp)\exp{(-i\rlattice\ionp)}d^3\ionp \\
    &= \frac{1}{|V|}\int_0^{2\pi}\int_0^{\pi}\int_0^\infty \frac{1}{\rho}\erf{\frac{\rho}{\sqrt{2}\sigma}}\exp{(-i\|\rlattice\|\rho\cos{(\phi)})}d\rho(\rho d\phi)(\rho\sin{\phi}d\theta) \\
    &= \frac{2\pi}{|V|}\int_0^\infty\rho \cdot \erf{\frac{\rho}{\sqrt{2}\sigma}}\int_0^\pi\exp{(-i\|\rlattice\|\rho\cos{(\phi)})}d\phi d\rho \\
    &= \frac{4\pi}{|V|\|\rlattice\|}\int_0^\infty\sin{(\|\rlattice\|\rho)}\erf{\frac{\rho}{\sqrt{2}\sigma}}d\rho 
\end{align*}
Next, we perform a variable change to replace $\|\rlattice\|\rho$ with some $x$ and $\frac{1}{\sqrt{2}\sigma}$ with $\alpha$
\begin{equation}
\begin{aligned}
    \label{eq:multVal}
    h(m) &= \frac{4\pi}{|V|\|\rlattice\|^2}\int_0^\infty\sin{(x)}\erf{\frac{\alpha x}{\|\rlattice\|}}dx \\
    &= \frac{4\pi}{|V|\|\rlattice\|^2}\exp{\bigg(-\frac{\|\rlattice\|^2}{4\alpha^2}\bigg)} 
\end{aligned}
\end{equation}
Putting everything together, we get
\begin{equation*}
\begin{aligned}
    \energy_{Coul}^L 
    &= \frac{k_e}{2}\sum_{i=1}^N\sum_{j=1}^N \sum_{m\in\mathbb{Z}^3}h(m)\exp{(i\rlattice\ionp)} \\
    &= \frac{k_e}{2}\sum_{i=1}^N\sum_{j=1}^N \sum_{m\in\mathbb{Z}^3}\frac{4\pi}{V\|\rlattice\|^2}\exp{\bigg(-\frac{\|\rlattice\|^2}{4\alpha^2}\bigg)}\exp{(i\rlattice\ionp)}
\end{aligned}
\end{equation*}
\end{proof}

\subsection{Proof of Proposition~\ref{prop:reciprocal} (See page~\pageref{prop:reciprocal})}
\begin{proof}
Let $f$ be defined as in the proof of Proposition~\ref{prop:reciprocalexp}. Let us also define a real function $g:\mathbb{R}^3\rightarrow\mathbb{R}$. We can express any vector $\ionp$ using the coordinate system constructed by the lattice vectors, such that
\begin{equation}
    \label{eq:vectorLatticeCoords}
    \ionp = x_1\frac{\vect_1}{\|\vect_1\|}+x_2\frac{\vect_2}{\|\vect_2\|}+x_3\frac{\vect_3}{\|\vect_3\|}
\end{equation} and then  $f(\ionp)=g(x_1,x_2,x_3)$. Hence, the periodicity of $f$ is now expressed through $g$ with 

\begin{equation*}
    g(x_1,x_2,x_3) = g(x_1+l'_1,x_2+l'_2,x_3 + l'_3), \, \mbox{ where } \, l'_j \in \{0,\|\vect_j\|\}, \, \mbox{ for } \, j \in \{1,2,3\}
\end{equation*}

Then, $g$ can be expanded to a Fourier series using a set of orthonormal functions. \\
Let $u_t = \frac{2\pi m_t}{\|\vect_t\|}, \ t\in\{1,2,3\}$ 
\begin{equation}
\begin{gathered}
    g(x_1,x_2,x_3) = \sum_{m_1=0}^{\infty}\sum_{m_2=0}^{\infty}\sum_{m_3=0}^{\infty} a_m\cos{(u_1x_1)}\cos{(u_2x_2)}\cos{(u_3x_3)} + \\
    b_m\sin{(u_1x_1)}\cos{(u_2x_2)}\cos{(u_3x_3)} +
    c_m\cos{(u_1x_1)}\sin{(u_2x_2)}\cos{(u_3x_3)} + \\
    d_m\sin{(u_1x_1)}\sin{(u_2x_2)}\cos{(u_3x_3)} +
    \alpha_m\cos{(u_1x_1)}\cos{(u_2x_2)}\sin{(u_3x_3)} + \\
    \beta_m\sin{(u_1x_1)}\cos{(u_2x_2)}\sin{(u_3x_3)} +
    \gamma_m\cos{(u_1x_1)}\sin{(u_2x_2)}\sin{(u_3x_3)} + \\
    \delta_m\sin{(u_1x_1)}\sin{(u_2x_2)}\sin{(u_3x_3)}     \label{eq:cosineSeries}
\end{gathered}
\end{equation}
with multipliers indexed by a triplet $m=(m_1,m_2,m_3), \ m_1,m_2,m_3\in\mathbb{Z}$. We can replace all terms with Euler's identity
\begin{align*}
\cos{(u_tx_t)} = \frac{1}{2}[ \exp{(iu_tx_t)} + \exp{(-iu_tx_t)}], \ t\in \{1,2,3\} \\
\sin{(u_tx_t)} = \frac{1}{2i}[ \exp{(iu_tx_t)} - \exp{(-iu_tx_t)} ], \ t\in \{1,2,3\}
\end{align*}
and get a series with terms comprising product combinations 
\[
    \exp{(\pm iu_1x_1)}\exp{(\pm iu_2x_2)}\exp{(\pm iu_3x_3)}=\exp{[i(\pm u_1x_1\pm u_2x_2\pm u_3x_3)]}. 
\] We can regroup the multipliers $a_m,b_m,c_m,d_m,\alpha_m,\beta_m,\gamma_m,\delta_m$ to accompany each unique exponential term of the form $\exp{[i(\pm u_1x_1\pm u_2x_2\pm u_3x_3)]}$ and get new multipliers of the form $\eta_{m(j)} = \frac{1}{8}(\pm a_m\pm ib_m\pm ic_m\pm d_m\pm i\alpha_m\pm \beta_m\pm \gamma_m\pm i\delta_m)$. Let us define the set, $\Theta$, of all such coefficients,
$\Theta =\{\frac{1}{8}(\oplus \, a_m\oplus ib_m\oplus ic_m\oplus d_m\oplus i\alpha_m\oplus \beta_m\oplus \gamma_m\oplus i\delta_m) \, : \, \oplus \in \{+, -\} \}$, so that
\begin{equation*}
\begin{gathered}
    g(x_1,x_2,x_3) = \sum_{m_1,m_2,m_3\in\mathbb{N}}\{ \\
    \eta_{m(1)}(m)\exp{[i(u_1x_1+u_2x_2+u_3x_3)]} + \eta_{m(2)}(m)\exp{[i(u_1x_1-u_2x_2+u_3x_3)]} + \\
    \eta_{m(3)}(m)\exp{[i(-u_1x_1+u_2x_2+u_3x_3)]} + \eta_{m(4)}(m)\exp{[i(-u_1x_1-u_2x_2+u_3x_3)]} + \\
    \eta_{m(5)}(m)\exp{[i(u_1x_1+u_2x_2-u_3x_3)]} + \eta_{m(6)}(m)\exp{[i(u_1x_1-u_2x_2-u_3x_3)]} + \\
    \eta_{m(7)}(m)\exp{[i(-u_1x_1+u_2x_2-u_3x_3)]} + \eta_{m(8)}(m)\exp{[i(-u_1x_1-u_2x_2-u_3x_3)]}\}  ,
\end{gathered} 
\end{equation*}
where $\eta_{m(j)}(m) \in \Theta$, are appropriate coefficients for $j \in \{1,2,\ldots,8\}$. We observe that
\begin{equation}
\begin{gathered}
    \label{eq:etaSeries}
    g(x_1,x_2,x_3) = \\
    \eta_{m(1)}\exp{[i(u_1x_1+u_2x_2+u_3x_3)]}+ [\eta_{m(1)}]^*\exp{[-i(u_1x_1+u_2x_2+u_3x_3)]} + \\
    \eta_{m(2)}\exp{[i(u_1x_1-u_2x_2+u_3x_3)]}+ [\eta_{m(2)}]^*\exp{[-i(u_1x_1-u_2x_2+u_3x_3)]} + \\
    \eta_{m(3)}\exp{[i(-u_1x_1+u_2x_2+u_3x_3)]}+ [\eta_{m(3)}]^*\exp{[-i(-u_1x_1+u_2x_2+u_3x_3)]} + \\
    \eta_{m(6)}\exp{[i(u_1x_1+u_2x_2-u_3x_3)]}+ [\eta_{m(6)}]^*\exp{[-i(u_1x_1+u_2x_2-u_3x_3)]}
\end{gathered}
\end{equation}

As a result of Equations~(\ref{eq:vectorLatticeCoords}),~(\ref{eq:expSeries}), we have that 
\[
    \exp{(i\rlattice\ionp)} = \exp{[i(u_1x_1 + u_2x_2 + u_3x_3)]}
\]
and

\begin{equation}
\begin{aligned}
    \label{eq:expSeriesMult}
    f(\ionp) &= g(x_1,x_2,x_3) = \\
    &\sum_{m\in\mathbb{N}^3}h(m_1,m_2,m_3)\exp{[i(m_1,m_2,m_3)\rlattice\ionp]} + \\
    &\sum_{m\in\mathbb{N}^3}h[-(m_1,m_2,m_3)]\exp{[-i(m_1,m_2,m_3)\rlattice\ionp]} + \\
    &\sum_{m\in\mathbb{N}^3}h(m_1,-m_2,m_3)\exp{[i(m_1,-m_2,m_3)\rlattice\ionp]} + \\
    &\sum_{m\in\mathbb{N}^3}h[-(m_1,-m_2,m_3)]\exp{[-i(m_1,-m_2,m_3)\rlattice\ionp]}+ \\
    &\sum_{m\in\mathbb{N}^3}h(-m_1,m_2,m_3)\exp{[i(-m_1,m_2,m_3)\rlattice\ionp]} + \\
    &\sum_{m\in\mathbb{N}^3}h[-(-m_1,m_2,m_3)]\exp{[-i(-m_1,m_2,m_3)\rlattice\ionp]} + \\
    &\sum_{m\in\mathbb{N}^3}h(m_1,m_2,-m_3)\exp{[i(m_1,m_2,m_3)\rlattice\ionp]} + \\
    &\sum_{m\in\mathbb{N}^3}h[-(m_1,m_2,-m_3)]\exp{[-i(m_1,m_2,-m_3)\rlattice\ionp]} 
\end{aligned}
\end{equation}

The uniqueness of the multipliers of trigonometric series and Equations (\ref{eq:etaSeries}), (\ref{eq:expSeriesMult}) show that the multipliers of the exponential and cosine-sine series are connected through the following relation
\begin{equation}
    \begin{aligned}
    \label{eq:multsRelation}
    h(m) &= h(m_1,m_2,m_3) \\
         &= \frac{1}{8}(a_m-ib_m-ic_m-d_m-i\alpha_m-\beta_m-\gamma_m+i\delta_m) \\
         &= \eta_{m1}
    \end{aligned}
\end{equation}
in accordance to the sign of each of the integers $m_1,m_2,m_3$. As a consequence of the orthogonality of the functions in Equation~(\ref{eq:cosineSeries}), we can exactly calculate its multipliers. Let $m_1',m_2',m_3' \in \mathbb{Z}$ and space $\mathcal{V} = [-\frac{\|\vect_1\|}{2},\frac{\|\vect_1\|}{2}]\times[-\frac{\|\vect_2\|}{2},\frac{\|\vect_2\|}{2}]\times[-\frac{\|\vect_3\|}{2},\frac{\|\vect_3\|}{2}]$ define a unit cell volume image. In order to obtain the value of multiplier $a_m$, we take advantage of the orthogonality of the summation terms and multiply both sides of Equation~(\ref{eq:cosineSeries}) with $\cos{(\frac{2\pi m_1'}{\|\vect_1\|}x_1)}\cos{(\frac{2\pi m_2'}{\|\vect_2\|}x_2)}\cos{(\frac{2\pi m_3'}{\|\vect_3\|}x_3)}$ and integrate over $\mathcal{V}$. This results into eliminating any terms with a multiplier different from $a_m$ and
\begin{gather*}
    \int_{-\frac{\|\vect_1\|}{2}}^{\frac{\|\vect_1\|}{2}}\int_{-\frac{\|\vect_2\|}{2}}^{\frac{\|\vect_2\|}{2}}\int_{-\frac{\|\vect_3\|}{2}}^{\frac{\|\vect_3\|}{2}} g(x_1,x_2,x_3)\cdot
    \cos{\left( \frac{2\pi m_1'}{\|\vect_1\|}x_1 \right)}\cos{\left( \frac{2\pi m_2'}{\|\vect_2\|}x_2 \right)}\cos{\left( \frac{2\pi m_3'}{\|\vect_3\|}x_3 \right)}dx_1dx_2dx_3
    \\
    = \int_{-\frac{\|\vect_1\|}{2}}^{\frac{\|\vect_1\|}{2}}\int_{-\frac{\|\vect_2\|}{2}}^{\frac{\|\vect_2\|}{2}}\int_{-\frac{\|\vect_3\|}{2}}^{\frac{\|\vect_3\|}{2}} a_m\cdot \left[ \cos{\left( \frac{2\pi m_1}{\|\vect_1\|}x_1 \right)}\cos{\left( \frac{2\pi m_2}{\|\vect_2\|}x_2 \right)}\cos{\left( \frac{2\pi m_3}{\|\vect_3\|}x_3 \right)} \right] \cdot
    \\
    \left[ \cos{\left( \frac{2\pi m_1'}{\|\vect_1\|}x_1 \right)}\cos{\left( \frac{2\pi m_2'}{\|\vect_2\|}x_2 \right)}\cos{\left( \frac{2\pi m_3'}{\|\vect_3\|}x_3 \right)} \right]
    dx_1dx_2dx_3 + 0
    \\
    = \int_{-\frac{\|\vect_1\|}{2}}^{\frac{\|\vect_1\|}{2}}\int_{-\frac{\|\vect_2\|}{2}}^{\frac{\|\vect_2\|}{2}}\int_{-\frac{\|\vect_3\|}{2}}^{\frac{\|\vect_3\|}{2}} \frac{a_m}{8}\cdot \left\{ \cos{\left[\frac{2\pi (m_1-m_1')}{\|\vect_1\|}x_1\right]}+\cos{\left[\frac{2\pi (m_1+m_1')}{\|\vect_1\|}x_1\right]} \right\} \cdot
    \\
     \left\{ \cos{\left[\frac{2\pi (m_2-m_2')}{\|\vect_2\|}x_2\right]}+\cos{\left[\frac{2\pi (m_2+m_2')}{\|\vect_2\|}x_2\right]} \right\} \cdot \\
    \left\{ \cos{\left[\frac{2\pi (m_3-m_3')}{\|\vect_3\|}x_3\right]}+\cos{\left[\frac{2\pi (m_3+m_3')}{\|\vect_3\|}x_3\right]} \right\} dx_1dx_2dx_3
    \\
    = \frac{1}{8}a_m\delta_{m_1m_1'}\delta_{m_2m_2'}\delta_{m_3m_3'}\|\vect_1\|\cdot\|\vect_2\|\cdot\|\vect_3\|
\end{gather*}
leading to
\begin{gather*}
    a_m = \frac{8}{V}\int_{-\frac{\|\vect_1\|}{2}}^{\frac{\|\vect_1\|}{2}}\int_{-\frac{\|\vect_2\|}{2}}^{\frac{\|\vect_2\|}{2}}\int_{-\frac{\|\vect_3\|}{2}}^{\frac{\|\vect_3\|}{2}}g(x_1,x_2,x_3)\cdot \\
    \cos{\left(\frac{2\pi m_1}{\|\vect_1\|}x_1\right)}\cos{\left(\frac{2\pi m_2}{\|\vect_2\|}x_2\right)}\cos{\left(\frac{2\pi m_3}{\|\vect_3\|}x_3\right)}dx_1dx_2dx_3
\end{gather*}
The multipliers $d_m,\beta_m,\gamma_m$ are calculated in the same fashion. However, for $b_m$ we find the following
\begin{align*}
    b_m &= \int_{-\frac{\|\vect_1\|}{2}}^{\frac{\|\vect_1\|}{2}}\int_{-\frac{\|\vect_2\|}{2}}^{\frac{\|\vect_2\|}{2}}\int_{-\frac{\|\vect_3\|}{2}}^{\frac{\|\vect_3\|}{2}} g(x_1,x_2,x_3)\cdot \sin{\left(\frac{2\pi m_1}{\|\vect_1\|}x_1\right)}\cos{\left(\frac{2\pi m_2}{\|\vect_2\|}x_2\right)}\cos{\left(\frac{2\pi m_3}{\|\vect_3\|}x_3\right)}dx_1dx_2dx_3 \\
    &= \int_{-\frac{\|\vect_1\|}{2}}^{0}\int_{-\frac{\|\vect_2\|}{2}}^{0}\int_{-\frac{\|\vect_3\|}{2}}^{0} g(x_1,x_2,x_3)\cdot \sin{\left(\frac{2\pi m_1}{\|\vect_1\|}x_1\right)}\cos{\left(\frac{2\pi m_2}{\|\vect_2\|}x_2\right)}\cos{\left(\frac{2\pi m_3}{\|\vect_3\|}x_3\right)}dx_1dx_2dx_3 + \\
    & \quad \int_{0}^{\frac{\|\vect_1\|}{2}}\int_{0}^{\frac{\|\vect_2\|}{2}}\int_{0}^{\frac{\|\vect_3\|}{2}} g(x_1,x_2,x_3)\cdot \sin{\left(\frac{2\pi m_1}{\|\vect_1\|}x_1\right)}\cos{\left(\frac{2\pi m_2}{\|\vect_2\|}x_2\right)}\cos{\left(\frac{2\pi m_3}{\|\vect_3\|}x_3\right)}dx_1dx_2dx_3
\end{align*}
with a change of variables $-u=x_1,\ -v=x_2, \ -w=x_3$ in the first triple integral
\begin{align*}
    b_m &= \int_{0}^{\frac{\|\vect_1\|}{2}}\int_{0}^{\frac{\|\vect_2\|}{2}}\int_{0}^{\frac{\|\vect_3\|}{2}} g(-u,-v,-w)\cdot \sin{\left(-\frac{2\pi m_1}{\|\vect_1\|}u\right)}\cos{\left(-\frac{2\pi m_2}{\|\vect_2\|}v\right)}\cos{\left(-\frac{2\pi m_3}{\|\vect_3\|}w\right)}dudvdw + \\
    & \quad \int_{0}^{\frac{\|\vect_1\|}{2}}\int_{0}^{\frac{\|\vect_2\|}{2}}\int_{0}^{\frac{\|\vect_3\|}{2}} g(x_1,x_2,x_3)\cdot \sin{\left(\frac{2\pi m_1}{\|\vect_1\|}x_1\right)}\cos{\left(\frac{2\pi m_2}{\|\vect_2\|}x_2\right)}\cos{\left(\frac{2\pi m_3}{\|\vect_3\|}x_3\right)}dx_1dx_2dx_3
\end{align*}
but since $\ionp = -u\frac{\vect_1}{\|\vect_1\|}-v\frac{\vect_2}{\|\vect_2\|}-w\frac{\vect_3}{\|\vect_3\|}$ represents the separation vector between two ions and $g$ depends only on their distance, meaning it only depends on $\|r\|$, we get the same value for $-\ionp = u\frac{\vect_1}{\|\vect_1\|}+v\frac{\vect_2}{\|\vect_2\|}+w\frac{\vect_3}{\|\vect_3\|}$ and
\begin{align*}
    b_m &= \int_{0}^{\frac{\|\vect_1\|}{2}}\int_{0}^{\frac{\|\vect_2\|}{2}}\int_{0}^{\frac{\|\vect_3\|}{2}} g(u,v,w)\cdot  \left[-\sin{\left(\frac{2\pi m_1}{\|\vect_1\|}u\right)}\right]\cos{\left(-\frac{2\pi m_2}{\|\vect_2\|}v\right)}\cos{\left(-\frac{2\pi m_3}{\|\vect_3\|}w\right)}dudvdw + \\
    & \quad \int_{0}^{\frac{\|\vect_1\|}{2}}\int_{0}^{\frac{\|\vect_2\|}{2}}\int_{0}^{\frac{\|\vect_3\|}{2}} g(x_1,x_2,x_3)\cdot \sin{\left(\frac{2\pi m_1}{\|\vect_1\|}x_1\right)}\cos{\left(\frac{2\pi m_2}{\|\vect_2\|}x_2\right)}\cos{\left(\frac{2\pi m_3}{\|\vect_3\|}x_3\right)}dx_1dx_2dx_3 \\
    &= -\int_{0}^{\frac{\|\vect_1\|}{2}}\int_{0}^{\frac{\|\vect_2\|}{2}}\int_{0}^{\frac{\|\vect_3\|}{2}} g(u,v,w)\cdot \sin{\left(\frac{2\pi m_1}{\|\vect_1\|}u\right)}\cos{\left(-\frac{2\pi m_2}{\|\vect_2\|}v\right)}\cos{\left(-\frac{2\pi m_3}{\|\vect_3\|}w\right)}dudvdw + \\
    & \quad \int_{0}^{\frac{\|\vect_1\|}{2}}\int_{0}^{\frac{\|\vect_2\|}{2}}\int_{0}^{\frac{\|\vect_3\|}{2}} g(x_1,x_2,x_3)\cdot \sin{\left(\frac{2\pi m_1}{\|\vect_1\|}x_1\right)}\cos{\left(\frac{2\pi m_2}{\|\vect_2\|}x_2\right)}\cos{\left(\frac{2\pi m_3}{\|\vect_3\|}x_3\right)}dx_1dx_2dx_3 \\
    &= 0
\end{align*}
Ultimately, we get that $ib_m=ic_m=i\alpha_m=i\delta_m=0$ and that Equation~(\ref{eq:etaSeries}) becomes
\begin{gather*}
    g(x_1,x_2,x_3) = \sum_{m_1=0}^{\infty}\sum_{m_2=0}^{\infty}\sum_{m_3=0}^{\infty} [a_m\cos{(u_1x_1)}\cos{(u_2x_2)}\cos{(u_3x_3)} + \\
    d_m\sin{(u_1x_1)}\sin{(u_2x_2)}\cos{(u_3x_3)} + \beta_m\sin{(u_1x_1)}\cos{(u_2x_2)}\sin{(u_3x_3)} + \\
    \gamma_m\cos{(u_1x_1)}\sin{(u_2x_2)}\sin{(u_3x_3)}]
\end{gather*}
which can be transformed into a series of cosine terms using trigonometric identities 
\begin{equation}
\begin{gathered}
    \label{eq:cosineSeriesFinal}
    g(x_1,x_2,x_3) = \sum_{m_1=0}^{\infty}\sum_{m_2=0}^{\infty}\sum_{m_3=0}^{\infty} \{ \\
    a_m[\cos{(u_1x_1+u_2x_2+u_3x_3)} + \cos{(u_1x_1+u_2x_2-u_3x_3)} + \\ \cos{(u_1x_1-u_2x_2+u_3x_3)} + \cos{(u_1x_1-u_2x_2-u_3x_3)}] + \\
    d_m[-\cos{(u_1x_1+u_2x_2+u_3x_3)} - \cos{(u_1x_1+u_2x_2-u_3x_3)} + \\ \cos{(u_1x_1-u_2x_2+u_3x_3)} + \cos{(u_1x_1-u_2x_2-u_3x_3)}] + \\ \beta_m[-\cos{(u_1x_1+u_2x_2+u_3x_3)} + \cos{(u_1x_1+u_2x_2-u_3x_3)} \\ -\cos{(u_1x_1-u_2x_2+u_3x_3)} + \cos{(u_1x_1-u_2x_2-u_3x_3)}] + \\
    \gamma_m[-\cos{(u_1x_1+u_2x_2+u_3x_3)} + \cos{(u_1x_1+u_2x_2-u_3x_3)} + \\ \cos{(u_1x_1-u_2x_2+u_3x_3)} - \cos{(u_1x_1-u_2x_2-u_3x_3)}] \}
\end{gathered}
\end{equation}
and, finally, because of Equations~(\ref{eq:expSeries}),(\ref{eq:expSeriesMult}),(\ref{eq:multsRelation}),(\ref{eq:multVal}),(\ref{eq:cosineSeriesFinal})
\begin{equation}
    f(\ionp) = \frac{4\pi}{V\|\rlattice\|^2}\sum_{m}\exp{\bigg(-\frac{\|\rlattice\|^2}{4\alpha^2}\bigg)}\cos{(\rlattice\ionp)}.
\end{equation}
\end{proof}

\subsection{Proof of Proposition~\ref{prop:ECoul} (See page~\pageref{prop:ECoul})}
\begin{proof}
    The proof of this Theorem follows from Propositions~\ref{prop:real}, ~\ref{prop:reciprocal} and the addition of one last summation term. Since we have added pair interactions to the summation of Equation~(\ref{eq:potfield_i}) so that
    \begin{equation*}
        \energy_{Coul}^L(\allionp,L) = \frac{k_e}{2}\sum_n\sum_{i=1}^N\sum_{j=1}^N\phi_j^L(\ionp_i+\lattice)
    \end{equation*}
    a self term $\energy^{self}_{Coul}$ is subtracted from the final result and becomes the third summand of the energy potential in Ewald form. With the help of the formula in Equation~(\ref{eq:potfieldL}) and the limit when $r_i$ approaches $r_j$ $\lim_{z\to 0}\erf{z} = \frac{2}{\sqrt{\pi}}z$, we have the direct evaluation of $\energy^{self}_{Coul}$ in real space
\begin{equation}
    \energy^{self}_{Coul} = -k_e\frac{\alpha}{\sqrt{\pi}}\sum_{i=1}^N q_i^2
\end{equation}
\end{proof}

\section{Proofs for the InflatedCellTruncation algorithm} \label{app2}
\subsection{Proof of Theorem~\ref{th:translation} (See page~\pageref{th:translation})}
\begin{proof}
    Let $t$ be the vector that translates $P$ to $P'$. We want $P$ and $P'$ to be parallel, so their distance is defined by $t$ and, as a result, $t$ is perpendicular to both $P$ and $P'$, or, in other words, 
    \begin{equation}\label{eq:tranparallel}
        t\parallel \vec{N}_{P'}
    \end{equation}
    where $\vec{N}_{P'}$ is the normal vector of plane $P'$. We want $P'$ tangent to sphere $(O,r_{off})$ and let $p$ be the point of $P'$ that touches $(O,r_{off})$. The distance between the centre $O$ and $P'$ is given by the length of vector $\vec{Op}$, which is perpendicular to $P'$. Then
    \begin{gather}\label{eq:opparallel}
        \vec{Op}\parallel \vec{N}_{P'} \\
        \|\vec{Op}\| = r_{off}
    \end{gather}
    Ultimately, from \ref{eq:tranparallel}, \ref{eq:opparallel}
    \begin{equation}
        t\parallel\vec{Op}\parallel\upsilon
    \end{equation}
    where $\upsilon$ is the vector of height of $\cell_n$ that is perpendicular to $P$. It follows that, since $O$ is the centre of gravity of the unit cell,
    \begin{equation}
        \vec{Op} = t - \frac{\upsilon}{2} \Rightarrow \|\vec{Op}\|v = \|t\|v - \frac{\|\upsilon\|}{2}-v \Rightarrow r_{off} = \|t\| + \frac{\|\upsilon\|}{2}
    \end{equation}
    where $v$ is a unit vector.
\end{proof}

\subsection{Proof of Corollary~\ref{cor:imgs} (See page~\pageref{cor:imgs})}
\begin{proof}
    Let $\cell_{(0,0,0)}$ be the central unit cell with lattice vectors $\vect_1, \vect_2, \vect_3$ and $\cell_{a_1}, \cell_{a_2}, ..., \cell_{a_m}, \ a_i\in\mathbb{N}^3, \ i\in[m]$ images of the unit cell surrounding $\cell_{(0,0,0)}$. Let also $\epsilon_1$ be a line parallel to $\vect_1$ that goes through the centre of gravity $O$ of $\cell_{(0,0,0)}$. Line $\epsilon_1$ cuts the faces of $\cell_{(0,0,0)}, \cell_{a_1}, \cell_{a_2}, ..., \cell_{a_m}$ at points $p_0, p_1, ..., p_m$. Let $\epsilon_2$ be the line in the direction of $t$ as defined in Theorem~\ref{th:translation} that goes through $O$ of $\cell_{(0,0,0)}$ and cuts adjacent images of the unit cell at points $p'_0, p'_1, ..., p'_k$. Since the faces of the images of the unit cells cut by $\epsilon_1, \epsilon_2$ are parallel, we have from Thales' theorem that $m=k$ and $Op_0 : p_0p_1 : p_1p_2 : ... : p_{m-1}p_m = Op'_0 : p'_0p'_1 : p'_1p'_2 : ... : p'_{m-1}p'_m$. However, since we deal with images of the unit cell and $\epsilon_1$ is parallel to $\vect_1$, we have that 
    \begin{gather}\label{eq:Thales}
        \frac{Op_0}{p_0p_1}=\frac{Op'_0}{p'_0p'_1} \\
        \text{ and } \nonumber \\
        \frac{p_{i-1}p_i}{p_ip_{i+1}}=1=\frac{p'_{i-1}p'_i}{p'_ip'_{i+1}}, \ \forall i\in{1,...,m-1}
    \end{gather}
    We know that, because $O$ is the center of gravity
    \begin{gather}\label{eq:height}
        \frac{Op_0}{p_0p_1} = \frac{Op'_0}{p'_0p'_1} \Leftrightarrow \frac{\|\vect_1\|/2}{\|\vect_1\|} = \frac{\|\upsilon\|/2}{p'_0p'_1} \Leftrightarrow p'_0p'_1 = \|\upsilon\|
    \end{gather}
    From \ref{eq:Thales} and \ref{eq:height} and Theorem~\ref{th:translation} we conclude that we can fit $(\|t\|+\frac{\|\upsilon\|}{2})/\|\upsilon\| = \frac{\|t\|}{\|\upsilon\|} + 1/2$ many images between $P$ and $P'$.
\end{proof}

\section{Proofs related to potential forces -- differentiation} \label{app3}
\subsection{Proof of Proposition~\ref{prop:Cforces}, the Coulomb first derivatives w.r.t. ion positions (See page~\pageref{prop:Cforces})}
\begin{proof}
    Let $\ionp_{i,j,n} = \ionp_i + \lattice - \ionp_j$ be the separation vector in 3-dimensional space between a pair of ions $i,j$. The distance $\|\ionp_{i,j,n}\|$ is a multivariate function $\mathbb{R}^3\rightarrow\mathbb{R}$. Let us define the two following functions that take a vector $r_t=(\ionp_{tx},\ionp_{ty},\ionp_{tz})$ as the free variable
    \begin{equation*}
        d_j(\ionp_t) = \|\ionp_t + \lattice - \ionp_j\|, \ d_j:\mathbb{R}^3\rightarrow\mathbb{R}.
    \end{equation*}
    If we take the gradient of $d$ with respect to the components of $\ionp_t$, we have
     \begin{align*}   
        \nabla d_j(\ionp_t) &= (\frac{\partial d_j}{\partial \ionp_{tx}}, \frac{\partial d_j}{\partial \ionp_{ty}}, \frac{\partial d_j}{\partial \ionp_{tz}})(\ionp_t) \\
        &= (\frac{\ionp_{tx}+n_1l_{1x}+n_2l_{2x}+n_3l_{3x}-\ionp_{jx}}{\|\ionp_t + \lattice - \ionp_j\|}\cdot\frac{\partial}{\partial\ionp_{tx}}(\ionp_{tx}\hat{x}+\ionp_{ty}\hat{y}+\ionp_{tz}\hat{z}), \\
        & \quad \quad \frac{\ionp_{ty}+n_1l_{1y}+n_2l_{2y}+n_3l_{3y}-\ionp_{jy}}{\|\ionp_t + \lattice - \ionp_j\|}\cdot\frac{\partial}{\partial\ionp_{ty}}(\ionp_{tx}\hat{x}+\ionp_{ty}\hat{y}+\ionp_{tz}\hat{z}), \\ 
        & \quad \quad \frac{\ionp_{tz}+n_1l_{1z}+n_2l_{2z}+n_3l_{3z}-\ionp_{jz}}{\|\ionp_t + \lattice - \ionp_j\|}\cdot\frac{\partial}{\partial\ionp_{tz}}(\ionp_{tx}\hat{x}+\ionp_{ty}\hat{y}+\ionp_{tz}\hat{z})) \\
        &= \frac{\ionp_{t,j,n}}{\|\ionp_{t,j,n}\|}
    \end{align*}
    in which $\hat{x},\hat{y},\hat{z}$ are the Cartesian unit vectors. Accordingly, we have
    \begin{equation*}
        d^i(\ionp_t) = \|\ionp_i + \lattice - \ionp_t\|, \ d^i:\mathbb{R}^3\rightarrow\mathbb{R}
    \end{equation*}
    with gradient
     \begin{align*}   
        \nabla d^i(\ionp_t) &= (\frac{\partial d^i}{\partial \ionp_{tx}}, \frac{\partial d^i}{\partial \ionp_{ty}}, \frac{\partial d^i}{\partial \ionp_{tz}})(\ionp_t) \\
        &= (\frac{\ionp_{ix}+n_1l_{1x}+n_2l_{2x}+n_3l_{3x}-\ionp_{tx}}{\|\ionp_i + \lattice - \ionp_t\|}\cdot\frac{\partial}{\partial\ionp_{tx}}(-\ionp_{tx}\hat{x}-\ionp_{ty}\hat{y}-\ionp_{tz}\hat{z}), \\
        & \quad \quad \frac{\ionp_{iy}+n_1l_{1y}+n_2l_{2y}+n_3l_{3y}-\ionp_{ty}}{\|\ionp_i+ \lattice - \ionp_j\|}\cdot\frac{\partial}{\partial\ionp_{ty}}(-\ionp_{tx}\hat{x}-\ionp_{ty}\hat{y}-\ionp_{tz}\hat{z}), \\ 
        & \quad \quad \frac{\ionp_{iz}+n_1l_{1z}+n_2l_{2z}+n_3l_{3z}-\ionp_{tz}}{\|\ionp_i + \lattice - \ionp_t\|}\cdot\frac{\partial}{\partial\ionp_{tz}}(-\ionp_{tx}\hat{x}-\ionp_{ty}\hat{y}-\ionp_{tz}\hat{z})) \\
        &= -\frac{\ionp_{i,t,n}}{\|\ionp_{i,j,n}\|}.
    \end{align*}

Let us define $f_S(x) = \frac{k_e}{2}\frac{\erfc{\alpha x}}{x}, \ f_S:\mathbb{R}\rightarrow \mathbb{R}$ and $d_j,d^i$ as before. Then we have:
      
    \begin{equation} \label{eq:fCoulS}
        f_S' = - \frac{k_e}{2}\left( \frac{2\alpha}{\sqrt{\pi}}\exp{(-\alpha^2x^2)+\frac{\erfc{\alpha x}}{x}}\right)\frac{1}{x}
    \end{equation}
    
    so that
    \begin{equation} 
    \energy^S_{Coul}(\allionp, L) = \sum_{i,j}^{N'}\sum_n  q_iq_jf_{S}(d^j(\ionp_i))
    \end{equation}
    
    and then the forces are given by
    \begin{equation}
    \begin{aligned}\label{eq:CgradS}
    \mathcal{F}^S_{Coul} &= -\nabla_\ionp\energy^S_{Coul}(\allionp, L) = -(\frac{\partial \energy^S_{Coul}}{\partial \ionp_1}, \frac{\partial \energy^S_{Coul}}{\partial \ionp_2}, ..., \frac{\partial \energy^S_{Coul}}{\partial \ionp_N}), \\
    \frac{\partial \energy^S_{Coul}}{\partial \ionp_t} &= \sum_{j=1}^{N'}\sum_n q_tq_jf_S'(d_j(\ionp_t))\nabla d_j(\ionp_t) + \\
    & \pushright{ \sum_{i=1}^{N'}\sum_n q_iq_tf_S'(d^i(\ionp_t))\nabla d^i(\ionp_t), \ t\in[N].}
    \end{aligned}
    \end{equation}

Similarly, let $f_L(x) = \sum_{m}\frac{2\pi k_e}{V\|\rlattice\|^2}\exp{(-\frac{\|\rlattice\|^2}{4\alpha^2})}\cos{(\rlattice x)}, \ f_L:\mathbb{R}^3\rightarrow \mathbb{R}$ and
    \begin{equation}
        \nabla_x f_L = -\sum_{m}\frac{2\pi k_e}{V\|\rlattice\|^2}\exp{\bigg(-\frac{\|\rlattice\|^2}{4\alpha^2}\bigg)}\rlattice\sin{(\rlattice x)}.
    \end{equation}

    Using $f_L$ we get
    \begin{equation*} 
    \energy^L_{Coul}(\allionp, L) = \sum_{i,j}^{N}  q_iq_jf_{L}(\ionp_{i,j})
    \end{equation*}

    and
    \begin{equation}
    \begin{aligned}\label{eq:gradL}
    \mathcal{F}^L_{Coul} &= -\nabla_\ionp\energy^L_{Coul}(\allionp, L) = -(\frac{\partial \energy^L_{Coul}}{\partial \ionp_1}, \frac{\partial \energy^L_{Coul}}{\partial \ionp_2}, ..., \frac{\partial \energy^L_{Coul}}{\partial \ionp_N}), \\
    \frac{\partial \energy^L_{Coul}}{\partial \ionp_t} &= \sum_{j=1}^{N} q_tq_j\nabla_\ionp f_L(\ionp_{t,j}) - \sum_{i=1}^N q_iq_t\nabla_\ionp f_L(\ionp_{i,t}), \ t\in[N] \\
    \end{aligned}
    \end{equation}

    where $\ionp_{t,j} = \ionp_t-\ionp_j$ and vice versa. Ultimately
    \begin{gather}
    -\nabla_\ionp\energy_{Coul} = \mathcal{F}_{Coul} = \mathcal{F}^S_{Coul} + \mathcal{F}^L_{Coul}.
    \end{gather}
\end{proof}

\subsection{Proof of Proposition~\ref{prop:Bforces}, the Buckingham first derivatives w.r.t. ion positions (See page~\pageref{prop:Bforces})}
\begin{proof}
Let 
\begin{equation*}
\begin{aligned}
    g_{S1}(x) &= \exp{(-\frac{x}{\rho})}, &\ g_{S2}:\mathbb{R}\rightarrow \mathbb{R} \\
    g_{S2}(x) &= - \frac{\exp{(-\alpha^2x^2)}}{x^6}\left( 1+\alpha^2x^2+\frac{\alpha^4x^4}{2} \right), &\ g_{S1}:\mathbb{R}\rightarrow \mathbb{R}
\end{aligned}
\end{equation*}
with
\begin{equation*}
    \energy^S_{Buck}(\allionp, L) = \frac{1}{2}\sum_{i,j}^{N'}\sum_n\bigg[ A_{ij}g_{S1}(d_j(\ionp_i))-C_{ij}g_{S2}(d_j(\ionp_i)) \bigg]
\end{equation*}
so that 
    \begin{equation}
    \begin{aligned}
        g_{S1}'(x) &= -\frac{1}{\rho}\exp{(-\frac{x}{\rho})}, \\
        g_{S2}'(x) &= \frac{\exp{(-\alpha^2x^2)}}{x^5}\left( \frac{6}{x^2}+6\alpha^2+\alpha^6x^4+3\alpha^4x^2 \right)
    \end{aligned}    
    \end{equation}
    then 
    \begin{equation}
    \begin{aligned}\label{eq:BgradS}
    \mathcal{F}^S_{Buck} &= -\nabla_\ionp\energy^S_{Buck}(\allionp, L) = -(\frac{\partial \energy^S_{Buck}}{\partial \ionp_1}, \frac{\partial \energy^S_{Buck}}{\partial \ionp_2}, ..., \frac{\partial \energy^S_{Buck}}{\partial \ionp_N}), \\
    \frac{\partial\energy^S_{Buck}}{\partial\ionp_t} &= \frac{1}{2}\sum_{j=1}^{N'}\sum_{n}  \bigg[ A_{tj}g_{S1}'(d_j(\ionp_t)) -C_{tj}g_{S2}'(d_j(\ionp_t)) \bigg]\nabla d_j(\ionp_t) + \\
    & \quad \frac{1}{2}\sum_{i=1}^{N'}\sum_{n}  \bigg[ A_{it}\nabla_\ionp g_{S1}(d^i(\ionp_t)) -C_{it}\nabla_\ionp g_{S2}(d^i(\ionp_t)) \bigg]\nabla d^i(\ionp_t))
    \end{aligned}
    \end{equation}
    with $t\in[N]$ and in the same fashion if 
    \begin{equation*}
    \begin{gathered}
        g_L(x) = -\frac{\pi^{3/2}}{12V}\sum_{m}\left[\sqrt{\pi}\cdot\erfc{\frac{\rlattice}{2\alpha}} + \left( \frac{4\alpha^3}{\|\rlattice\|^3}-\frac{2\alpha}{\rlattice}\right)\exp{\left(-\frac{\|\rlattice\|^2}{4\alpha^2}\right)}\right]\cdot\cos{(\rlattice x)}\|\rlattice\|^3
    \end{gathered}
    \end{equation*}     
    and  $g_L:\mathbb{R}^3\rightarrow\mathbb{R}$ with
    \begin{equation*}
        \energy^L_{Buck}(\allionp, L) = \frac{1}{2}\sum_{i,j}^{N}C_{ij}g_L(\ionp_{i,j})
    \end{equation*}
    then
     \begin{equation*}
    \begin{gathered}
        \nabla_x g_{L} = \frac{\pi^{3/2}}{12V}\sum_{m} \left[ \sqrt{\pi}\cdot\erfc{\frac{\|\rlattice\|}{2\alpha}} + \left( \frac{4\alpha^3}{\|\rlattice\|^3}-\frac{2\alpha}{\|\rlattice\|}\right) \exp{\bigg(-\frac{\|\rlattice\|^2}{4\alpha^2}\bigg)} \right]\cdot \rlattice\sin(\rlattice x).
    \end{gathered}
    \end{equation*}
    So, eventually 
    \begin{equation}
    \begin{aligned}
    \mathcal{F}^L_{Buck} &= -\nabla_\ionp\energy^L_{Buck}(\allionp, L) = -(\frac{\partial \energy^L_{Buck}}{\partial \ionp_1}, \frac{\partial \energy^L_{Buck}}{\partial \ionp_2}, ..., \frac{\partial \energy^L_{Buck}}{\partial \ionp_N}), \\
    \frac{\partial\energy^L_{Buck}}{\partial\ionp_t} &= \frac{1}{2}\sum_{j=1}C_{tj}\nabla_\ionp g_L(\ionp_{t,j}) - \frac{1}{2}\sum_{i=1}C_{it}\nabla_\ionp g_L(\ionp_{i,t})), \ t\in[N].
    \end{aligned}
    \end{equation}

    Thus
    \begin{gather}
    -\nabla_\ionp\energy_{Buck} = \mathcal{F}_{Buck} = \mathcal{F}^S_{Buck} + \mathcal{F}^L_{Buck}.
    \end{gather} 
\end{proof}

\subsection*{Proof of Proposition~\ref{prop:stresstheo}, stress and strain (See page~\pageref{prop:stresstheo})}
\begin{proof}
    Let us assume that at the start of the relaxation, the strain is zero. If we form a matrix $L_{0}\in \mathbb{R}^{3\times 3}$ using the lattice vectors $L$ in the initial state, the transformation of the vectors with strains would be
\begin{equation}
\label{eq:strainUpd}
    L_1 = (I+\epsilon_1)^TL_{0}, \ L,\epsilon_1\in\mathbb{R}^{3\times 3}
\end{equation}
where $I\in \mathbb{R}^{3\times 3}$ stands for the identity matrix. By continuously updating the lattice vectors' matrix as in Equation~(\ref{eq:strainUpd}), we can retrieve a sequence $L_0,L_{1}, L_{2}, .., L_{n}$ such that $L_{n}$ represents the lattice in equilibrium. This fractional change is called strain and involves two kinds of deformation; the normal strain $(\epsilon_{xx},\epsilon_{yy},\epsilon_{zz} \in \mathbb{R})$, that expresses changes to the length of the lattice vectors and shear strain $(\epsilon_{xy},\epsilon_{yz},\epsilon_{zx} \in \mathbb{R})$, which is the tangent of the angular change between two axes. We will now proceed to present the origins and physical interpretation of strain.

The deformation of the lattice represents motion of lattice points relative to each other. It is a linear transformation~\cite{Schlenker1978StraintensorParameters} of the lattice points and can be described by a displacement vector
\begin{equation}
    u = u_x\cdot\hat{x}+u_y\cdot\hat{y}+u_z\cdot\hat{z}, \ u\in\mathbb{R}^3, \ u_x,u_y,u_z\in\mathbb{R}
\end{equation}
where the components $u_x, u_y, u_z$ of vector $\mathfrak{u}$ are in reality continuous scalar functions of time $t\in \mathbb{R}^+$ and the position vector $p\in \mathbb{R}^3$ of a point, so that $u_x = u_x(p,t):\mathbb{R}^2\rightarrow\mathbb{R}$.

 By exploiting the partial derivatives of the displacement functions $u$ with respect to the position's $p$ each component, we can derive the fractional change per direction along the Cartesian coordinate axes as follows
\begin{equation}
\begin{gathered}\label{eq:strainEqs}
    \epsilon_{xx} = \frac{\partial u_x}{\partial p_x}, \quad \epsilon_{yy} = \frac{\partial u_y}{\partial p_y}, \quad \epsilon_{zz} = \frac{\partial u_z}{\partial p_z} \\
    \epsilon_{xy} = \frac{1}{2}\left(\frac{\partial u_y}{\partial p_x} + \frac{\partial u_x}{\partial p_y}\right) = \epsilon_{yx}, \quad \epsilon_{yz} = \frac{1}{2}\left(\frac{\partial u_z}{\partial p_y} + \frac{\partial u_y}{\partial p_z}\right) = \epsilon_{zy}, \\
    \epsilon_{xz} = \frac{1}{2}\left(\frac{\partial u_z}{\partial p_x} + \frac{\partial u_x}{\partial p_z}\right) = \epsilon_{zx}
\end{gathered}    
\end{equation}

Then, the second rank symmetrical tensor of strain at iteration $i$ is
\begin{equation}
    \epsilon_{i} = 
    \begin{bmatrix}
    \epsilon_{xxi} & \epsilon_{xyi} & \epsilon_{xzi} \\
    \epsilon_{xyi} & \epsilon_{yyi} & \epsilon_{yzi} \\
    \epsilon_{xzi} & \epsilon_{yzi} & \epsilon_{zzi} 
    \end{bmatrix}
\end{equation}

Since the strains are used as the external coordinates of the unit cell and are directly related to the lattice vectors, we can differentiate the energy function $\energy$ with respect to the strains. Under symmetrical, infinitesimal strain $\epsilon$, the derivative of the energy potential function with respect to strains expresses stress $\sigma$
    \begin{equation*}
        \sigma_{\lambda\mu} = \frac{1}{V}\frac{\partial \energy}{\partial \epsilon_{\lambda\mu}}, \ \forall \lambda,\mu \in \{x,y,z\}
    \end{equation*}
whereby very small, compared to the dimensions of the unit cell with volume $V$, changes to strain are measured. Stress is a symmetrical tensor comprising the components $\sigma_{\lambda\mu}\in \mathbb{R}, \ \forall \lambda,\mu \in \{x,y,z\}$ which are divided into normal and shear stresses as in the case o strains. Next, we will explain the tensor's physical interpretation that also accounts for its symmetry and properties.

Let us consider point $O\in \mathbb{R}^3$ in the interior of the crystal. When the lattice is deformed, forces start to act along the volume of the crystal per unit area, such that, internal forces are transmitted across a plane that passes through $O$ and separates the crystal into two halves. The cohesive forces crossing the plane with normal $n\in \mathbb{R}^3$ can be described by a force vector $f\in \mathbb{R}^3$ acting at point $O$, which belongs to a small area $A\in \mathbb{R}$ of the plane, as follows
\begin{equation}
    \overset{(n)}{T} = \lim_{\Delta A\to 0}\frac{f}{\Delta A}, \  \overset{(n)}{T}\in \mathbb{R}^3, \ f=(f_x,f_y,f_z)\in \mathbb{R}^3
\end{equation}
called stress vector. The stress vector is has arbitrary direction compared to the normal vector $n$. In order to be able to better investigate the phenomenon, we assume that the stress vector is acting on planes whose normal vectors are parallel to the Cartesian coordinate axes. After such an assumption, we can analyse the force vector $F$ into three components, one normal to each plane and two components parallel to the plane. For example, if we declare $\Delta A$ the piece of surface with normal vector $n$ oriented parallel to the x Cartesian coordinate axis and passing through $O$, then we can express $f$ as $f = f_x\cdot\hat{x} + f_y\cdot\hat{y} + f_z\cdot\hat{z}$, where $\hat{x},\hat{y},\hat{z}$ the unit vectors parallel to the Cartesian coordinate axes. The limits with respect to the confinement of $\Delta A$ express the average force exerted on the point $O$ of the face with area $\Delta A$ and normal vector parallel to the x axis. These are the components of the stress vector $\overset{(n)}{T}$ acting on the corresponding plane. It becomes apparent that each stress component is associated with two directions, the direction of $n$ and the direction of the component of $f$, and are respectively referred to by two subscripts. Hence, for a face with area $\Delta A$ and normal vector parallel to x axis the stress components are 
\begin{align*}
    \sigma_{xx} &= \lim_{\Delta A\to 0}\frac{f_x}{\Delta A} \\
    \sigma_{xy} &= \lim_{\Delta A\to 0}\frac{f_y}{\Delta A} \\
    \sigma_{xz} &= \lim_{\Delta A\to 0}\frac{f_z}{\Delta A}
\end{align*}
Consequently, by considering each plane whose normal vector direction is parallel to the Cartesian axes and passes through point $O$, we get three faces $\Delta A_x,\Delta A_y,\Delta A_z$ and nine stress components
\begin{align*}
    \sigma_{xx} = \lim_{\Delta A_x\to 0}\frac{f_x}{\Delta A_x} \quad \sigma_{xy} = \lim_{\Delta A_x\to 0}\frac{f_y}{\Delta A_x} \quad \sigma_{xz} = \lim_{\Delta A_x\to 0}\frac{f_z}{\Delta A_x} \\
    \sigma_{yx} = \lim_{\Delta A_y\to 0}\frac{f_x}{\Delta A_y} \quad \sigma_{yy} = \lim_{\Delta A_y\to 0}\frac{f_y}{\Delta A_y} \quad \sigma_{yz} = \lim_{\Delta A_y\to 0}\frac{f_z}{\Delta A_y}\\
    \sigma_{zx} = \lim_{\Delta A_z\to 0}\frac{f_x}{\Delta A_z} \quad \sigma_{zy} = \lim_{\Delta A_z\to 0}\frac{f_y}{\Delta A_z} \quad \sigma_{zz} = \lim_{\Delta A_z\to 0}\frac{f_z}{\Delta A_z}\\
\end{align*}
The stress components whose direction is parallel to the planes are called shear stresses, and the components normal to the planes are called normal stresses. We can express all of the stress components as in 
\begin{equation}
    \sigma_{\lambda\mu} = \lim_{\Delta A_\lambda\to 0}\frac{f_\mu}{\Delta A_\lambda}, \ \lambda,\mu \in \{x,y,z\}
\end{equation}
Furthermore, we impose equilibrium conditions such that the unit cell does not perform rigid-body movements and rotations. As explained in~\cite{Crandall1960AnSolids}, this results in the shear stresses to be equal $\sigma_{\lambda\mu}=\sigma_{\mu\lambda}, \ \lambda\neq \mu$. Similarly to the case of strains, the stresses are represented by a second rank symmetrical tensor 
\begin{equation*}
    \sigma = 
    \begin{bmatrix}
    \sigma_{xx} & \sigma_{xy} & \sigma_{xz} \\
    \sigma_{yx} & \sigma_{yy} & \sigma_{yz} \\
    \sigma_{zx} & \sigma_{zy} & \sigma_{zz} 
    \end{bmatrix}, \ \sigma_{\lambda\mu}\in\mathbb{R} \quad \forall \lambda,\mu \in \{x,y,z\}
\end{equation*}
which, because of the equal components, it comes to six independent components
\begin{equation}
    \label{eq:stress6}
    \sigma = \begin{pmatrix}
    \sigma_{xx} \\ \sigma_{yy} \\ \sigma_{zz} \\ \sigma_{yz} \\ \sigma_{xz} \\ \sigma_{xy}
    \end{pmatrix}, \ \sigma_{\lambda\mu}\in\mathbb{R} \quad \forall \lambda,\mu \in \{x,y,z\}
\end{equation}
which are used to update strain tensor $\epsilon_i$ as instructed in Equation~(\ref{eq:strainEqs}) and more elaborately shown in Equation~(\ref{eq:epsilonup})
\begin{equation} \label{eq:epsilonup}
    \epsilon_{i+1} =  
    \begin{bmatrix}
    1+(\epsilon_{xxi}+a\sigma_{xxi}) & \frac{1}{2}(\epsilon_{xyi}+a\sigma_{xyi}) & \frac{1}{2}(\epsilon_{xzi}+a\sigma_{xzi}) \\
    \frac{1}{2}(\epsilon_{xyi}+a\sigma_{xyi}) & 1+(\epsilon_{yyi}+a\sigma_{yyi}) & \frac{1}{2}(\epsilon_{yzi}+a\sigma_{yzi}) \\
    \frac{1}{2}(\epsilon_{xzi}+a\sigma_{xzi}) & \frac{1}{2}(\epsilon_{yzi}+a\sigma_{yzi}) & 1+(\epsilon_{yzi}+a\sigma_{yzi})
    \end{bmatrix}
\end{equation}
\end{proof}
where $a\in \mathbb{R}$. 

\subsection{Proof of Lemma~\ref{lem:stresscalc}, the derivatives w.r.t. lattice strain (See page~\pageref{prop:stresstheo})}
\begin{proof}
Because of the chain rule, we can obtain each stress component $\sigma_{\lambda\mu}$ as a sum of the derivatives of $\energy$ -- with respect to each parameter affected by deformation -- multiplied by the derivatives of the parameters with respect to strains. More precisely, we use the partial derivatives of $\energy$ next to the following equations
    \begin{equation}
    \begin{gathered}\label{eq:allprtdrvs}
        \frac{\partial\ionp_{t\psi}}{\partial\epsilon_{\lambda\mu}} = \delta_{\psi\lambda}\ionp_{t\mu}, \quad   
        \frac{\partial\vect_{t\psi}}{\partial\epsilon_{\lambda\mu}} = \delta_{\psi\lambda}\vect_{t\mu}, \quad 
        \frac{\partial\rvect_{t\psi}}{\partial\epsilon_{\lambda\mu}} = -\delta_{\psi\mu}\rvect_{t\lambda} \\
        \frac{\partial\lattice_\psi}{\partial\epsilon_{\lambda\mu}} = \delta_{\psi\lambda}\lattice_\mu, \quad
        \frac{\partial\rlattice_\psi}{\partial\epsilon_{\lambda\mu}} = -\delta_{\psi\mu}\rlattice_\lambda, \quad
        \frac{\partial V}{\partial\epsilon_{\lambda\mu}} = \delta_{\lambda\mu}V.
    \end{gathered}
\end{equation}
In the above latin letters as subscripts denote a part of the corresponding vector's name, e.g. for $\ionp_1$ we have $t=1$, and $\psi,\lambda,\mu$ refer to the Cartesian coordinates of the vector with subscript $t$, or to the components of the strain tensor $\epsilon$; $\delta$ is the Kronecker delta. We note that, while the component $\lambda$ of the real cell lattice vector $\vect_t(0)$ after distortion $\vect_t(\epsilon)$ becomes
\begin{equation*}
    \vect_{t\lambda}(\epsilon) = \sum^3_{\beta}(\delta_{\lambda\beta}+\epsilon_{\lambda\beta})\vect_{t\beta}(0)
\end{equation*}
the components of reciprocal vectors transform in the following way
\begin{equation*}
    \rvect_{t\lambda}(\epsilon) = \sum_{\beta}^3(\delta_{\beta\lambda}-\epsilon_{\beta\lambda})\rvect_{t\beta}(0).
\end{equation*}
\end{proof}

\subsection{Proof of Proposition~\ref{prop:Cstress}, the derivatives w.r.t. lattice strain (See page~\pageref{prop:Cstress})}
\begin{proof}
Using Proposition~\ref{prop:Cforces}, we can differentiate each of the Coulomb energy terms $\energy_{Coul}^S$, $\energy_{Coul}^L$, $\energy_{Coul}^{self}$ with respect to ion positions $\ionp_i, \ i \in [N]$. The stress tensor due to Coulomb forces is calculated with the help of Lemma~\ref{lem:stresscalc}, in which case we will also need the derivatives of the energy with respect to the lattice vectors $\lattice$,$\rlattice$ and the unit cell volume $V$. Then, these can be combined with the partial derivatives of Equation~(\ref{eq:allprtdrvs}) as shown in the next. From Equation~(\ref{eq:CgradS}) we get
\begin{gather}
     \frac{\partial \energy^S_{Coul}}{\partial \ionp_t} \frac{\partial\ionp_t}{\partial\epsilon_{\lambda\mu}} = \nonumber \\
     \sum_{j=1}^{N'}\sum_n q_tq_jf_S'(d_j(\ionp_t))\nabla d_j(\ionp_t)  \frac{\partial\ionp_t}{\partial\epsilon_{\lambda\mu}} + \sum_{i=1}^{N'}\sum_n q_iq_tf_S'(d^i(\ionp_t))\nabla d^i(\ionp_t)  \frac{\partial\ionp_t}{\partial\epsilon_{\lambda\mu}}  \nonumber\\
     = \sum_{j=1}^{N'}\sum_n q_tq_jf_S'(d_j(\ionp_t)) \left[\nabla d_j(\ionp_t) \right]_\lambda  \ionp_{t\mu} + \sum_{i=1}^{N'}\sum_n q_iq_tf_S'(d^i(\ionp_t)) \left[ \nabla d^i(\ionp_t) \right]_\lambda \ionp_{t\mu} \nonumber\\
     = \sum_{j=1}^{N'}\sum_n q_tq_jf_S'(d_j(\ionp_t)) \frac{[\ionp_{t,j,n}]_\lambda}{\|\ionp_{t,j,n}\|}  \ionp_{t\mu} - \sum_{i=1}^{N'}\sum_n q_iq_tf_S'(d^i(\ionp_t)) \frac{[\ionp_{i,t,n}]_\lambda}{\|\ionp_{i,t,n}\|} \ionp_{t\mu},  \nonumber\\
     t\in [N] \label{eq:CSstression}
\end{gather}
and
\begin{equation}
    \begin{gathered}
     \frac{\partial \energy^S_{Coul}}{\partial \lattice} \frac{\partial\lattice}{\partial\epsilon_{\lambda\mu}} = \sum_{i,j}^{N'}\sum_n q_iq_jf_S'(d_j(\ionp_i))\nabla d_j(\ionp_i)  \frac{\partial\lattice}{\partial\epsilon_{\lambda\mu}} \\
     = \sum_{i,j}^{N'}\sum_n q_iq_jf_S'(d_j(\ionp_i)) \left[\nabla d_j(\ionp_i)\right]_\lambda  \lattice_\mu  \\
     = \sum_{i,j}^{N'}\sum_n q_iq_jf_S'(d_j(\ionp_i)) \frac{[\ionp_{i,j,n}]_\lambda}{\|\ionp_{i,j,n}\|}  \lattice_\mu.
\end{gathered}
\end{equation}
which together give
\begin{equation}
\begin{gathered} \label{eq:CSstresslation}
    \sum_t^N \frac{\partial \energy^S_{Coul}}{\partial \ionp_t} \frac{\partial\ionp_t}{\partial\epsilon_{\lambda\mu}} + \frac{\partial \energy^S_{Coul}}{\partial \lattice} \frac{\partial\lattice}{\partial\epsilon_{\lambda\mu}} = \\ \sum_{i,j}^{N'}\sum_n q_iq_j f_S'(d_j(\ionp_i)) \frac{[\ionp_{i,j,n}]_\lambda}{\|\ionp_{i,j,n}\|}\cdot(\ionp_{i\mu}-\ionp_{j\mu}) + \sum_{i,j}^{N'}\sum_n q_iq_j f_S'(d_j(\ionp_i)) \frac{[\ionp_{i,j,n}]_\lambda}{\|\ionp_{i,j,n}\|}\cdot\lattice_\mu \\
    = \sum_{i,j}^{N'}\sum_n q_iq_j f_S'(d_j(\ionp_i)) \frac{[\ionp_{i,j,n}]_\lambda}{\|\ionp_{i,j,n}\|}\cdot[\ionp_{i,j,n}]_\mu
\end{gathered}
\end{equation}
As for the volume $V$ we need to use Equations~(\ref{eq:alphadrv}),~(\ref{eq:allprtdrvs}) like in the following
\begin{equation}\label{eq:CSstressV}
    \begin{gathered}
     \frac{\partial \energy^S_{Coul}}{\partial V}\frac{\partial V}{\partial\epsilon_{\lambda\mu}} = \left( \frac{k_e}{2}\sum_{i,j}^{N'}\sum_n q_iq_j\frac{1}{\|\ionp_{i,j,n}\|}\frac{\partial \erfc{(\alpha\|\ionp_{i,j,n}\|)}} {\partial V}\right) \frac{\partial V}{\partial\epsilon_{\lambda\mu}} \\
     = -\left( \frac{k_e}{\sqrt{\pi}}\sum_{i,j}^{N'}\sum_n q_iq_j\exp{(-\alpha^2\|\ionp_{i,j,n}\|^2)}\alpha'\right) \delta_{\lambda\mu}V
\end{gathered}
\end{equation}
and eventually, from Equations~(\ref{eq:CSstresslation}), ~(\ref{eq:CSstressV}) we get
\begin{equation}
\begin{gathered}\label{eq:CSstress}
    \frac{\partial \energy^S_{Coul}}{\partial\epsilon_{\lambda\mu}} =  \sum_t^{N} \frac{\partial \energy^S_{Coul}}{\partial \ionp_t} \frac{\partial\ionp_t}{\partial\epsilon_{\lambda\mu}} + \sum_n\frac{\partial \energy^S_{Coul}}{\partial \lattice} \frac{\partial\lattice}{\partial\epsilon_{\lambda\mu}} + \frac{\partial \energy^S_{Coul}}{\partial V}\frac{\partial V}{\partial\epsilon_{\lambda\mu}} \\
     = \sum_{i,j}^{N'}\sum_n q_iq_jk_e\bigg[ \exp{(-\alpha^2\|\ionp_{i,j,n}\|^2)} \alpha'\delta_{\lambda\mu}V - \Big. \\ \left.
     \frac{1}{2}\left( \frac{2\alpha}{\sqrt{\pi}}\frac{\exp{(-\alpha^2\|\ionp_{i,j,n}\|)}}{\|\ionp_{i,j,n}\|} + \frac{\erfc{\alpha\|\ionp_{i,j,n}\|}}{\|\ionp_{i,j,n}\|} \right)\frac{[\ionp_{i,j,n}]_\lambda[\ionp_{i,j,n}]_\mu}{\|\ionp_{i,j,n}\|}  \right].
\end{gathered}
\end{equation}
In the same fashion
\begin{equation}
\begin{gathered}\label{eq:CLstression}
     \frac{\partial \energy^L_{Coul}}{\partial \ionp_t} \frac{\partial\ionp_t}{\partial\epsilon_{\lambda\mu}} = \\
     \sum_{j=1}^{N} q_tq_j\nabla_\ionp f_L(\ionp_{t,j})\frac{\partial\ionp_t}{\partial\epsilon_{\lambda\mu}} - \sum_{i=1}^N q_iq_t\nabla_\ionp f_L(\ionp_{i,t})\frac{\partial\ionp_t}{\partial\epsilon_{\lambda\mu}}  \\
     = \sum_{j=1}^{N} q_tq_j [\nabla_\ionp f_L(\ionp_{t,j})]_\lambda\ionp_{t\mu} - \sum_{i=1}^N q_iq_t [\nabla_\ionp f_L(\ionp_{i,t})]_\lambda\ionp_{t\mu}  \\
     = - \frac{2\pi k_e}{V}\sum_{j=1}^{N} q_tq_j\sum_{m} \frac{\exp{(-\frac{\|\rlattice\|^2}{4\alpha^2})}}{\|\rlattice\|^2} \sin{(\rlattice \ionp_{t,j})}\rlattice_\lambda\ionp_{t\mu} +  \\
     \pushright{ \frac{2\pi k_e}{V}\sum_{i=1}^{N} q_iq_t\sum_{m}\frac{\exp{(-\frac{\|\rlattice\|^2}{4\alpha^2})}}{\|\rlattice\|^2} \sin{(\rlattice \ionp_{i,t})}\rlattice_\lambda\ionp_{t\mu} }  \\
     t\in [N] 
\end{gathered}
\end{equation}
and also
\begin{equation}
\begin{gathered}\label{eq:CLstresslat}
    \frac{\partial \energy^L_{Coul}}{\partial \ionp_t} \frac{\partial\rlattice}{\partial\epsilon_{\lambda\mu}} = \\
    \frac{2\pi k_e}{V}\sum_{i,j}^N\sum_n q_iq_j\frac{\exp{(-\frac{\rlattice^2}{4\alpha^2})}}{\|\rlattice\|^2} \bigg[ \cos{[\rlattice(\ionp_i-\ionp_j)]} \bigg(-\frac{1}{2\alpha^2}-\frac{2}{\rlattice^2} \bigg)\rlattice_\mu - \\
    \sin{[\rlattice(\ionp_i-\ionp_j)]\cdot(\ionp_{i\mu}-\ionp_{j\mu})}  \bigg]\rlattice_\lambda.
\end{gathered}
\end{equation}
Ultimately, Equation~(\ref{eq:CLstression}) is cancelled out because of Equation~(\ref{eq:CLstresslat}) as in
\begin{equation}
\begin{gathered}\label{eq:CLstresslation}
      \sum_t^N \frac{\partial \energy^L_{Coul}}{\partial \ionp_t} \frac{\partial\ionp_t}{\partial\epsilon_{\lambda\mu}} + \frac{\partial \energy^L_{Coul}}{\partial \ionp_t} \frac{\partial\rlattice}{\partial\epsilon_{\lambda\mu}} = \\
      - \frac{2\pi k_e}{V}\sum_{i,j}^{N} q_iq_j\sum_{m} \frac{\exp{(-\frac{\|\rlattice\|^2}{4\alpha^2})}}{\|\rlattice\|^2} \bigg[ \sin{(\rlattice \ionp_{i,j})}\rlattice_\lambda\cdot(\ionp_{i\mu}-\ionp_{j\mu}) + \\ 
      \cos{[\rlattice(\ionp_i-\ionp_j)]} \bigg(-\frac{1}{2\alpha^2}-\frac{2}{\rlattice^2} \bigg)\rlattice_\mu\rlattice_\lambda - \sin{[\rlattice(\ionp_{i,j})]\rlattice_\lambda\cdot(\ionp_{i\mu}-\ionp_{j\mu})} \bigg] \\
      = - \frac{2\pi k_e}{V}\sum_{i,j}^{N} q_iq_j\sum_{m} \frac{\exp{(-\frac{\|\rlattice\|^2}{4\alpha^2})}}{\|\rlattice\|^2}\cos{[\rlattice(\ionp_i-\ionp_j)]}\cdot \bigg(-\frac{1}{2\alpha^2}-\frac{2}{\rlattice^2} \bigg)\rlattice_\mu\rlattice_\lambda
\end{gathered}
\end{equation}
and as for the volume $V$, we again invoke Equation~(\ref{eq:alphadrv}) to address the change of parameter $\alpha$
\begin{equation}
\begin{gathered}\label{eq:CLstressV}
    \frac{\partial \energy^L_{Coul}}{\partial V}\frac{\partial V}{\partial\epsilon_{\lambda\mu}} = \frac{\partial}{\partial V}\bigg( \frac{2\pi k_e}{V}\sum_{i,j}^N\sum_n q_iq_j\frac{\exp{(-\frac{\|\rlattice\|^2}{4\alpha^2})}}{\|\rlattice\|^2} \cos{(\rlattice\ionp_{i,j})} \bigg) \frac{\partial V}{\partial\epsilon_{\lambda\mu}} \\
    = \frac{\partial}{\partial V}\bigg( \frac{2\pi k_e}{V}\sum_{i,j}^N\sum_n q_iq_j\frac{\exp{(-\frac{\|\rlattice\|^2}{4\alpha^2})}}{\|\rlattice\|^2} \cos{(\rlattice\ionp_{i,j})} \bigg) \delta_{\lambda\mu}V \\
    = \bigg[ -\frac{2\pi k_e}{V}\sum_{i,j}^N\sum_n q_iq_j\frac{\exp{(-\frac{\|\rlattice\|^2}{4\alpha^2})}}{V\|\rlattice\|^2} \cos{(\rlattice\ionp_{i,j})}\cdot V \cdot\bigg( 1-\frac{\|\rlattice\|^2}{2\alpha^3}\alpha'V \bigg) \bigg]\delta_{\lambda\mu} \\
    = -\delta_{\lambda\mu}\frac{2\pi k_e}{V}\sum_{i,j}^N\sum_n q_iq_j\frac{\exp{(-\frac{\|\rlattice\|^2}{4\alpha^2})}}{\|\rlattice\|^2} \cos{(\rlattice\ionp_{i,j})}\bigg( 1-\frac{\|\rlattice\|^2}{2\alpha^3}\alpha'V\bigg).
\end{gathered}
\end{equation}
The combination of Equations~(\ref{eq:CLstresslation}) and~(\ref{eq:CLstressV}) gives
\begin{equation}
\begin{gathered}\label{eq:CLstress}
    \frac{\partial \energy^L_{Coul}}{\partial\epsilon_{\lambda\mu}} = \sum_t^N \frac{\partial \energy^L_{Coul}}{\partial \ionp_t} \frac{\partial\ionp_t}{\partial\epsilon_{\lambda\mu}} + \frac{\partial \energy^L_{Coul}}{\partial \ionp_t} \frac{\partial\rlattice}{\partial\epsilon_{\lambda\mu}} + \frac{\partial \energy^L_{Coul}}{\partial V} \\
    = \frac{2\pi k_e}{V}\sum_{i,j}^N\sum_n q_iq_j\frac{\exp{(-\frac{\|\rlattice\|^2}{4\alpha^2})}}{\|\rlattice\|^2} \cos{(\rlattice\ionp_{i,j})}\cdot \\
    \bigg[\bigg( \frac{1}{2\alpha^2}+\frac{2}{\|\rlattice\|^2} \bigg)\rlattice_\mu\rlattice_\lambda - \delta_{\lambda\mu}\bigg( 1 - \frac{\|\rlattice\|^2}{2\alpha^3}\alpha'V \bigg) \bigg].
\end{gathered}
\end{equation}
Finally, we need to differentiate the self term with respect to $V$ using Equation~(\ref{eq:alphadrv}) and get
\begin{equation}
\begin{gathered}\label{eq:Cselfstress}
     \frac{\partial \energy^{self}_{Coul}}{\partial\epsilon_{\lambda\mu}} = \frac{\partial \energy^{self}_{Coul}}{\partial V}\frac{\partial V}{\partial\epsilon_{\lambda\mu}} =  \frac{\partial \energy^{self}_{Coul}}{\partial V}\cdot \delta_{\lambda\mu}V = -\frac{\alpha'k_e}{\sqrt{\pi}}\sum_i^Nq_i^2\cdot \delta_{\lambda\mu}V
\end{gathered}
\end{equation}
The calculation of Equations~(\ref{eq:CSstress}),~(\ref{eq:CLstress}) and~(\ref{eq:Cselfstress}) reveals the value of lattice stress because of the Coulomb forces, as shown by Lemma~\ref{lem:stresscalc}, hence it gives us the stress component with
\begin{equation}
    \sigma_{(Coul)\lambda\mu} = \frac{\partial \energy_{Coul}}{\partial\epsilon_{\lambda\mu}} = \frac{\partial \energy^S_{Coul}}{\partial\epsilon_{\lambda\mu}} + \frac{\partial \energy^L_{Coul}}{\partial\epsilon_{\lambda\mu}} + \frac{\partial \energy^{self}_{Coul}}{\partial\epsilon_{\lambda\mu}}.
\end{equation}
\end{proof}

\subsection{Proof of Proposition~\ref{prop:Bstress}, the derivatives w.r.t. lattice strain (See page~\pageref{prop:Bstress})}
\begin{proof}
Using Proposition~\ref{prop:Bforces}, we can differentiate each of the Buckingham energy terms $\energy_{Buck}^S$, $\energy_{Buck}^L$, $\energy_{Buck}^{self}$ with respect to ion positions $\ionp_i, \ i \in [N]$. According to the  Lemma~\ref{lem:stresscalc}, the stress onto the lattice due to 
Buckingham forces include the partial derivatives of the energy with respect to the lattice vectors $\lattice$,$\rlattice$ and the unit cell volume $V$. Then, these can be combined with the partial derivatives of Equation~(\ref{eq:allprtdrvs}) as shown in the next. From Equation~(\ref{eq:BgradS}) we get
\begin{equation}
\begin{gathered}
    \frac{\partial\energy^S_{Buck}}{\partial\ionp_t}\frac{\partial\ionp_t}{\partial\epsilon_{\lambda\mu}} = \frac{1}{2}\sum_{j=1}^{N'}\sum_{n}  \left[ A_{tj}g_{S1}'(d_j(\ionp_t)) -C_{tj}g_{S2}'(d_j(\ionp_t)) \right] \nabla d_j(\ionp_t)\frac{\partial\ionp_t}{\partial\epsilon_{\lambda\mu}} + \\
    \frac{1}{2}\sum_{i=1}^{N'}\sum_{n}  \left[ A_{it}\nabla_\ionp g_{S1}(d^i(\ionp_t)) -C_{it}\nabla_\ionp g_{S2}(d^i(\ionp_t)) \right]\nabla d^i(\ionp_t))\frac{\partial\ionp_t}{\partial\epsilon_{\lambda\mu}} \\
    =  \frac{1}{2}\sum_{j=1}^{N'}\sum_{n}  \left[ A_{tj}g_{S1}'(d_j(\ionp_t)) -C_{tj}g_{S2}'(d_j(\ionp_t)) \right] [\nabla d_j(\ionp_t)]_\lambda\ionp_{t\mu} + \\
    \frac{1}{2}\sum_{i=1}^{N'}\sum_{n}  \left[ A_{it}\nabla_\ionp g_{S1}(d^i(\ionp_t)) -C_{it}\nabla_\ionp g_{S2}(d^i(\ionp_t)) \right] [\nabla d^i(\ionp_t))]_\lambda\ionp_{t\mu} \\
    = \frac{1}{2}\sum_{j=1}^{N'}\sum_{n}  \left[ A_{tj}g_{S1}'(d_j(\ionp_t)) -C_{tj}g_{S2}'(d_j(\ionp_t)) \right] \frac{[\ionp_{t,j,n}]_\lambda}{\|\ionp_{t,j,n}\|}\ionp_{t\mu} + \\
    \frac{1}{2}\sum_{i=1}^{N'}\sum_{n}  \left[ A_{it}\nabla_\ionp g_{S1}(d^i(\ionp_t)) -C_{it}\nabla_\ionp g_{S2}(d^i(\ionp_t)) \right] \frac{[\ionp_{i,t,n}]_\lambda}{\|\ionp_{i,t,n}\|}\ionp_{t\mu} \\
    t\in [N]
\end{gathered}
\end{equation}
also
\begin{equation}
\begin{gathered}
    \frac{\partial\energy^S_{Buck}}{\partial\lattice}\frac{\partial\lattice}{\partial\epsilon_{\lambda\mu}} = \frac{1}{2}\sum_{i,j}^{N'}\sum_{n}  \left[ A_{ij}\nabla_\ionp g_{S1}(d_j(\ionp_i)) -C_{ij}\nabla_\ionp g_{S2}(d_j(\ionp_i)) \right] \frac{[\ionp_{i,j,n}]_\lambda}{\|\ionp_{i,j,n}\|}\lattice_{\mu}
\end{gathered}
\end{equation}
so, altogether
\begin{equation}
\begin{gathered}\label{eq:BSstresslation}
    \sum_t^N \frac{\partial\energy^S_{Buck}}{\partial\ionp_t}\frac{\partial\ionp_t}{\partial\epsilon_{\lambda\mu}} + \frac{\partial\energy^S_{Buck}}{\partial\lattice}\frac{\partial\lattice}{\partial\epsilon_{\lambda\mu}} = \\
     \frac{1}{2}\sum_{i,j}^{N'}\sum_{n}  \left[ A_{ij}\nabla_\ionp g_{S1}(d_j(\ionp_i)) -C_{ij}\nabla_\ionp g_{S2}(d_j(\ionp_i)) \right] \frac{[\ionp_{i,j,n}]_\lambda}{\|\ionp_{i,j,n}\|}\cdot(\ionp_{i\mu}-\ionp_{j\mu}) + \\
     \frac{1}{2}\sum_{i,j}^{N'}\sum_{n}  \left[ A_{ij}\nabla_\ionp g_{S1}(d_j(\ionp_i)) -C_{ij}\nabla_\ionp g_{S2}(d_j(\ionp_i)) \right] \frac{[\ionp_{i,j,n}]_\lambda}{\|\ionp_{i,j,n}\|}\lattice_{\mu} \\
     =  \frac{1}{2}\sum_{i,j}^{N'}\sum_{n}  \left[ A_{ij}\nabla_\ionp g_{S1}(d_j(\ionp_i)) -C_{ij}\nabla_\ionp g_{S2}(d_j(\ionp_i)) \right] \frac{[\ionp_{i,j,n}]_\lambda}{\|\ionp_{i,j,n}\|}[\ionp_{i,j,n}]_\mu.
\end{gathered}
\end{equation}
Then, with respect to $V$ we have
\begin{equation}
\begin{gathered}
    \frac{\partial\energy^S_{Buck}}{\partial V} = -\frac{1}{2}\sum_{i,j}^{N'}\sum_n \frac{C_{ij}}{\|\ionp_{i,j,n}\|^6}\exp{(-\alpha^2\|\ionp_{i,j,n}\|^2)}(-\alpha'\alpha^5\|\ionp_{i,j}\|^6)
\end{gathered}
\end{equation}
and, eventually
\begin{equation}
\begin{gathered}\label{eq:BSstress}
    \frac{\partial\energy^S_{Buck}}{\partial \epsilon_{\lambda\mu}} = -\frac{1}{2}\sum_{i,j}^{N'}\sum_n \frac{C_{ij}}{\|\ionp_{i,j,n}\|^6}\exp{(-\alpha^2\|\ionp_{i,j,n}\|^2)}\cdot \\
    \bigg[\bigg(\frac{6}{\|\ionp_{i,j,n}\|}+6\alpha^2+\alpha^6\|\ionp_{i,j,n}\|^4+3\alpha^4\|\ionp_{i,j,n}\|^2 \bigg)[\ionp_{i,j,n}]_\lambda[\ionp_{i,j,n}]_\mu + \\
    \alpha'\alpha^5\|\ionp_{i,j,n}\|^6V\delta_{\lambda\mu} \bigg]
\end{gathered}
\end{equation}
Moving to the long range interactions, we have
\begin{equation}
\begin{gathered}\label{eq:BLstression}
    \frac{\partial\energy^L_{Buck}}{\partial\ionp_t}\frac{\partial\ionp_t}{\partial\epsilon_{\lambda\mu}} = \frac{1}{2}\sum_{j=1}C_{tj}\nabla_\ionp g_L(\ionp_{t,j})\frac{\partial\ionp_t}{\partial\epsilon_{\lambda\mu}} - \frac{1}{2}\sum_{i=1}C_{it}\nabla_\ionp g_L(\ionp_{i,t}))\frac{\partial\ionp_t}{\partial\epsilon_{\lambda\mu}} \\
    = \frac{1}{2}\sum_{j=1}C_{tj}[\nabla_\ionp g_L(\ionp_{t,j})]_\lambda\ionp_{t\mu} - \frac{1}{2}\sum_{i=1}C_{it}[\nabla_\ionp g_L(\ionp_{i,t}))]_\lambda\ionp_{t\mu} \\
    = - \frac{1}{2}\sum_{j=1}C_{tj}\frac{\pi^{3/2}}{12V}\sum_{m} \left[ \sqrt{\pi}\cdot\erfc{\frac{\|\rlattice\|}{2\alpha}} + \right. \\ \left.
    \left( \frac{4\alpha^3}{\|\rlattice\|^3}-\frac{2\alpha}{\|\rlattice\|}\right) \exp{\bigg(-\frac{\|\rlattice\|^2}{4\alpha^2}\bigg)} \right]\cdot \rlattice_\lambda\sin(\rlattice \ionp_{t,j}) + \\ 
    \frac{1}{2}\sum_{j=1}C_{it}\frac{\pi^{3/2}}{12V}\sum_{m} \left[ \sqrt{\pi}\cdot\erfc{\frac{\|\rlattice\|}{2\alpha}} + \right. \\ \left.
    \left( \frac{4\alpha^3}{\|\rlattice\|^3}-\frac{2\alpha}{\|\rlattice\|}\right) \exp{\bigg(-\frac{\|\rlattice\|^2}{4\alpha^2}\bigg)} \right]\cdot \rlattice_\lambda\sin(\rlattice \ionp_{i,t}), \\
    t\in[N].
\end{gathered}
\end{equation}
and then
\begin{equation}
\begin{gathered}\label{eq:BLstresslat}
    \frac{\partial\energy^L_{Buck}}{\partial\rlattice}\frac{\partial\rlattice}{\partial\epsilon_{\lambda\mu}} = - \frac{1}{2}\sum_{i,j}^N C_{ij}\frac{\pi^{3/2}}{12V}\sum_{m} \cos{(\rlattice\ionp_{i,j})}\cdot \\
    \bigg[3\sqrt{\pi}\|\rlattice\|\erfc{\frac{\|\rlattice\|}{2\alpha}} - \|\rlattice\|^2 \frac{\exp{(-\frac{\|\rlattice\|^2}{4\alpha^2})}}{\alpha} -4\alpha\exp{\bigg(-\frac{\|\rlattice\|^2}{4\alpha^2}\bigg)} - \\
    (4\alpha^3-2\alpha\|\rlattice\|^2)\frac{\exp{(-\frac{\|\rlattice\|^2}{4\alpha^2})}}{2\alpha^2} \bigg]\rlattice_\mu\rlattice_\lambda.
\end{gathered}
\end{equation}
Once more, the combination of Equations~(\ref{eq:BLstression}) and~(\ref{eq:BLstresslat}) results into the removal of some terms because of mutual cancellation and we get
\begin{equation}
\begin{gathered}\label{eq:BLstressionlat}
    \sum_t^N\frac{\partial\energy^L_{Buck}}{\partial\ionp_t}\frac{\partial\ionp_t}{\partial\epsilon_{\lambda\mu}} + \frac{\partial\energy^L_{Buck}}{\partial\rlattice}\frac{\partial\rlattice}{\partial\epsilon_{\lambda\mu}} = \\
    - \frac{1}{2}\sum_{i,j}^N C_{ij}\frac{\pi^{3/2}}{12V}\sum_{m} \cos{(\rlattice\ionp_{i,j})}\bigg[ 3\sqrt{\pi}\|\rlattice\|\erfc{\frac{\|\rlattice\|}{2\alpha}} \\
    - 6\alpha\exp{\bigg(-\frac{\|\rlattice\|^2}{4\alpha^2}\bigg)} \bigg]\rlattice_\mu\rlattice_\lambda.
\end{gathered}
\end{equation}
then we also have
\begin{equation}
\begin{gathered}\label{eq:BLstressV}
    \frac{\partial \energy^L_{Buck}}{\partial V}\frac{\partial V}{\partial\epsilon_{\lambda\mu}} = \\
    - \frac{1}{2}\sum_{i,j}^N C_{ij}\frac{\pi^{3/2}}{12V}\sum_{m} \cos{(\rlattice\ionp_{i,j})}\delta_{\lambda\mu}\bigg[ -\sqrt{\pi}\erfc{\frac{\|\rlattice\|}{2\alpha}}\|\rlattice\|^3 + \\
    2\alpha(-2\alpha^2+ \|\rlattice\|^2+6\alpha V\alpha')\exp{\bigg( -\frac{\|\rlattice\|^2}{4\alpha^2} \bigg)} \bigg]
\end{gathered}
\end{equation}
Finally, by employing Equations~(\ref{eq:BLstressionlat}) and~(\ref{eq:BLstressV}),
\begin{equation}
\begin{gathered}\label{eq:BLstress}
    \frac{\partial \energy^L_{Buck}}{\partial\epsilon_{\lambda\mu}} = \frac{1}{2}\sum_{i,j}^N C_{ij}\frac{\pi^{3/2}}{12V}\sum_{m} \cos{(\rlattice\ionp_{i,j})}\bigg[\bigg( 3\sqrt{\pi}\|\rlattice\|\erfc{\frac{\|\rlattice\|}{2\alpha}} - \\
    6\alpha\exp{\bigg(-\frac{\|\rlattice\|^2}{4\alpha^2}\bigg)}\bigg)\rlattice_\lambda\rlattice_\mu - \delta_{\lambda\mu}\bigg( -\sqrt{\pi}\erfc{\frac{\|\rlattice\|}{2\alpha}}\|\rlattice\|^3 + \\
    (-4\alpha^3+ 2\alpha\|\rlattice\|^2+12\alpha^2 V\alpha')\exp{\bigg( -\frac{\|\rlattice\|^2}{4\alpha^2} \bigg)} \bigg) \bigg]
\end{gathered}
\end{equation}
The final differentiation still missing to complete all first order derivative calculations is the self term differentiation of the Buckingham potential. This is presented in Equation~(\ref{eq:Bstresself})
\begin{equation}\label{eq:Bstresself}
    \frac{\partial \energy^{self}_{Buck}}{\partial\epsilon_{\lambda\mu}} = \frac{\partial \energy^{self}_{Buck}}{\partial V}\frac{\partial V}{\partial\epsilon_{\lambda\mu}} = \bigg( - \frac{1}{2}\sum_{i,j}^N C_{ij}\frac{\pi^{3/2}}{3}\cdot\frac{3V\alpha'\alpha^2-\alpha^3}{V^2} + \frac{1}{2}\sum_i^NC_{ii}\alpha^5\alpha' \bigg)V\delta_{\lambda\mu}
\end{equation}
and then by using Equations~(\ref{eq:BSstress}),~(\ref{eq:BLstress}),~(\ref{eq:Bstresself}), Lemma~\ref{lem:stresscalc} can give us the stress component
\begin{equation}
\begin{gathered}
    \sigma_{(Buck)\lambda\mu} = \frac{\partial \energy_{Buck}}{\partial\epsilon_{\lambda\mu}} = \frac{\partial \energy^S_{Buck}}{\partial\epsilon_{\lambda\mu}} + \frac{\partial \energy^L_{Buck}}{\partial\epsilon_{\lambda\mu}} + \frac{\partial \energy^{self}_{Buck}}{\partial\epsilon_{\lambda\mu}}.
\end{gathered}
\end{equation}

\end{proof}

\section{Buckingham Catastrophe} \label{buckcat}
One of the limitations of the widely used Buckingham-Coulomb energy potential is called Buckingham catastrophe~\cite{Angyan2020CHAPTERUnits} and is a hard situation to recover from for minimization techniques. It refers to the Buckingham potential, whose form causes the Buckingham-Coulomb model to have neighbourhoods of deep wells that mathematically tend to negative infinity.

Intuitively, in terms of chemical components, a small distance between ions can cause the attraction $\|\ionp_{i,j,n}\|^{-6}$ term of Equation~(\ref{eq:EBuck}) to overpower the repulsion term and continuously push the ions together until they start to merge. In terms of the mathematical approach, the $\|\ionp_{i,j,n}\|^{-6}$ term diverges as $\|\ionp_{i,j,n}\|^{6} \rightarrow 0$ and a minimization algorithm sets out to reach the infimum of $\energy$ in a neighbourhood of the potential energy surface that stretches to $-\infty$. Obviously, this creates an infinite loop that can only be salvaged by creating an iteration deadline, or placing constraints to Problem~\ref{eqn:Problem}. For our experiments we examined the possibility of arriving to such a catastrophe using unconstrained optimization and limiting the running time by iteration number.
\end{document}